\documentclass[a4paper,12pt]{article}
\usepackage[LGR, T1]{fontenc}
\usepackage[english]{babel}
\usepackage[french=guillemets]{csquotes}
\RequirePackage[colorlinks=true]{hyperref}
\hypersetup{
    linkcolor = {blue},
    citecolor = {red},
    urlcolor = {blue},
    }
\usepackage{lmodern}
\usepackage{mathtools}
\usepackage{amssymb}
\usepackage{mathrsfs}
\usepackage{color}
\usepackage{soul}
\usepackage{xspace}
\usepackage[normalem]{ulem}
\usepackage{enumitem}
\usepackage{graphicx} % Utile pour redimensionner les choses, par exemple le triangle au dessus des "égal par définition"
\usepackage{stmaryrd} % Pour llbracket et rrbracket
\SetSymbolFont{stmry}{bold}{U}{stmry}{m}{n} % Evite le warning lié au conflit entre \boldsymbol (pour rendre les caractères grecs gras) et stmaryrd.
\usepackage{subfig}
\usepackage{multicol}
\usepackage{authblk}
\usepackage{listings}
\usepackage{tikz}
\usepackage[style=alphabetic,backref=true]{biblatex}
\usepackage[a4paper,inner=2cm,outer=2cm,top=2cm, bottom = 2cm]{geometry}
\newcounter{compteur}[section] % Numérotation des items du type Définition, Proposition, etc.
\renewcommand\thecompteur{\thesection.\arabic{compteur}}

\newenvironment{Definition}
    {\refstepcounter{compteur} \vspace{0.5cm} \noindent \textbf{Definition \thecompteur~:}  }{\normalfont\vspace{0.2cm}}

\newenvironment{Proposition}
    {\refstepcounter{compteur} \vspace{0.5cm} \noindent \textbf{Proposition \thecompteur~:} \itshape }{\normalfont \vspace{0.1cm}}

\newenvironment{Lemma}
    {\refstepcounter{compteur} \vspace{0.2cm} \noindent \textbf{Lemma \thecompteur~:} \itshape }{\normalfont}

\newenvironment{Theoreme}
    {\refstepcounter{compteur} \vspace{0.5cm} \noindent \textbf{Theorem \thecompteur~:} \itshape }{\normalfont}

\newenvironment{Remarque}{\refstepcounter{compteur} \vspace{0.2cm} \noindent \textit{Remark \thecompteur. }}{\vspace{0.1cm}}

\newenvironment{Demo}[1][]{\vspace{0.2cm} \noindent  \textbf{Proof~#1~: }}{\hfill $\square$\vspace{0.1cm}}

\newenvironment{Hyp}{\refstepcounter{compteur} \vspace{0.2cm} \noindent \textbf{Hypothesis \thecompteur. }}{\vspace{0.1cm}}

\usepackage{fancyhdr}
\newcommand{\R}{\mathbb{R}}
\newcommand{\Q}{\mathbb{Q}}

\newcommand{\N}{\mathbb{N}}

\newcommand{\E}{\mathbb{E}}
\renewcommand{\P}{\mathbb{P}}
\newcommand{\D}{\mathrm{D}}
\renewcommand{\d}{\mathrm{d}}
\newcommand{\norm}[2]{{\left\Vert #1 \right\Vert}_{#2}}

\newcommand{\defeq}{:=}
\newcommand{\Cc}{\mathcal{C}}
\newcommand{\Ff}{\mathcal{F}}
\newcommand{\Gg}{\mathcal{G}}
\newcommand{\Hh}{\mathcal{H}}
\newcommand{\eps}{\varepsilon}

\newcommand{\ms}{m_{\sigma}}
\newcommand{\mm}{m_{\mu}}
\newcommand{\mx}{m_{X}^{(p)}}
\newcommand{\mxx}{m_{X}^{(p,2)}}
\newcommand{\mdx}{m_{\partial X}^{(p,q)}}
\newcommand{\mdxx}{m_{\partial X}^{(p,q,2)}}
\newcommand{\mddx}{m_{\partial^2 X}^{(p,q)}}
\newcommand{\mddxx}{m_{\partial^2 X}^{(p,q,2)}}

\newcommand{\msx}{m_{X, \sigma}^{(p)}}
\newcommand{\msxx}{m_{X, \sigma}^{(p,2)}}
\newcommand{\msdx}{m_{\partial X, \sigma}^{(p,q)}}
\newcommand{\msdxx}{m_{\partial X, \sigma}^{(p,q,2)}}
\newcommand{\msddx}{m_{\partial^2 X, \sigma}^{(p,q)}}
\newcommand{\msddxx}{m_{\partial^2 X, \sigma}^{(p,q,2)}}

\newcommand{\Blue}[1]{{\color{blue}{#1}\color{black}\xspace}}

\begin{document}

\title{A Poisson--Alekseev--Gröbner formula through Malliavin calculus for Poisson random integrals}
\author[1]{Paul Maurer\footnote{paul.maurer@inria.fr}}
\author[2]{Jérémy Zurcher\footnote{ jzurcher3@gatech.edu}}
\affil[1]{Universit{\'e} C{\^o}te d'Azur, Inria, CNRS, France}
\affil[2]{School of Mathematics, Georgia Institute of Technology, Atlanta, GA, 30332, USA}
\maketitle

\begin{abstract}
    In this paper, we establish an Alekseev--Gröbner formula for stochastic differential equations (SDEs) driven by a Poisson random measure, which express the global error between a functional of two processes solution of SDEs started at the same initial condition, in terms of the infinitesimal error (i.e, the difference between the SDEs coefficients). In particular, we consider the situation where the flow process is only assumed to be jointly stochastically continuous with respect to space and time. Our proof relies on a new approach for the definition of the Skorohod--Poisson integral to treat the anticipating term appearing in the formula, based on a definition of the Malliavin derivative on a class of smooth random variables, instead of the more standard polynomial chaos approach.
\end{abstract}

\tableofcontents

\section{Introduction}

In this article, we establish an Itô--Alekseev--Gröbner (or Skorohod--Alekseev--Gröbner) formula in the case of SDEs driven by a (compensated) Poisson random measure. In the meantime, we introduce a new approach to Malliavin calculus and in particular Skorohod integrals in the Poisson setting, and prove an Itô formula for independent random functionals and general càdlàg semimartingales.

The Alekseev--Gröbner formula, as established in Alekseev \cite{alekseev1961} and Gröbner \cite{grobner1960} for ordinary differential equations, as a tool to express the global error in terms of the infinitesimal error between solutions of equations. It states that for a jointly continuous flow $X_{s,t}^{x} = x + \int_s^t b(r,X_{s,r}^x) \ \d r$ and a continuously differentiable function $Y$ satisfying $Y_t = \int_0^t A_r \ \d r$, we have that
\begin{equation} \label{eq:deterministicAG}
    X_{0,T}^{Y_0} - Y_T = \int_0^T \partial_x X_{r,T}^{Y_r} \ \big( b(r,Y_r) - A_r \big) \ \d r.
\end{equation}
This result, particularly useful to analyse perturbations of deterministic ODEs, has been generalized to the case of stochastic differential equations driven by a Brownian motion in the independent works from Del Moral and Singh \cite[Theorem 1.2]{delmoral2022} and from Hudde, Hurtznthaler, Jentzen and Mazzonetto \cite{ItoAG}.
The first work firstly demonstrates an backward Itô--Ventzell formula which is later used to show the stochastic Alekseev--Gröbner formula, while the second work relies on a time discretisation argument to demonstrate the stochastic Alekseev--Gröbner formula directly. 

The main difficulty in the extension of \eqref{eq:deterministicAG} to the stochastic setting is the presence of the anticipating term $\partial_x X_{r,T}^{Y_r}$, which may not treated by standard Itô integration. In both of the aforementioned works, this obstacle is tackled by the use of Malliavin calculus through the Skorohod integral, which is known to extend the Itô's one. To address the case of SDEs with jumps, we rely on the time discretisation argument. This method of proof requires to construct a version of the Skorohod integral that respects the Chasles relation, which is not possible with standard definitions in the Poisson setting (e.g. using the polynomial chaos approach). In addition, for the time discretisation to converges when the step size goes to zero, standard methods rely on the joint \textit{almost sure} continuity of the flow process (in time and in space). This argument does not hold for SDEs with jumps which are not almost surely continuous in time. Instead, we show that this convergence still holds with the weaker condition of joint \textit{stochastic} continuity.

\subsection{Main result and motivations}

 The main result of the paper is Theorem \ref{maintheo}. Informally speaking, considering a Poisson random measure $N$ on a filtered probability space $(\Omega, \Ff, (\Ff_t)_{t \in \R_+} , \P)$ with associated compensator $\nu$, Theorem \ref{maintheo} states that for an adapted flow process 
 \[ X_{s,t}^x = x + \int_s^t \mu \left( r,X_{s,r}^x \right) \ \d r + \int_s^t \int_{\R \backslash \{ 0 \}} \sigma \left( r,z,X_{s,{r^-}}^x \right) \ \tilde{N}(\d r, \d z), \] 
 an adapted stochastic process 
 \[ Y_t = Y_0 + \int_0^t A_r \ \d r + \int_0^t \int_{\R \backslash \{ 0 \}} B_{r,z} \ \tilde{N}(\d r, \d z), \] and a test function $f: \R \longrightarrow \R$, under suitable moment and regularity conditions, we have
\begin{multline} \label{eq:informal_itoAG}
      f \left( X_{0,T}^{Y_0} \right) - f(Y_T) = \int_{0}^{T} f' \left( X_{r,T}^{Y_{r}} \right) \partial_x X_{r,T}^{Y_{r}} \left(  \mu \left( r, Y_r \right) - A_r \right) \ \d r   \\
    + \int_{0}^{T} \int_{\R \backslash \{ 0 \}} \left\{ f \left( X_{r,T}^{Y_r + \sigma \left( r,Y_r,z \right)} \right) - f \left( X_{r,T}^{Y_{r}+B_{r,z}} \right) - f' \left( X_{r,T}^{Y_r} \right) \partial_x X_{r,T}^{Y_r} \Big(    \sigma \left( r,Y_r,z \right) - B_{r,z} \Big) \right\} \nu( \d z) \d r \\
     + \int_0^T \int_{\R \backslash \{ 0 \}} \left( f \left( X_{r,T}^{Y_{r} + \sigma \left( r,Y_{r},z \right)} \right) - f \Big( X_{r,T}^{Y_{r}+B_{r,z}} \Big) \right) \ \tilde{N}(\delta r, \delta z),
\end{multline}  
where the stochastic integral is a Skorohod one (denoted by the symbol $\delta$ instead of $\d$). Although it might not be clear at first glance that \eqref{eq:informal_itoAG} express the global error in term of the infinitesimal error, a Taylor formula with integral leftover applied to the function $x \mapsto f(X_{r,T}^{Y_r + x})$ in the second and third terms of the right hand side of \eqref{eq:informal_itoAG} makes appear the expected local error $\sigma \left( r,Y_r,z \right) - B_{r,z}$ in factor of these two terms leading to 
\begin{multline} \label{eq:informal_itoAG_Taylored}
    f \left( X_{0,T}^{Y_0} \right) - f(Y_T) = \int_{0}^{T} f' \left( X_{r,T}^{Y_{r}} \right) \partial_x X_{r,T}^{Y_{r}} \left(  \mu \left( r, Y_r \right) - A_r \right) \ \d r   \\
  + \int_{0}^{T} \int_{\R \backslash \{ 0 \}} F_{r,T,z} \left( \sigma \left( r,Y_r,z \right) - B_{r,z} \right) \nu( \d z) \d r \\
   + \int_0^T \int_{\R \backslash \{ 0 \}} G_{r,T,z} \left( \sigma \left( r,Y_r,z \right) - B_{r,z} \right) \ \tilde{N}(\delta r, \delta z).
\end{multline}
with
\begin{equation*}
  F_{r,T,z} \coloneqq \int_0^1 \left\{ f' \left( X_{r,T}^{Y_r + \lambda \sigma ( r,Y_r,z) + (1-\lambda) B_{r,z}} \right) \partial_x X_{r,T}^{Y_r + \lambda \sigma ( r,Y_r,z) + (1-\lambda) B_{r,z}} - f' \left( X_{r,T}^{Y_r} \right) \partial_x X_{r,T}^{Y_r} \right\} \d \lambda
\end{equation*}
and
\begin{equation*}
  G_{r,T,z} \coloneqq \int_0^1 f' \left( X_{r,T}^{Y_r + \lambda \sigma ( r,Y_r,z) + (1-\lambda) B_{r,z}} \right) \partial_x X_{r,T}^{Y_r + \lambda \sigma ( r,Y_r,z) + (1-\lambda) B_{r,z}} \d \lambda,
\end{equation*}
The formula \eqref{eq:informal_itoAG_Taylored} can be applied afterwards to derive probabilistic rates of convergence for approximations of SDEs driven by Poisson random measures. In particular,
\begin{itemize}
    \item We may obtain \textit{strong} rates of convergence from \eqref{eq:informal_itoAG} by taking $f(x) = x$ and then the $L^p$-norm on both sides of the equation. This requires, however, to have proper estimates for the $L^p$-norm of the Skorohod integral. In the diffusive case, this is possible using e.g \cite[Lemma 4]{alos1998} when $p=2$ or \cite[Lemma 5]{alos1998} when $p \in (2,4)$. In the Poisson case, an answer was recently given in \cite[Theorem 4.2]{Tara2025}.
    \item We may obtain \textit{weak} rates of convergence by directly taking the expectation in both sides of \eqref{eq:informal_itoAG}. In this manner, the Skorohod term vanishes and the analysis reduces to the decoupling and the estimation of the multiplicative stochastic terms, as well as the estimation of the infinitesimal errors themselves.
\end{itemize}

The precise hypotheses for equation \eqref{eq:informal_itoAG} to hold are given in Section \ref{sec:Hypotheses}. In particular, we need the Lévy measure $\nu$ and the flow process $X_{s,t}^x$ to have at least a finite fourth-order moment. While this requirement does not allow for an $\alpha$-stable driving Lévy process, it however includes any \textit{tempered} stable Lévy processes as driving SDE process. Tempered stable processes, introduced by Koponen \cite{koponen1995} in the context of turbulence in physics and later on formalized by Rosinski in \cite{rosinski2007}, are obtained by tempering the jump size of an $\alpha$-stable process at the level of the Lévy measure $\nu$. The tempering may be exponential ($\nu(\d z) = \frac{\d z}{|z|^{\alpha+1}} e^{-\beta z}$ for some parameter $\beta > 0$) or a truncation ($\nu(\d z) = \frac{\d z}{|z|^{\alpha+1}} \mathbf{1}_{ \{|z| \leq \beta \}}$). In any case, the resulting process has finite moments of all orders, and behaves like an $\alpha$-stable process for small times, while approximating a Brownian motion at large times (see \cite[Theorem 3.1]{rosinski2007}). The class of tempered stable processes provides stochastic models well-fited for turbulence (see \cite{novikov1994} and \cite{mantegna1994}, or also \cite[Section 6]{bossymaurer2024}), often more relevant than $\alpha$-stable processes themselves since many non-diffusive or non-Gaussian physical processes actually still have moments of any orders. Tempered stable processes are also widely used in the mathematical finance literature, through the name of CGMY processes, in reference to the authors Carr, Geman, Madan and Yor of the seminal article \cite{carr2002}. In fact, as discussed in \cite{geman2002}, they appear to provide more accurate models for option pricing than Brownian diffusions, as financial assets prices often exhibit large probability tails for the marginals and jumps in the dynamics (think about abrupt financial crisis).

An important motivation for our theorem is the possibility to obtain weak convergence rate for a tempered Euler--Maruyama scheme for SDEs which coefficients are only locally Lipschitz but satisfy a Lyapunov condition. In \cite{chak2024}, the Ito--Alekseev--Gröbner formula is effectively used to obtain a first order convergence of a stopped increment-tamed Euler--Maruyama approximation for an Itô SDE driven by a Brownian motion where the coefficients satisfy a locally Lipschitz assumption along with a Lyapunov condition. To our knowledge, there is no result in the litterature about the approximation of Poisson-driven SDEs with locally Lipschitz coefficients and infinite activity Lévy measure (which requires to approximate or neglect the jumps that are smaller than a certain treshold in addition to the time discretisation as in \cite{bossymaurer2024}).
To demonstrate such a result using the Poisson--Alekseev--Gröbner formula, one would need to go beyond the setting of \cite{BretonPrivault} and show the regularity and integrability of the flow in the case of locally Lipschitz coefficients. This will be the subject for a future work.

\subsection{Method of proof}

\paragraph{Treatment of anticipating integrands with extended filtrations. }
We start by giving a short heuristic as a motivation for the strategy of the proof. 
Consider a Lévy process $(L_t)_{t \in [0,3]}$. For $i \in \{0,1\}$, let's say that we would like to apply the "Itô formula" to the (random) function $f_i: x \ni \R \mapsto \sin(x  (L_{i+2} - L_{i+1}))$ and the process $L$ between times $i$ and $i+1$. By formally applying the Itô formula, the result would contain the stochastic integral 
\begin{equation} \label{eq:anticipating_levy}
    \int_i^{i+1} (L_{i+2} - L_{i+1}) \cos( L_{s^-} (L_{i+2} - L_{i+1})) \ \d L_s.
\end{equation}
Obviously, \eqref{eq:anticipating_levy} seems ill-defined due to the presence of the anticipating term $L_{i+2} - L_{i+1}$. Nevertheless, although \eqref{eq:anticipating_levy} is in fact ill-defined in the standard filtration $\mathbf{F}^i$ defined by 
\begin{equation*}
\mathbf{F}_t^i \defeq \sigma( L_s, \ i \leqslant s \leqslant t)
\end{equation*}
for $t \in [i,i+1]$, it is actually perfectly well-defined in the extended filtration $\mathbf{H}^i$ defined by 
\begin{equation*}
    \mathbf{H}^i_t = \mathbf{F}_t^i \vee \sigma(L_{i+2}-L_{i+1}) 
\end{equation*}
for $t \in [i,i+1]$. To point out the dependence of the extended filtration, we often employ the notation
\begin{equation*}
    \int_i^{i+1} (L_{i+2} - L_{i+1})\cos( L_{s^-}  (L_{i+2} - L_{i+1})) \ \d L_s^{\mathbf{H}^i}.
\end{equation*}
Moreover, the process $(L_t)_{t \in [i,i+1]}$ remains a semimartingale with respect to the filtration $\mathbf{H}^i$, due to the independence between the $\sigma$-algebras $\mathbf{F}_t^i$ and $\sigma(L_{i+2}-L_{i+1})$. This is the key observation that enables us to show the Itô formula for independent random functionals, which is Proposition \ref{prop:ito_formula_indep_cadlag} in the Appendix. 
However, we cannot combine the sum 
\begin{equation} \label{eq:combination_levy}
    \int_0^1 (L_{2} - L_{1}) \cos \big( L_{s^-} (L_2 - L_1) \big) \ \d L_s^{\mathbf{H}^1} + \int_1^2 (L_{3} - L_{2}) \cos \big( L_{s^-} (L_3 - L_2) \big) \ \d L_s^{\mathbf{H}^2}
\end{equation}
into a single Itô integral on the form 
\[ \int_0^2 (L_{\lceil s \rceil +1} - L_{\lceil s \rceil}) \cos \big( L_{s^-} (L_{\lceil s \rceil +1} - L_{\lceil s \rceil}) \big) \ \d L_s, \] 
because the process $((L_{\lceil s \rceil +1} - L_{\lceil s \rceil})\cos(L_{s^-} (L_{\lceil s \rceil +1} - L_{\lceil s \rceil})), s \in [0, 2])$ is not adapted to any filtration. Instead, we follow the idea of \cite{ItoAG} in the Brownian setting, and construct in the Poisson setting a Skorohod integral that respect the Chasles relation, allowing to combine the sum \eqref{eq:combination_levy} into a single Skorohod integral from $0$ to $2$.

\paragraph{Malliavin calculus for Poisson random integrals.}
In both the Brownian and Poisson settings, three ways emerge to define the Malliavin derivative: the smooth random variables approach (see \cite{Nualart} for the Brownian case), the polynomial chaos approach (see \cite{Last2016} for the definition in the Poisson case using the Itô integral), and the construction in a canonical space (see \cite{NualartVives}). In \cite{ItoAG}, in the Brownian case, the authors essentially rely on the first one, defining the Malliavin derivative on a class of simple functions with \textit{additional information}. However in the Poisson setting, this approach is often abandoned in favour of polynomial chaos in the literature, which is not suitable for the concept of additional information. Note however the work of Bally and Qin in \cite{bally2022} in this direction, with a different approach in the definition of the Skorohod integral. We introduce in Section \ref{sec:Malliavin} a new definition of the Malliavin derivative for smooth random variables based on Poisson random integrals, and prove that the standard results of Malliavin calculus hold with this definition. To do that, we make a slight change in the definition of smooth random variables: instead of using functions with growth at most polynomial, we define it with trigonometrical polynomials. A Fourier transform argument concludes that we constructed a dense set inside $L^2 (\Omega)$, for the $\sigma$-algebra generated by the underlying Poisson measure (see Lemma \ref{lem:DensiteSmoothRV}). This allows to define a Skorohod--Poisson stochastic integral (see Definition \ref{def:Skorohod_integral}) that respects the Chasles relation (see Theorem \ref{thmMalliavin:Chasles}). 

\paragraph{Sketch of the proof of Theorem \ref{maintheo}}
Now equipped with an efficient integration framework to deal with anticipating Poisson integrands, we provide a proof of Theorem \ref{maintheo} based on a time discretization, the use of an Itô formula with random, independent functional for general semimartingales, followed by a convergence in $L^1$ for the Lebesgue integrals and in the weak topology of $L^2$ for the Skorohod integral. These convergences are proven respectively by the means of Vitali's Theorem for $\sigma$-finite measures, and using a weak convergence Theorem in Banach--Saks spaces from \cite{jakszto2010another}.

\paragraph{On the space-time joint continuity of the flow. }
In contrast with the setting of  \cite{delmoral2022} and \cite{ItoAG}, the flow process $((s,t),x) \mapsto X_{s,t}^x$ and its derivatives are \textit{not almost surely jointly continuous}, since the Poisson random measure introduces jumps in the dynamics. Trying to extend the definition of being "right-continuous with left limits" for a function with variable in more than one dimension is quite technical, and such a hypothesis would be challenging to check in the SDE setting.
Fortunately, it turns out that both Vitali's Lemma and \cite[Theorem 4]{jakszto2010another} do not require a.s continuity, but only continuity \textit{in probability}. Thanks to this observation and combined with the fact that Lévy process (and more generally, SDE driven by Poisson random measures) are continuous in probability with respect to the time variable, we manage to rely on the hypothesis that for every $t \in [0,T]$, the process $(s,x) \mapsto X_{s,t}^x$ is continuous in probability. This is a quite mild hypothesis to check when dealing with SDEs driven by a tempered stable process, where the joint continuity in probability can be directly derived from the continuity for the $L^2$-norm (see Section \ref{sec:Example} for an example). 

\subsection{Organization of the paper}

In Section \ref{sec:Poisson_AG_Formula}, we formulate our hypotheses and the precise statement of our main result, Theorem \ref{maintheo}. We give an example of application, with SDEs driven by truncated $\alpha$-stable processes and smooth coefficients. 
Then, in Section \ref{sec:Malliavin} we define the appropriate tools to construct a Skorohod--Poisson integral with additional information. We also provide conditions for the convergence of sequences of Skorohod--Poisson integrals.
In Section \ref{sec:Poisson_AG_Formula} we derive an Itô formula with respect to Poisson stochastic integrals with a random test function, which is assumed to be independent from the underlying stochastic process. 
Finally, Section \ref{sec:proof_maintheo} is devoted to the proof of Theorem \ref{maintheo}. The proof of the technical lemmas and other estimates are written in Section \ref{sec:Technical_Lemmas}.

\subsection*{Notations}

\begin{itemize}
    \item $\R$ is the space of all real numbers.
    \item $\N = \{0,1,2, \ldots\}$ is the space of all non-negative integers.
    \item If $(\Lambda,\mathcal{A})$ is a measured space and $\mathcal{A}_1,\mathcal{A}_2 \subset \mathcal{A}$, we note $\mathcal{A}_1 \vee \mathcal{A}_2 \defeq \sigma(\mathcal{A}_1 \cup \mathcal{A}_2)$.
    \item If $(\Lambda,\mathcal{A},\mu)$ is a measure space and $p \in [1, +\infty)$, we denote by $L^p(\Lambda,\mathcal{A},\mu)$ the space of every real-valued measurable functions $f:\Lambda \to \R$ satisfying 
    \begin{equation*}
        \int_{\Lambda} |f|^p \ \d \mu < + \infty, 
    \end{equation*}
    which is a Banach space equipped with the norm $\| \cdot \|_{L^p(\Lambda,\mathcal{A},\mu)}$, defined by
    \begin{equation*}
        \| f \|_{L^p(\Lambda,\mathcal{A},\mu)} \defeq \left( \int_{\Lambda} |f|^p \ \d \mu\right)^{\frac{1}{p}}
    \end{equation*}
    for every $f \in L^p(\Lambda,\mathcal{A},\mu)$.
    \item If $(\Lambda,\mathcal{A},\mu)$ is a measure space we denote by $L^\infty(\Lambda,\mathcal{A},\mu)$  the space of every real-valued measurable functions $f:\Lambda \to \R$ satisfying
    \begin{equation*}
       \sup \text{ess}  |f|  < + \infty, 
    \end{equation*}
    which is a Banach space equipped with the norm $\| \cdot \|_{L^{\infty}(\Lambda,\mathcal{A},\mu)}$, defined by
     \begin{equation*}
        \| f \|_{L^\infty(\Lambda,\mathcal{A},\mu)} \defeq \sup \text{ess}  |f|
    \end{equation*}
    for every $f \in L^\infty(\Lambda,\mathcal{A},\mu)$.
    \item When the confusion is not possible, we may denote $L^p(\Lambda,\mathcal{A},\mu)$ with $p \in [1, +\infty]$ by $L^p(\mu)$, or $L^p(\Lambda)$, or $L^p(\mathcal{A})$ depending on the situation. 
    \item For $d \in \N \backslash \{ 0 \}$ and $A \subset \R^d$, $\mathcal{B}(A)$ designates the Borel $\sigma$-algebra on $A$.
    \item We say that a real function $f: \R \to \R$ is of class $C^k$ for $k \in \N$ if $f$ is differentiable $k$ times and if its $k$-th derivative is continuous. We sometimes also write $f \in C^k(\R \to \R)$. 
\end{itemize}

\section{The Poisson--Alekseev--Gröbner formula} \label{sec:Poisson_AG_Formula}

Let $T > 0$, $(\Omega,\Ff,\P)$ be a probability space, $N$ a Poisson random measure on $\R \backslash \{0 \} \times [0,T]$ with Lévy measure $\nu$, so that we have
\begin{equation*}
    \forall A \in \mathcal{B}(\R \backslash \{0 \}), \ \forall t \in [0,T], \quad \E[ N([0,t], A) ] = t \ \nu(A).
\end{equation*}
We denote by $\tilde{N}$ the compensated Poisson measure associated to $N$, defined by $\tilde{N}(\d t, \d z) = N(\d t, \d z) - \nu (\d z) \d t$. Let $\mathcal{N}$ be the set of the negligible sets. We consider the filtration $(\mathbf{F}_t)_{t \geq 0}$ generated by $N$ and the negligible sets. Precisely, we set $ \mathbf{F}_t = \sigma( N([0,s], A), \ s \in [0,t], \ A \in \mathcal{B}(\R \backslash \{0 \})) \vee \mathcal{N}$ for every $t \in (0,T]$ and assume that $\mathbf{F}_0$ is independent from $\sigma( N([0,t], A), \ t\in [0,T], \ A \in \mathcal{B}(\R \backslash \{0 \}))$.  
We consider measurable deterministic functions $\mu: [0,T] \times \R \to \R$, $\sigma: [0,T] \times \R \times \R \backslash \{0 \} \to \R$, a $\mathbf{F}_0$-measurable random variable $Y_0$, and stochastic processes $(A_t)_{t \in [0,T]}$ and $(B_{t,z})_{t \in [0,T], z \in \R \backslash \{0 \}}$. We assume that $A$ and $B_{\cdot,z}$ are $\mathbf{F}$-predictable for almost every $z \in \R \backslash \{0 \}$. 

We consider the flow process $X$ defined, for $0 \leq s \leq t \leq T$ and $x \in \R$, by
\begin{equation*}
    X_{s,t}^x = x + \int_s^t \mu \left( r,X_{s,r}^x \right) \ \d r + \int_s^t \int_{\R \backslash \{ 0 \}} \sigma \left( r, X_{s,r^-}^x,z \right) \ \tilde{N}(\d r , \d z), 
\end{equation*}
so that for any $s \in [0,T]$ and $x \in \R$, the stochastic process $(X_{s,t}^x)_{t \in [s,T]}$ is $\mathbf{F}$-adapted.

We also consider the $\mathbf{F}$-adapted stochastic process $Y$, defined for every $t \in [0,T]$ by
\begin{equation*}
    Y_t = Y_0 + \int_0^t A_r \ \d r + \int_0^t \int_{\R \backslash \{ 0 \}} B_{r,z} \ \tilde{N}(\d r , \d z).  
\end{equation*}
Finally, we consider a measurable function $f: \R \to \R$, which we refer to as the test function. 

\subsection{Hypotheses} \label{sec:Hypotheses}

Let $p > 4$, $q \in (0,\frac{p}{2}-2)$ and $k \geqslant 1$ be real numbers. 
We consider the following real constants:
\begin{itemize}
    \item $m_{\mu},\mx, \mxx, \mdx, \mdxx, \mddx, \mddxx \geqslant 1$,
    \item  $\ms, \msx, \msxx, \msdx, \msdxx, \msddx, \msddxx \geqslant 1$, and
    \item $C_{\mu}, C_f, C_\sigma, C_{X,p}, C_{X,{\mx}}, C_{X,{\mdx}}, C_{X,\mddx}, C_{X,p, \sigma}, C_{X,{\mx}, \sigma}, C_{X,{\mdx}, \sigma}, C_{X,\mddx, \sigma} > 0$.
\end{itemize}

\subsubsection{On the coefficients and test functions}

We start by giving the hypothesis on the SDEs coefficients and the test function $f$. 

\begin{Hyp}{\Blue{\bf \textrm Regularity of the coefficients. (\ref{hyp:coeffs}).}}
\makeatletter\def\@currentlabel{{\bf\textrm H}$_{\mbox{\scriptsize\bf\textrm{Coeffs}}}$}\makeatother
\label{hyp:coeffs}
We assume that there exists a positive predictable stochastic process $(\overline{B}_r )_{r \in [0,T]}$ such that
\begin{equation*}
    \forall r\in [0,T], \ \forall z \in \R\backslash\{0\}, \ |B_{r,z}| \leqslant |z|^k \ \overline{B}_r \quad \text{and} \quad \sup_{r \in [0,T]} \E \left[  (\overline{B}_r)^{p \vee \mx \vee \mdx \vee \mddx } \right]  < +\infty.
\end{equation*}
We also assume that
\begin{equation*}
    \forall r\in [0,T], \ \sup_{r \in [0,T]} \E \Big[ |A_r|^p \Big] < + \infty,
\end{equation*}
\begin{equation*} 
    \forall r\in [0,T], \ \forall z \in \R\backslash\{0\}, \ \forall x \in \R, \ |\sigma(r,x,z)| \leqslant C_\sigma |z|^k (1+|x|^{\ms}),
\end{equation*}
and
\begin{equation*}
     \forall r\in [0,T], \  \forall x \in \R, \ |\mu(r,x)| \leqslant C_{\mu} ( 1+ |x|^{\mm}).
\end{equation*}
\end{Hyp}

\begin{Hyp}{\hspace{-0.15cm}\Blue{\bf \textrm Regularity of the test function $f$. (\ref{hyp:testfunction}).}}
\makeatletter\def\@currentlabel{{\bf\textrm H}$_{\mbox{\scriptsize\bf\textrm{TestFunction}}}$}\makeatother
\label{hyp:testfunction}
We assume that $f$ is of class $C^2$ and satisfies
\begin{equation*}
    \forall x \in \R, \ \frac{|f(x)|}{1+|x|} + |f'(x)| + |f''(x)| \leqslant C_f (1+|x|^q).
\end{equation*}
\end{Hyp}

\subsubsection{On the flows and processes}

\begin{Hyp}{\hspace{-0.15cm}\Blue{\bf \textrm Flow property. (\ref{hyp:PushTheFlow}).}}
\makeatletter\def\@currentlabel{{\bf\textrm H}$_{\mbox{\scriptsize\bf\textrm{FlowProperty}}}$}\makeatother
\label{hyp:PushTheFlow}
    We suppose that for almost every $0 \leqslant s \leqslant t \leqslant T$, $x \mapsto X_{s, t}^x$ is $C^2$ on $\R$. We also suppose that for every $t \in [0, T]$ and $(s,x) \in [0,t] \times \R$, the following equality holds $\P$-almost surely:
\begin{equation*} 
    X_{t,T}^{X^x_{s,t}} = X_{s,T}^x.
\end{equation*}
\end{Hyp}

\begin{Hyp}{\hspace{-0.15cm}\Blue{\bf \textrm Stochastic continuity of the flow.(\ref{hyp:FlowContinuity}).}}
\makeatletter\def\@currentlabel{{\bf\textrm H}$_{\mbox{\scriptsize\bf\textrm{FlowContinuity}}}$}\makeatother
\label{hyp:FlowContinuity}
We suppose that for $\beta \in \{0, 1\}$, for every $t \in [0, T]$, the stochastic process of two variables $(s, x) \in [0, t] \times \R  \mapsto  \partial_x^{\beta} X^x_{s, t}$ is continuous in probability: for every $s \in [0, t]$ and $x \in \R$, for every $\varepsilon > 0$
\[ \P \Big( \big| \partial_x^{\beta} X_{s+\delta, t}^{x+h} - \partial_x^{\beta} X_{s, t}^x \big| > \varepsilon \Big) \xrightarrow[(\delta, h) \to (0, 0)]{} 0. \]  
\end{Hyp}

\begin{Hyp}{\hspace{-0.15cm}\Blue{\bf \textrm Moments control of the flow. (\ref{hyp:MomentsX}).}}
\makeatletter\def\@currentlabel{{\bf\textrm H}$_{\mbox{\scriptsize\bf\textrm{Moments($\partial_x^{\beta}X$)}}}$}\makeatother
\label{hyp:MomentsX}
For every triple $(M, N, C)$ in
\begin{multline*}
    \left\{ \left( p,\mx,C_{X,p} \right),\left( \mx,\mxx,C_{X,\mx} \right), \left( \mdx,\mdxx,C_{\partial X, \mdx} \right), \right. \\ \left. \left( \mddx,\mddxx, C_{\partial^2 X, \mddx} \right) , \left( \ms p, \msx, C_{X, p, \sigma}  \right), \left( \ms \msx, \msxx, C_{X, \mx, \sigma}  \right),  \right. \\
    \left. \left( \ms \msdx, \msdxx, C_{X, \mdx, \sigma}  \right), \left(\ms \msddx,\msddxx, C_{\partial^2 X, \msddx} \right) \right\},
\end{multline*}
we assume that
\begin{equation*} 
    \forall x \in \R ,  \ \sup_{0 \leqslant s \leqslant t \leqslant T} \E \Big[ | X_{s,t}^x|^{M} \Big] \leqslant C \ \left( 1 + |x|^{N} \right), \\
\end{equation*}
\begin{equation*} 
\forall x \in \R , \ \sup_{0 \leqslant s \leqslant t \leqslant T} \E \left[ 
 | \partial_x X_{s,t}^x|^{\frac{4p}{p-2(q+1)} } \right] \leqslant C_{\partial_X} \ \left( 1 + |x|^{\mdx} \right),
\end{equation*}
\begin{equation*}
\forall x \in \R , \ \sup_{0 \leqslant s \leqslant t \leqslant T} \E \left[ 
 | \partial^2_{x} X_{s,t}^x|^{\frac{2p}{p-2(q+1)} } \right] \leqslant C_{\partial^2_X} \ \left( 1 + |x|^{\mddx} \right),
\end{equation*}
\end{Hyp}

\begin{Hyp}{\hspace{-0.15cm}\Blue{\bf \textrm Moments control of the process $Y$. (\ref{hyp:MomentsY}).}}
\makeatletter\def\@currentlabel{{\bf\textrm H}$_{\mbox{\scriptsize\bf\textrm{Moments($Y$)}}}$}\makeatother
\label{hyp:MomentsY}
For 
\begin{multline*}
    M = \max \left\{\mx,\mxx,\mdx,\mdxx,\mddx,\mddxx, \right. \\ \left. \msx,\msxx,\msdx,\msdxx,\msddx,\msddxx \right\},
\end{multline*}
we suppose
\begin{equation*} 
    \sup_{0 \leq r \leq T} \E \Big[ |Y_r|^{M} \Big]   < + \infty.
\end{equation*}
\end{Hyp}

\subsubsection{On the Lévy measure}

Finally, we make the following assumptions on the moments of the Lévy measure $\nu$.

\begin{Hyp}{\hspace{-0.15cm}\Blue{\bf \textrm Moments of $\nu$. (\ref{hyp:mesure_levy}).}}
\makeatletter\def\@currentlabel{{\bf\textrm H}$_{\mbox{\scriptsize\bf\textrm{Moments($\nu$)}}}$}\makeatother
\label{hyp:mesure_levy}
For any \begin{equation*}
    \eta \in \left[ 2, \ k \ \bigg( 2 + \frac{1}{p} \Big(
  2 \mx q + (\mdx \vee \mddx)(p-2(1+q) \Big) \bigg)  \right],
\end{equation*} we assume that
\begin{equation*}
    \int_{\R \backslash \{ 0 \}} |z|^{\eta} \nu( \d z) + \int_{\R \backslash \{ 0 \}} |z|^{2\eta} \nu( \d z) \ < + \infty.
\end{equation*}
\end{Hyp}
\begin{Remarque}
    Although this requirement appears technical, it reduces to 
\begin{equation*}
    \int_{\R \backslash \{ 0 \}} |z|^{2} \nu( \d z) + \int_{\R \backslash \{ 0 \}} |z|^{4} \nu( \d z) \ < + \infty
\end{equation*}
for the case of an SDE driven by a Lévy process with coefficients of class $C^2$ in space with uniformly bounded derivatives. In fact, in this situation we have $\sigma(r,x,z) = z \ \tilde{\sigma}(r,x)$, leading to $k = 1$, and due to \cite[Theorem 5.1]{BretonPrivault}, every inequality in \ref{hyp:MomentsX} holds with $N = \mx = \mdx = \mddx = 0$. Moreover, as mentioned in the introduction we are particularly interest in the case where the Lévy measure has a tempered $\alpha$-stable shape, and in this case, we have 
\begin{equation*}
   \forall \eta \geqslant 2, \ \int_{\R \backslash \{ 0 \}}  |z|^{\eta} < + \infty,
\end{equation*}
so \ref{hyp:mesure_levy} is obviously satisfied.
\end{Remarque}

\subsection{Main result} \label{sec:MainResult}

We may now state our main result, which we call the Poisson--Alekseev--Gröbner formula. The proof is postponed in Section \ref{sec:proof_maintheo}.

\begin{Theoreme} \label{maintheo}
Assume that hypotheses \ref{hyp:coeffs}, \ref{hyp:testfunction}, \ref{hyp:PushTheFlow}, \ref{hyp:FlowContinuity}, \ref{hyp:MomentsX}, \ref{hyp:MomentsY} and \ref{hyp:mesure_levy} are satisfied. Then the following identity holds $\P$-almost surely:
\begin{multline} \label{eq:ito_AG}
      f \left( X_{0,T}^{Y_0} \right) - f(Y_T) = \int_{0}^{T} f' \left( X_{r,T}^{Y_{r}} \right) \partial_x X_{r,T}^{Y_{r}} \left(  \mu \left( r, Y_r \right) - A_r \right) \ \d r   \\
    + \int_{0}^{T} \int_{\R \backslash \{ 0 \}} \left\{ f \left( X_{r,T}^{Y_r + \sigma \left( r,Y_r,z \right)} \right) - f \left( X_{r,T}^{Y_{r}+B_{r,z}} \right) - f' \left( X_{r,T}^{Y_r} \right) \partial_x X_{r,T}^{Y_r} \Big(    \sigma \left( r,Y_r,z \right) - B_{r,z} \Big) \right\} \nu( \d z) \d r \\
     + \int_0^T \int_{\R \backslash \{ 0 \}} \left( f \left( X_{r,T}^{Y_{r} + \sigma \left( r,Y_{r},z \right)} \right) - f \Big( X_{r,T}^{Y_{r}+B_{r,z}} \Big) \right) \ \tilde{N}(\delta r, \delta z),
\end{multline}
where the integral against $\tilde{N}(\delta r, \delta z)$ stands for the Poisson--Skorohod integral, defined in Definition \ref{def:Skorohod_integral} in Section \ref{sec:Malliavin} below (precisely, $\tilde{N}(\delta r, \delta z)$ designates $\tilde{N}(\delta r, \delta z)^{\mathbf{F}_S}$ with $S = 0$).
\end{Theoreme}

\begin{Remarque}
Theorem \ref{maintheo} is stated in one dimension, however we don't see any difficulty with its generalization to multiple dimensions, by modifying accordingly \ref{hyp:coeffs} with the adequate norms. The equation \eqref{eq:ito_AG} would remain exactly the same, with the replacement of $\R \backslash \{0 \}$ by $\R^d \backslash \{0 \}$ for a $d$-dimensional Poisson random measure $N$. We chose to stay in a one-dimensional setting only to lighten the proof of the estimates in Section \ref{sec:Technical_Lemmas} which are already technical. 
\end{Remarque}

\subsection{A simple example} \label{sec:Example}

Let $X_0 \in \bigcap_{p \geqslant 2} L^p(\Omega)$, consider measurable functions $b, \overline{b}, \sigma, \overline{\sigma}: \R \to \R$, $\alpha \in (0,2)$ and a $1$-truncated $\alpha$-stable Lévy process $L_t$, i.e
\begin{equation*}
    L_t = \int_0^t \int_{\R \backslash \{ 0 \}} z \tilde{N}(\d s,\d z)
\end{equation*}
where $N$ is a random Poisson measure with intensity $\nu(\d z) = |z|^{-(\alpha+1)} \mathbf{1}_{ \{ |z| \leq 1 \} } \d z$. 
Consider the flow process $(X_{s,t}^x)_{0 \leq s \leq t \leq T, x \in \R}$ defined\footnote{Observe that the stochastic integral $\int_0^t \sigma(X_{s^-,r}^x) \d L_r$ is just equivalent to $\int_0^t \sigma(X_{s^-,r}^x) z \tilde{N}(\d r, \d z)$} by
\begin{equation} \label{eq:SdeFlowExample}
    X_{s,t}^x = x + \int_0^t b \left( X_{s,r}^x \right) \ \d r + \int_0^t \sigma \left( X_{s^-,r}^x \right) \ \d L_r,
\end{equation}
and the stochastic process $(Y_t)_{t \in [0,T]}$ defined as
\begin{equation} \label{eq:SdeApproxExample}
    Y_t = X_0 + \int_0^t \overline{b}(Y_r) \ \d r + \int_0^t \overline{\sigma}(Y_{r^-}) \ \d L_r.  
\end{equation}
Assume that the maps $b, \overline{b}, \sigma, \overline{\sigma}$ are uniformly bounded, and $C^2$-differentiable with first and second derivatives uniformly bounded. Let $f: \R \to \R$ a $C^2$-differentiable function such that 
\begin{equation*}
    \forall x \in \R, \frac{|f(x)|}{1+|x|} + |f'(x)| + |f''(x)| \leqslant C(1 + |x|^q) 
\end{equation*}
for some $C >0 $ and $q > 0$. In particular, \eqref{hyp:coeffs} is satisfied with $k=1$ and $m_{\mu} = m_{\sigma} = 0$ and \ref{hyp:testfunction} is satisfied. Then, using e.g \cite[Theorem 2.11]{kunita}, the SDEs \eqref{eq:SdeFlowExample} and \eqref{eq:SdeApproxExample} have a unique strong solution which lies in $L^2(\Omega)$. Thanks to \cite[Proposition 4]{BretonPrivault}, we have that $x \mapsto X_{s,t}^x$ is a function of class $C^2$, and the derivative process $(\partial_x^\alpha X_{s,t}^x)_{t \in [s,T]}$ are solution of stochastic differential equations. 
Moreover, due to \cite[Theorem 5]{BretonPrivault}, for any $p > 2q+4$ and $\beta \in \{ 0, 1, 2 \}$, we have
\begin{equation*}
    \sup_{x \in \R} \sup_{0 \leqslant s \leqslant T} \E \left[ \sup_{s \leqslant t \leqslant T} \left| \partial_x^{\beta} X_{s,t}^x \right|^p \right] < + \infty.
\end{equation*}
Therefore,  \ref{hyp:MomentsX} is satisfied with $\mx = \mxx = \mdx = \mdxx = \mddx = \mddxx = 0$ and $\msx = \msxx= \msdx = \msdxx=\msddx = \msddxx = 0$, and \ref{hyp:MomentsY} is trivially satisfied. 
Since we have $\int_{\R \backslash \{0 \}} |z|^\eta \nu (\d z) < + \infty$ for every $\eta \geqslant 2$, \ref{hyp:mesure_levy} is also satisfied. 
In addition, we have the following standard result, that we prove in Section \ref{subsec:ProofL2Continuity} for the sake of completeness.

\begin{Lemma} \label{lm:L2ContinuityExampleProcess}
    For any $t \in [0,T]$ the bi-variate stochastic processes $(s,x) \mapsto X_{s,t}^x$ and $(s,x) \mapsto \partial_x X_{s,t}^x$ are continuous for the $L^2$-norm.
\end{Lemma}

This implies that \ref{hyp:FlowContinuity} holds. 
Finally, we obtain \ref{hyp:PushTheFlow} using the strong uniqueness of the solution. Specifically, for any $0 \leqslant s \leqslant t \leqslant T$ and $x \in \R$ we have that
\begin{align*}
    X_{s,T}^x  & = x + \int_s^T b(X_{s,r}^x)   \d r + \int_s^T \sigma(X_{s,r^{-}}^x)  \d L_r  \\
    & = x + \left(\int_s^t b(X_{s,r}^x)  \d r + \int_s^t \sigma(X_{s,r^{-}}^x)  \d L_r \right) + \int_t^T  b(X_{s,r}^x)  \d r + \int_t^T  \sigma(X_{s,r^{-}}^x)  \d L_r \\ 
    & = X_{s,t}^x + \int_t^T  b(X_{s,r}^x)  \d r + \int_t^T  \sigma(X_{s,r^{-}}^x)  \d L_r,
\end{align*}
while by definition, 
\begin{equation*}
    X_{t,T}^{X_{s,t}^x} = X_{s,t}^x + \int_t^T  b \left( X_{s,r}^{X_{s,t}^x} \right) \ \d r + \int_t^T  \sigma \left( X_{s,r^{-}}^{X_{s,t}^x} \right) \ \d L_r,
\end{equation*}
so the processes $(X_{s,r}^x)_{r \in [t,T]}$ and $(X_{t,T}^{X_{s,t}^x})_{r \in [t,T]}$ satisfy the same SDE with same initial condition, hence they are almost surely equal by strong uniqueness. 
We have checked all the conditions to apply the Poisson--Alekseev--Gröbner formula, which gives
\begin{multline*}
    f \left( X_{0,T}^{Y_0} \right) - f(Y_T) = \int_0^T f' \left( X_{r,T}^{Y_r} \right) \ \partial_x X_{r,T}^{Y_r} \ \left( b(Y_r) - \overline{b}(Y_r) \right) \ \d r \\  + \int_0^T \int_{\R \backslash \{0 \}} \left\{ f \left( X_{r,T}^{Y_r + \sigma(Y_r) z} \right) - f \left( X_{r,T}^{Y_r + \overline{\sigma}(Y_r) z } \right) - f' \left( X_{r,T}^{Y_r} \right) \ \partial_x X_{r,T}^{Y_r} \ \left( \sigma(Y_r) - \overline{\sigma}(Y_r) \right) z \right\} \nu(\d z) \d r  \\
     + \int_0^T \int_{\R \backslash \{ 0 \}} \left\{ f \left( X_{r,T}^{Y_r + \sigma(Y_r) z} \right) - f \left( X_{r,T}^{Y_r + \overline{\sigma}(Y_r) z} \right) \right\} \ \tilde{N}(\delta r , \delta z).
\end{multline*}
In particular, taking the expectation and applying Taylor's formula we obtain that
\begin{multline*}
    \E[f(X_{0,T}^{Y_0})] - \E[f(Y_T)] = \int_0^T \E[f'(X_{r,T}^{Y_r}) \partial_x X_{r,T}^{Y_r} \ (b(Y_r) - \overline{b}(Y_r))] \d r \\  +\int_0^T \int_{\R \backslash \{0 \}} z \ \E \bigg[ F_{r, T, z} ( \sigma(Y_r) - \overline{\sigma}(Y_r) ) \bigg] \nu(\d z) \d r,
\end{multline*}
with
\[ F_{r, T, z} = \int_0^1 \bigg\{ f'(X_{r,T}^{Y_r + z \overline{\sigma}(Y_r) + \lambda z (\sigma(Y_r) - \overline{\sigma}(Y_r))}) \partial_x X_{r,T}^{Y_r + z \overline{\sigma}(Y_r) + \lambda z (\sigma(Y_r) - \overline{\sigma}(Y_r))} - f'(X_{r,T}^{Y_r}) \partial_x X_{r,T}^{Y_r} \bigg\} \  \d \lambda. \]
This allows to express the weak global error $\E[f(X_{0,T}^{Y_0})] - \E[f(Y_T)]$ in terms of the local error between the coefficients $b-\overline{b}$, $\sigma - \overline{\sigma}$, as long as one can decouple the latter from the multiplicative stochastic terms (e.g via Hölder inequality, or Malliavin calculus).
%It would be tempting to apply Taylor's formula a second time to the integral with respect to $\nu(\d z) \d r$, but that would unfortunately not enhance a convergence result since 
%\begin{multline*}
%    f'(X_{r,T}^{Y_r + z \overline{\sigma}(Y_r) + \lambda z (\sigma(Y_r) - \overline{\sigma}(Y_r))}) \partial_x X_{r,T}^{Y_r + z \overline{\sigma}(Y_r) + \lambda z (\sigma(Y_r) - \overline{\sigma}(Y_r))} - f'(X_{r,T}^{Y_r}) \partial_x X_{r,T}^{Y_r} \\ = z ( \overline{\sigma}(Y_r) + \lambda (\sigma(Y_r) - \overline{\sigma}(Y_r)) \ \int_0^1 \frac{\partial^2}{\partial x^2} f(X_{r,T}^x) \bigg|_{x = Y_r + z \overline{\sigma}(Y_r) + \lambda z (\sigma(Y_r) - \overline{\sigma}(Y_r)) + \mu \lambda z ( \overline{\sigma}(Y_r) + \lambda (\sigma(Y_r) - \overline{\sigma}(Y_r))} \d \mu.
%\end{multline*}

\section{An Itô formula for random independent functional of Poisson stochastic integrals}

Let $S < T$ and $(\mathbf{H}_t, t \in [S, T])$ a filtration such that for every $t \in [S, T]$, $\mathbf{F}_{[S,t]} \subset \mathbf{H}_t$, where we recall from \eqref{def:filtrations_bizarres} in Section \ref{sec:Malliavin} that 
\[ \mathbf{F}_{[S,t]} = \mathcal{N} \vee \sigma \big( N([S, s], A), S \leqslant s \leqslant t, A \in \mathcal{B}(\R \backslash \{ 0 \}) \big), \]
where $\mathcal{N}$ is the set of $\mathcal{F}$-null sets. We let $\mathcal{P}_2^{\mathbf{H}}(S,T)$ be the space of $(\mathbf{H}_t)_{t \in [S, T]}$-predictable mappings $F: [S,T] \times \R \backslash \{ 0 \} \times \Omega \longrightarrow \R$ such that,
\begin{equation}\label{eqIto:defP2}
    \int_S^T \int_{\R \backslash \{ 0 \}} \E \left[ F(s,z)^2 \right] \ \nu(\d z) \ \d s < +\infty.
\end{equation}

We denote its Itô integral with respect to the filtration $\mathbf{H} = (\mathbf{H}_t, t \in [S, T])$ by
\[ \int_{S}^T \int_{\R \backslash \{ 0 \}} F(t, x) \ \tilde{N}^{\mathbf{H}} (\d t, \d z). \]
If $\mathbf{H} =\mathbf{F}_{[S,\cdot]}$, then we only denote it $\int_{S}^T \int_{\R \backslash \{ 0 \}} F(t, z) \ \tilde{N} (\d t, \d z)$. If $u$ is $\mathbf{F}_{[S,\cdot]}$-adapted, then it is $\mathbf{H}$-adapted and 
\[ \int_{S}^T \int_{\R \backslash \{ 0 \}} F(t, z) \ \tilde{N}^{\mathbf{H}} (\d t, \d z) = \int_{S}^T \int_{\R \backslash \{ 0 \}} F(t, x) \ \tilde{N} (\d t, \d z). \]
We now assume that $\mathbf{H}$ satisfies
\begin{equation*}
    \forall t \in [S,T], \quad \mathbf{H}_t = \mathbf{F}_{[S,t]}\vee \mathbf{G}_T
\end{equation*}
where $\mathbf{G}_T$ is a $\sigma$-algebra independent from $\mathbf{F}_{[S,t]}$ for every $t \in [S,T]$. We consider a $(\mathbf{F}_{[S,t]})_{t \in [S,T]}$-predictable stochastic process $A:[S,T] \times \Omega \xrightarrow{} \R$ such that $\P$-a.s, $\int_{S}^T |A_r| \ \d r < \infty$, a stochastic process $B \in \mathcal{P}_2^{\mathbf{F}_{[S,\cdot]}}(S,T)$ and a real $x\in \R$. We introduce the stochastic process $(X_t)_{t \in [S,T]}$ defined by
\begin{equation*}
    X_t \defeq x + \int_S^t A_r \ \d r + \int_S^t \int_{\R \backslash \{ 0 \}} B_{r,z} \ \tilde{N}(\d r, \d z).
\end{equation*}
Let $f: \Omega \times \R \to \R$ be a random functional such that:
\begin{itemize}
    \item for every $x \in \R$, $f(\cdot,x)$ is $\mathbf{G}_T$-measurable,
    \item for almost every $\omega \in \Omega$, $f(\omega,\cdot) \in \mathcal{C}^2(\R \to \R)$. 
\end{itemize}
Then we have the following Itô formula:

\begin{Proposition} \label{prop:ito_formula_random_functions_poisson} For every $t \in [S,T]$, 
\begin{multline*}
    f(X_t) - f(X_S) =  \int_S^t f'(A_r) \ \d r + \int_S^t \int_{\R \backslash \{0 \}} \Big( f(X_r + B_{r,z}) - f(X_r) - B_{r,z} f'(X_r) \Big) \ \nu (\d z) \d r \\ + \int_S^t\int_{\R \backslash \{0 \}} \Big( f(X_{r^-} + B_{r,z}) - f(X_{r^-}) \Big) \ \tilde{N}^{\mathbf{H}}(\d r , \d z ).
\end{multline*}
\end{Proposition}

\begin{Demo}
    This is a direct consequence of Proposition \ref{prop:ito_formula_indep_cadlag}, since $X$ is a $(\mathbf{F}_{[S,t]})_{t \in [S,T]}$-càdlàg semimartingale, and $\mathbf{G}_T$ is independent of $\mathbf{F}_{[S,t]}$ for every $t \in [S,T]$.
\end{Demo}

\section{Malliavin calculus for Poisson random integrals with additional information} \label{sec:Malliavin}

The aim of this section is to construct a Skorohod integral with respect to a Poisson process which satisfies a Chasles relation. To do this, we introduce an modification of the known Malliavin derivative, see for instance the set of papers in \cite{PeccatiReitzner}. We follow the construction in \cite{ItoAG}. We construct the Skorohod integral by being the adjoint of a derivative operator, which will be built with the so called smooth random variables.

\subsection{Smooth random variables}

Let $(\Omega, \mathcal{F}, \P)$ a probability space. Let $m \geqslant 1$ an integer, $S < T$ two real numbers (which could eventually be negative) and $N$ a Poisson measure on $E \defeq \R \backslash \{ 0 \} \times (S, T)$, with intensity $\nu$. We denote by $\tilde{N}$ the compensated Poisson measure associated, defined by $\tilde{N}(\d s, \d z) \defeq N(\d s, \d z) - \nu(\d z) \d s$. For $\phi \in L^2(E, \d s \otimes \nu(\d z))$, we define $\tilde{N}(\phi)$ as the Poisson stochastic integral of $\phi$ with respect to $\tilde{N}$, that is
\begin{equation*}
    \tilde{N}(\phi) = \int_S^T \int_{\R \backslash \{0 \}} \phi(s,z) \tilde{N}(\d s, \d z). 
\end{equation*}

We consider $\mathbf{F}_S \subset \mathcal{F}$ a $\sigma$-algebra which is independent of $\sigma(N([S, t], A), t \in [S, T], A \in \mathcal{B}(\R \backslash \{ 0 \}))$. Let $\mathcal{N}$ the set consisting of all the $\P$-negligible sets. Consider the following $\sigma$-algebra which will be the ambient one: $\mathbf{F}_T$ is the $\sigma$-algebra generated by $\mathcal{N}$, $\mathbf{F}_S$ and the family of random variables $N([S, t], A)$ for every $t \in [S, T]$ and $A \in \mathcal{B}(\R \backslash \{ 0 \})$.

When $\mathcal{G} \subset \mathcal{F}$ is a $\sigma$-algebra, and for $n \geqslant 1$, we define $\mathcal{P}^{\mathcal{G}} (\R^n)$ as the set of functions $f: \Omega \times \R^n \longrightarrow \R$ such that
\begin{enumerate}
    \item For almost every $\omega \in \Omega$, $x \in \R^n \longmapsto f(\omega, x)$ is a trigonometric polynomial; 
    \item For every $x \in \R^n$, $f_x: \omega \longmapsto f(\omega, x)$ is $\mathcal{G}$-measurable; 
    \item The map $(\omega, x) \longmapsto f(\omega, x)$ is $\mathcal{G}\otimes \mathcal{B}(\R^n)$-measurable.
\end{enumerate}

\begin{Definition}\label{DefSmoothST}
    We define the space of \emph{smooth random vectors with additional information in $S$} as the set of random vectors $F \in L^2 (\mathbf{F}_T)$ such that there exists $n \geqslant 1$, $f \in \mathcal{P}^{\mathbf{F}_S} (\R^n)$ and $\phi_1, \cdots, \phi_n \in L^2 (E)$ such that there exists $\Omega^* \in \mathbf{F}_T$ depending on all those parameters such that $\P(\Omega^*) = 1$ and for every $\omega \in \Omega^*$,
    \begin{equation}\label{EqDefSmoothST}
        F (\omega) = f \left( \omega, \left[ \int_S^T \int_{\R \backslash \{ 0 \}} \phi_1 (s, z) \ \tilde{N}(\d s, \d z) \right] (\omega), \cdots, \left[ \int_S^T \int_{\R \backslash \{ 0 \}} \phi_n (s, z) \ \tilde{N}(\d s, \d z) \right] (\omega) \right).
    \end{equation}
    We denote by $\mathcal{S} (\mathbf{F}_S, \tilde{N})$ those random vectors.
\end{Definition}

For $S \leqslant a \leqslant b \leqslant T$, we note 
\begin{equation} \label{def:filtrations_bizarres}
    \mathbf{F}_{[a, b]} \defeq \mathcal{N} \vee \sigma \Big( N([a, r], A), a \leqslant r \leqslant b, A \in \mathcal{B}(\R \backslash \{ 0 \}) \Big).
\end{equation}
For $S <s<t<T$, we can build a notion of smooth random vectors with additional information on $[S, s] \cup [t, T]$. We introduce the following $\sigma$-algebra
\[ \mathbf{F}_{[S, s] \cup [t, T]} \defeq \mathbf{F}_S \vee \mathbf{F}_{[S, s]} \vee \mathbf{F}_{[t, T]}. \]
In this case, we can define $\mathcal{S} (\mathbf{F}_{[S, s] \cup [t, T]}, \tilde{N}_{\left|_{(s, t)} \right.})$, where $\tilde{N}_{\left|_{(s, t)} \right.}$ is the restriction of $\tilde{N}$ on $\mathcal{B}(s, t) \otimes \mathcal{B}(\R \backslash \{ 0 \})$. Hence, $F \in \mathcal{S} (\mathbf{F}_{[S, s] \cup [t, T]}, \tilde{N}_{\left|_{(s, t)} \right.})$ if there exists $n \geqslant 1$, $f \in \mathcal{P}^{\mathbf{F}_{[S, s] \cup [t, T]}} (\R^n)$ and $\phi_1, \cdots, \phi_n \in L^2 ((s, t) \times \R \backslash \{ 0 \})$ such that there exists $\Omega^* \in \mathbf{F}_T$ with $\P(\Omega^*) = 1$ and for every $\omega \in \Omega^*$,
\begin{equation}\label{EqDefSmoothSstT}
        F (\omega) = f \left( \omega, \left[ \int_s^t \int_{\R \backslash \{ 0 \}} \phi_1 (r, z) \ \tilde{N}(\d r, \d z) \right] (\omega), \cdots, \left[ \int_s^t \int_{\R \backslash \{ 0 \}} \phi_n (r, z) \ \tilde{N}(\d r, \d z) \right] (\omega) \right).
\end{equation}

In the following, when we will express a statement about $\mathcal{S} (\mathbf{F}_S, N)$, it is in fact also true with $\mathcal{S} (\mathbf{F}_{[S, s] \cup [t, T]}, N_{|_{[s, t]}})$, with adapted notations.

\subsection{Malliavin derivative}

\begin{Definition}\label{DefDerivee}
    Let $F \in \mathcal{S} (\mathbf{F}_S, N)$ given by \eqref{EqDefSmoothST}. We define its \emph{Malliavin derivative} as the following real random process: for almost every $t \in [S, T]$, $\nu$-almost every $z \in \R \backslash \{ 0 \}$ and $\P_{|_{\mathbf{F}_T}}$-almost every $\omega \in \Omega$ 
    \begin{equation}\label{EqDefDeriveeST}
        \left[ \D_{(t, z)} (\mathbf{F}_S, N) \right] (\omega) \defeq F^+_{(t, z)} (\omega) - F (\omega),
    \end{equation}
    where 
    \[ F^+_{(t, z)} (\omega) \defeq f \left( \omega, \left[ \tilde{N}(\phi_1) \right] (\omega) + \phi_1(t, z), \cdots, \left[ \tilde{N}(\phi_n) \right] (\omega) + \phi_n (t, z) \right) . \]
\end{Definition}

Note that $F^+$ is well-defined : if $F = f(\tilde{N}(\phi_1), \cdots, \tilde{N}(\phi_n)) = g(\tilde{N}(\psi_1), \cdots, \tilde{N}(\psi_m))$, then $F^+$ defined for the $f$-definition is equal in $L^2$ to the $F^+$ given by the $g$-definition, by Mecke's formula (see \cite[Theorem 4.1]{last_lectures_2017}). Since our Definition of smooth random variable in this case is quite new, we will develop some proofs about it. They will enlighten why we chose those kind of smooth random variables.

\begin{Proposition}\label{PropIPP}
    Let $F \in \mathcal{S} (\mathbf{F}_S, N)$ like in \eqref{DefSmoothST} and $\phi \in L^2 (E)$. Then
    \begin{equation}\label{EqIPP1}
        \E \left[ F \tilde{N}(\phi) \right] = \E \left[ \langle \D F, \phi \rangle_{L^2 (E)} \right].
    \end{equation}
\end{Proposition}

\begin{Demo} Suppose first that $F$ has the form
\begin{equation}\label{EqDemIPPPolTrigo}
     F = e^{\mathrm{i} \gamma \tilde{N}(\phi)} \prod_{k=1}^n e^{\mathrm{i} \gamma_k \tilde{N}(\phi_k)},
\end{equation}
with $\gamma, \gamma_1, \cdots, \gamma_n \in \R$ which are deterministic. Then, we have in fact
\[ \E \left[ F \tilde{N}(\phi) \right] = - \mathrm{i} \frac{\partial}{\partial \gamma} \left\{ \E \left[ \exp \left( \mathrm{i} \tilde{N} \left( \gamma \phi + \sum_{k=1}^n \gamma_k \phi_k \right) \right) \right] \right\}. \]
Using the Lévy--Khintchine formula to express the characteristic function of the Poisson stochastic integrals with respect to $\tilde{N}$, we conclude that
\[ \E \left[ F \tilde{N}(\phi) \right] = \E[F] \int_{E} \phi(t, z) \left( e^{\mathrm{i} \gamma \phi (t, z)} \prod_{k=1}^n e^{\mathrm{i} \gamma_k \phi_k (t, z)} - 1 \right) \ \d t  \ \nu(\d z). \]
This conclude the proof in this case, since for almost every $(t, z) \in E$, 
\[ \D_{(t, z)} F = e^{\mathrm{i} \gamma \tilde{N}(\phi)} \prod_{k=1}^n e^{\mathrm{i} \gamma_k \tilde{N}(\phi_k)} \left( e^{\mathrm{i} \gamma \phi (t, z)} \prod_{k=1}^n e^{\mathrm{i} \gamma_k \phi_k (t, z)} - 1 \right). \]
This proves \eqref{EqIPP1} when $F$ is a trigonometrical polynomial with deterministic coefficients and frequencies. Now, suppose that $F$ has the form \eqref{EqDemIPPPolTrigo}, with $\gamma, \gamma_1, \cdots, \gamma_n$ being random variables which are $\mathbf{F}_S$-measurable. Then those random variables are independent of $\tilde{N}(\psi)$ for every $\psi \in L^2 (E)$. Hence, we have for $\alpha, \alpha_1, \cdots, \alpha_n \in \R$
\[  \E \left[ \left. F \tilde{N}(\phi) \right| \gamma = \alpha, \gamma_1 = \alpha_1, \cdots, \gamma_n = \alpha_n \right] = \E \left[ \tilde{N} (\phi) e^{\mathrm{i} \alpha \tilde{N}(\phi)} \prod_{k=1}^n e^{\mathrm{i} \alpha_k \tilde{N}(\phi_k)} \right].  \]
By the deterministic case, 
\[ \E \left[ \left. F \tilde{N}(\phi) \right| \gamma = \alpha, \gamma_1 = \alpha_1, \cdots, \gamma_n = \alpha_n \right] = \E \left[ \left\langle \D \left( e^{\mathrm{i} \alpha \tilde{N}(\phi)} \prod_{k=1}^n e^{\mathrm{i} \alpha_k \tilde{N}(\phi_k)} \right) ,  \phi \right\rangle_{L^2(E)} \right]. \]
By using once again the independence between the $\gamma$ variables and the $\tilde{N}(\phi)$ ones, we conclude that for every $\alpha, \alpha_1, \cdots, \alpha_n \in \R$,
\[ \E \left[ \left. F \tilde{N}(\phi) \right| \gamma = \alpha, \gamma_1 = \alpha_1, \cdots, \gamma_n = \alpha_n \right] = \E \left[ \left. \langle \D F, \phi \rangle_{L^2 (E)} \right| \gamma = \alpha, \gamma_1 = \alpha_1, \cdots, \gamma_n = \alpha_n \right]. \]
We conclude this case by taking the expectation. Consequently, we proved \eqref{EqIPP1} for every $F$ being a trigonometric monomial with random frequencies. For a general $F \in \mathcal{S} (\mathbf{F}_S, N)$, $F$ is linear combination of those monomials, with random coefficients being $\mathbf{F}_S$-measurable. By using the conditional expectation as before, we conclude in the fact that \eqref{EqIPP1} is true for every $F \in \mathcal{S} (\mathbf{F}_S, N)$. \end{Demo}

\begin{Lemma}\label{lem:DensiteSmoothRV}
    The space $\mathcal{S} (\mathbf{F}_S, N)$ is dense into $L^2 (\mathbf{F}_T)$ for the $L^2$ norm.
\end{Lemma}

\begin{Demo}
    We use the criterion by orthogonality. Let $G \in L^2 (\mathbf{F}_T)$ such that for every $F \in \mathcal{S} (\mathbf{F}_S, N)$, $\E[F G] = 0$. We want to prove that $G = 0$ $\P$-almost surely. It is enough to prove that $\E[G | \mathbf{F}_T] = 0$. 
    
    First, since every $\mathbf{F}_S$-measurable random variables belongs to $\mathcal{S} (\mathbf{F}_S, N)$, it is immediate that $\E[G|\mathbf{F}_S] = 0$. 
    
    Second, it is also true that $\E[G|\mathbf{F}_{[S, T]}] = 0$, where we recall that $\mathbf{F}_{[S, T]}$ is the $\sigma$-algebra generated by $N_t(A) - N_S (A)$, where $t \in [S, T]$ and $A \in \mathcal{B}(\R \backslash \{ 0 \})$. Indeed, by hypothesis, for every $n \geqslant 1$, for every (deterministic) $\gamma_1, \cdots, \gamma_n \in \R$ and $\phi_1, \cdots, \phi_n \in L^2 (E)$, we have 
    \[ \E \left[ G \prod_{k=1}^n e^{\mathrm{i} \gamma_k \tilde{N}(\phi_k)} \right] = 0. \]
    This means that the signed finite measure $\eta$ defined on $\R^n$ by
    \[ \forall B \in \mathcal{B}(\R^n), \quad \eta(B) \defeq \E \left[ G \mathbf{1}_{B} \left( \tilde{N}(\phi_1), \cdots, \tilde{N} (\phi_n) \right) \right] \]
    admits a null Fourier transform. By injectivity, $\eta$ is the zero measure. This shows that for every $C \in \mathbf{F}_{[S, T]}$ (since the $\sigma$-algebra generated by the variables $N_t(A) - N_S(A)$ is the same as the one generated by $\tilde{N}(\phi)$ for $\phi \in L^2(E)$), $\E[G \mathbf{1}_C] = 0$, and so that $\E[G | \mathbf{F}_{[S, T]}] = 0$.

    Finally, by doing the same procedure as before but with $G \mathbf{1}_A$, with $A \in \mathbf{F}_S$ instead of $G$, we conclude that for every $B \in \mathbf{F}_{[S, T]}$, $\E[G \mathbf{1}_A \mathbf{1}_B] = 0$. 
    
    We conclude that for every $A \in \mathbf{F}_S$ and $B \in \mathbf{F}_{[S, T]}$, 
    \[ \E \left[ G \mathbf{1}_{A \cup B} \right] = \E \left[ G \mathbf{1}_A \right] + \E \left[ G \mathbf{1}_B \right] - \E \left[ G \mathbf{1}_{A \cap B} \right] = 0. \]
    Since the $\sigma$-algebra generated by $\mathbf{F}_S$ and $\mathbf{F}_{[S, T]}$ is also the one generated by $A \cup B$ with $A \in \mathbf{F}_S$ and $B \in \mathbf{F}_{[S, T]}$, we conclude that $\E[G|\mathbf{F}_S, \mathbf{F}_{[S, T]}] = 0$. It is not a big deal to add that we also have $\E[G|\mathbf{F}_S, \mathbf{F}_{[S, T]}, \mathcal{N}] = 0$, meaning that $\E[G|\mathbf{F}_T] = 0$. 
\end{Demo}

\begin{Proposition}\label{PropFermabilitéDérivée}
    The operator 
    \[ \D: \mathcal{S} (\mathbf{F}_S, N) \subset L^2 (\mathbf{F}_T) \longrightarrow L^2 (E \times \Omega) \]
    is closable. We denote by $\mathbb{D}^{1,2}_{\mathrm{J}} (\mathbf{F}_S, N)$ the domain of the closed operator $\D$. In particular, $\mathcal{S} (\mathbf{F}_S, N)$ is dense into $\mathbb{D}^{1,2}_{\mathrm{J}} (\mathbf{F}_S, N)$ for the graph norm:
    \[ \forall F \in \mathbb{D}^{1,2}_{\mathrm{J}} (\mathbf{F}_S, N), \norm{F}{\mathbb{D}^{1,2}_{\mathrm{J}} (\mathbf{F}_S, N)}^2 \defeq \E \left[ F^2 \right] + \E \left[ \norm{\D F}{L^2(E)}^2 \right]. \]
\end{Proposition}

\begin{Demo}
    Let $(F_n)_{n \geqslant 1}$ a sequence of random variables in $\mathcal{S} (\mathbf{F}_S, N)$ such that $(F_n)_{n \geqslant 1}$ goes to zero in $L^2 (\mathbf{F}_T)$ and $(\D F_n)_{n \geqslant 1}$ goes to a process $u \in L^2 (E \times \Omega)$ in the associated $L^2$ norm. The operator $\D$ is closable when $u = 0$.

    To prove it, we use several times that the orthogonal of the ambient Hilbert space is zero. By the Integration by Parts formula (Proposition \ref{PropIPP}), for every $G \in \mathcal{S} (\mathbf{F}_S, N)$, $\phi \in L^2 (E)$ and $n \geqslant 1$,
    \[ \E \left[ G \langle \D F_n, \phi \rangle_{L^2 (E)} \right] = \E \left[ F_n G  \tilde{N}(\phi) \right] - \E \left[ F_n \langle G, \phi \rangle_{L^2(E)} \right] - \E \left[ \langle \D F_n \ \D G, \phi \rangle_{L^2 (E)} \right]. \]
    We take the limit as $n$ goes to $\infty$ and we get for every $G \in \mathcal{S} (\mathbf{F}_S, N)$ and $\phi \in L^2 (E)$
    \begin{equation}\label{eqMalliavin:Fermabilité}
        \E \left[ \langle u G^+, \phi \rangle_{L^2 (E)} \right] = 0,
    \end{equation}
    where $G^+$ is defined in Definition \ref{DefDerivee}. For every $G \in \mathcal{S} (\mathbf{F}_S, N)$, for fixed $(t, z) \in E$, let $G^-_{(t, z)}$ defined as $G^{+}_{(t, z)}$ by replacing "${}+ \phi_k (t, z)$" by "${}- \phi_k (t, z)$". Then, we have
    \[ G = \left( G_{(t, z)}^- \right)^+_{(t, z)}. \]
    This is well defined, since $G_{(t, z)}^- \in \mathcal{S} (\mathbf{F}_S, N)$ : if $f \in \mathcal{P}^{\mathbf{F}_S} (\R^n)$ then $f^-_{t,z}$ defined by 
    \[ f_{t, z}^- (x_1, \cdots, x_n) := f(\cdot - \phi_1(t, z), \cdots, \cdot - \phi_n (t, z)) \]
    is still in $\mathcal{P}^{\mathbf{F}_S} (\R^n)$. Hence, we conclude that for every $G \in \mathcal{S} (\mathbf{F}_S, N)$ and $\phi \in L^2 (E)$, by \eqref{eqMalliavin:Fermabilité}
    \[ \int_E \E[u(t, z) G] \ \phi(t, z) \ \d t \ \nu(\d z) = \int_E \E \left[ u(t, z) \left( G_{(t, z)}^- \right)_{(t, z)}^+ \right] \ \phi(t, z) \ \d t \ \nu(\d z) = 0. \]
    In particular, $\langle u, \phi \rangle_{L^2 (E)} \in \mathcal{S}(\mathbf{F}_S, N)^{\perp}$, so is equal to zero in $L^2 (\P)$, by Lemma \ref{lem:DensiteSmoothRV}. Since $L^2 (E)$ is separable ($\mathrm{Leb} \otimes \nu$ is $\sigma$-finite and $\mathcal{B}([0,T] \times \R \backslash \{ 0 \})$ is countably generated), this means that there exists $\Omega^* \in \mathbf{F}_T$ (not depending on $\phi$) such that $\P(\Omega^*) = 1$ and for all $\omega \in \Omega^*$ and $\phi \in L^2 (E)$, $\langle u(\omega), \phi \rangle_{L^2 (E)} = 0$. This concludes that $u = 0$ in $L^2 (\Omega \times E)$. % L'argument "L^2(E) séparable" résout le problème de dépendance de Omega^* par rapport à phi, j'espère. Peut-être qu'on aurait pu faire plus élémentaire ?
\end{Demo}

\begin{Proposition}
    Let $S \leqslant s \leqslant t \leqslant T$. We have the following inclusion
    \begin{equation}\label{eq:inclusionDomainePetiteDerivee}
         \mathbb{D}^{1,2}_{\mathrm{J}} (\mathbf{F}_S, N) \subset \mathbb{D}^{1,2}_{\mathrm{J}} \left( \mathbf{F}_{[S, s] \cup [t, T]}, N_{|_{[s, t]}} \right) .
    \end{equation}
    Moreover, we have $\mathrm{Leb}(s, t) \otimes \nu$-almost everywhere and $\P_{|_{\mathbf{F}_T}}$-almost surely
    \begin{equation}\label{eq:PetiteGrandeDerivee}
        \D (\mathbf{F}_S, N)_{\left|_{[s, t] \times \R \backslash \{ 0 \}} \right.}  = \D \left( \mathbf{F}_{[S, s] \cup [t, T]}, N_{|_{[s, t]}} \right) .
    \end{equation}
\end{Proposition}

\begin{Demo}
    By the same arguments as the proof of Lemma A.6 in \cite{ItoAG}, we can show that 
    \begin{equation}\label{EqDemPetiteDerivee}
       \mathcal{S} (\mathbf{F}_S, N) \subset \mathcal{S} \left( \mathbf{F}_{[S, s] \cup [t, T]}, N_{|_{[s, t]}} \right). 
    \end{equation}
    In particular, if $F \in \mathcal{S} (\mathbf{F}_S, N)$ is given by \eqref{EqDefSmoothST}, then we prove that $F \in \mathcal{S} \left( \mathbf{F}_{[S, s] \cup [t, T]}, N_{|_{[s, t]}} \right)$ by defining for $\omega \in \Omega$ and $x_1, \cdots, x_n \in \R$
    \begin{equation*}
        \begin{split}
            \\ & g(\omega, x_1, \cdots, x_n) \\ 
            \defeq & \ f \left( \omega, x_1 + \left[  \int_S^s \int_{\R \backslash \{ 0 \}} \phi_1 (r, z) \ \tilde{N}(\d r, \d z) \right] (\omega) + \left[  \int_t^T \int_{\R \backslash \{ 0 \}} \phi_1 (r) \ \tilde{N}(\d r, \d z) \right] (\omega),   \right. \\
            & \left. \cdots, x_n + \left[  \int_S^s \int_{\R \backslash \{ 0 \}} \phi_n (r, z) \ \tilde{N}(\d r, \d z) \right] (\omega) + \left[ \int_t^T \int_{\R \backslash \{ 0 \}} \phi_n (r, z) \ \tilde{N}(\d r, \d z) \right] (\omega) \right).
        \end{split}
    \end{equation*}

    In this case, the derivative of $F$ with respect to $\mathbf{F}_{[S, s] \cup [t, T]}$ is given for almost every $(\tau, x) \in (s, t) \times \R \backslash \{ 0 \}$ and $\P_{|_{\mathbf{F}_T}}$-almost surely by

    \begin{equation*}
        \begin{split}
            & \D_{(\tau, x)} \left( \mathbf{F}_{[S, s] \cup [t, T]}, N_{|_{[s, t]}} \right) F \\ 
            = & \quad g \left( \int_s^t \int_{\R \backslash \{ 0 \}} \phi_1 (r, z) \ \tilde{N}(\d r, \d z) + \phi_1(\tau, x), \cdots, \int_s^t \int_{\R \backslash \{ 0 \}} \phi_n (r, z) \ \tilde{N}(\d r, \d z) + \phi_n(\tau, x) \right) \\
            & - g \left( \int_s^t \int_{\R \backslash \{ 0 \}} \phi_1 (r, z) \ \tilde{N}(\d r, \d z), \cdots, \int_s^t \int_{\R \backslash \{ 0 \}} \phi_n (r, z) \ \tilde{N}(\d r, \d z)  \right).
        \end{split}
    \end{equation*}
    With the definition of $g$, it is in fact immediate that we have \eqref{eq:PetiteGrandeDerivee} when the operators are applied for variables in $\mathcal{S} \left( \mathbf{F}_{S}, N \right)$. Moreover, we get the following continuity inequality of the graph norms: for every $F \in \mathcal{S} \left( \mathbf{F}_{S}, N \right)$:
    \[ \norm{F}{\mathbb{D}^{1,2}_{\mathrm{J}} \left( \mathbf{F}_{[S, s] \cup [t, T]}, N_{|_{[s, t]}} \right)} \leqslant \norm{F}{\mathbb{D}^{1,2}_{\mathrm{J}} \left( \mathbf{F}_{S}, N \right)}.  \]
    This means that if we take the closure of both sets in the inclusion \eqref{EqDemPetiteDerivee} for the graph norm of $\D(\mathbf{F}_S, N)$, we conclude in \eqref{eq:inclusionDomainePetiteDerivee} and \eqref{eq:PetiteGrandeDerivee} in $\mathbb{D}^{1,2}_{\mathrm{J}} (\mathbf{F}_S, N)$. \end{Demo}

\subsection{Skorohod integral}

\begin{Definition} \label{def:Skorohod_integral}
    We define the \emph{domain of the operator $\delta$} as the following subset of $L^2(\Omega \times E)$:
    \begin{equation*}
        \begin{split}
            & \mathrm{Dom}_{\mathrm{J}}(\delta(\mathbf{F}_S, N)) \\
           \defeq  & \left\{ u \in L^2 \big( (\Omega, \mathbf{F}_T) \times E \big), \exists C > 0, \forall F \in \mathbb{D}^{1,2}_{\mathrm{J}} (\mathbf{F}_S, N), \E \left[ \langle \D F, u \rangle_{L^2 (E)} \right] \leqslant C \E[F^2]^{1/2} \right\}.
        \end{split}
    \end{equation*}

    For $u \in \mathrm{Dom}_{\mathrm{J}} (\delta)$, we define its \emph{Skorohod integral} $\delta(\mathbf{F}_S, N)u$ by being the unique (almost surely) real random variable such that for every $F \in \mathbb{D}^{1,2}_{\mathrm{J}} (\mathbf{F}_S, N)$, 
    \begin{equation}\label{eq:RelationDualité}
        \E \left[ \langle u, \D F \rangle_{L^2 (E)} \right] = \E[F \delta(\mathbf{F}_S, N)u].
    \end{equation}
    We will consider the following notation:
    \[ \delta(\mathbf{F}_S, N) u \defeq \int_S^T \int_{\R \backslash \{ 0 \}}  u(t, z) \ \tilde{N}^{\mathbf{F}_S} (\delta t, \delta z). \]
\end{Definition}

\begin{Remarque}
    Every occurrence of $\mathbb{D}^{1,2}_{\mathrm{J}} (\mathbf{F}_S, N)$ in this Definition can be replaced by $\mathcal{S} (\mathbf{F}_S, N)$.
\end{Remarque}

\begin{Remarque}\label{rq:skorohod_without_minus}
Let $u \in L^2 (\Omega \times E)$ be a measurable function such that $r \mapsto u(r,z)$ is càdlàg, and assume in addition that $u \in L^2( \Omega \times E)$. Then the process $u' \in L^2(\Omega \times E)$ defined by $u'(r,z) = u(r^-,z)$ belongs to the same equivalence class than $u$ in $L^2(\Omega,E)$, since they only differ from a the countable set of discontinuity of $r \mapsto u(r,z)$, which has Lebesgue measure zero. In particular, if $u \in \mathrm{Dom}_{\mathrm{J}} \big( \delta(\mathbf{F}_{S}, N) \big) $ we have that  $\P$-a.s:
\begin{equation*}
    \int_E u(r^{-},z) \ \tilde{N}^{\mathbf{F}_S} (\delta r, \delta z) = \int_E u(r,z) \ \tilde{N}^{\mathbf{F}_S} (\delta r, \delta z),
\end{equation*}
which can easily be seen by using the duality formula. Indeed, for every $F \in \mathbb{D}^{1, 2}_{\mathrm{J}} (\mathbf{F}_S, N)$ we have
\begin{align*}
        \E \left[ F \int_E u(r^{-},z) \tilde{N}^{\mathbf{F}_S} (\delta r, \delta z)   \right]  & = \E \left[ \int_E \D_{(r,z)} F \ u(r^-,z) \ \nu(\d z) \ \d r \right]
        \\ & = \E \left[ \int_E \D_{(r,z)} F \ u(r,z) \ \nu( \d z) \ \d r \right] \\ 
        & = \E \left[ F \int_E u(r,z) \ \tilde{N}^{\mathbf{F}_S} (\delta r, \delta z)   \right].
    \end{align*}
\end{Remarque}

\begin{Lemma}\label{LemmeDeriveeNulle}
    Let $F \in \mathbb{D}^{1,2}_{\mathrm{J}} (\mathbf{F}_S, N)$, and $S \leqslant s \leqslant t \leqslant T$. Then, $\E[F|\mathbf{F}_{[S, s] \cup [t, T]}] \in \mathbb{D}^{1,2}_{\mathrm{J}} (\mathbf{F}_S, N)$ and for almost every $(\tau, \xi) \in E$,
    \begin{equation}\label{eq:lemmeDeriveeNulle}
        \D_{(\tau, \xi)} \E \left. \big[ F \right| \mathbf{F}_{[S, s] \cup [t, T]} \big] = \E \left[ \left. \D_{(\tau, \xi)} F \right| \mathbf{F}_{[S, s] \cup [t, T]} \right] \mathbf{1}_{[S, s] \cup [t, T]} (\tau).
    \end{equation}
    In particular, if $F \in \mathbb{D}^{1,2}_{\mathrm{J}} (\mathbf{F}_S, N)$ is $\mathbf{F}_{[S, s] \cup [t, T]}$-measurable, then $\D_{(\tau, \xi)} F = 0$ for almost every $\xi \in \R \backslash \{ 0 \}$ and $\tau \in [s, t]$. 
\end{Lemma}

\begin{Demo} By density of step functions, we only have to deal with variables being random linear combination of variables of the type $\prod_{k=1}^n e^{\mathrm{i} \gamma_k \tilde{N}([t_k, t_{k+1}), A_k)}$, with $\mathbf{F}_S$-measurable coefficients in the linear combination, $\gamma_1, \cdots, \gamma_n$ being $\mathbf{F}_S$-measurable real random variables, $S \leqslant t_1 < t_2 < \cdots < t_{n+1} = T$ and $A_1, \cdots, A_n \in \mathcal{B}(\R \backslash \{ 0 \})$. By simplicity of notations, we will denote in the following of this proof $\phi_k \defeq \mathbf{1}_{[t_k, t_{k+1})} \otimes \mathbf{1}_{A_k}$. We consider 
\[ F:= \prod_{k=0}^n e^{\mathrm{i} \gamma_k \tilde{N}\big( [t_k, t_{k+1}), A_k \big)}.  \]

Let $1 \leqslant j \leqslant n+1$ such that $t_j \leqslant s < t_{j+1}$ and $j \leqslant J \leqslant n+1$ such that $t_J \leqslant t < t_{J+1}$. Then
\[ \begin{split}
      \E \left[ F | \mathbf{F}_{[S, s] \cup [t, T]} \right] = & \ \left( \prod_{k=1}^{j-1} e^{\mathrm{i} \gamma_k \tilde{N}(\phi_k)} \right) e^{\mathrm{i} \gamma_{j} \tilde{N}\big( [t_j, s), A_j \big)} \E \left[ e^{\mathrm{i} \gamma_j \tilde{N} \big( [s, t_{j+1}), A_j \big)} \left( \prod_{k=j+1}^{J-1} e^{\mathrm{i} \gamma_k \tilde{N}(\phi_k)} \right) e^{\mathrm{i} \gamma_J \tilde{N} \big( [t_J, t), A_J \big)} \right] \\
      & \quad e^{\mathrm{i} \gamma_J \tilde{N} \big( [t, t_{J+1}), A_J \big)} \left( \prod_{k=J+1}^{n} e^{\mathrm{i} \gamma_k \tilde{N}(\phi_k)} \right).
\end{split} \]
We compute the derivative of that variable (which belongs to $\mathcal{S} (\mathbf{F}_S, N)$)
\begin{equation*}
    \begin{split}
        & \D \E \left[ F | \mathbf{F}_{[S, s] \cup [t, T]} \right] \\
        = & \ \left( \prod_{k=1}^{j-1} e^{\mathrm{i} \gamma_k \tilde{N}(\phi_k)} \right) e^{\mathrm{i} \gamma_{j} \tilde{N} \big( [t_j, s), A_j \big)} \left( \left( \prod_{k=1}^{j-1} e^{\mathrm{i} \gamma_k \phi_k} \right) e^{\mathrm{i} \gamma_{j} \mathbf{1}_{[t_j, s)} \otimes \mathbf{1}_{A_j} } - 1 \right) \\
        & \E \left[ e^{\mathrm{i} \gamma_j \tilde{N} \big( [s, t_{j+1}), A_j \big)} \left( \prod_{k=j+1}^{J-1} e^{\mathrm{i} \gamma_k \tilde{N}(\phi_k)} \right) e^{\mathrm{i} \gamma_J \tilde{N} \big( [t_J, t), A_J \big)} \right] \\
      & \quad e^{\mathrm{i} \gamma_J \tilde{N} \big( [t, t_{J+1}), A_J \big)} \left( \prod_{k=J+1}^{n} e^{\mathrm{i} \gamma_k \tilde{N}(\phi_k)} \right) \left( \left( \prod_{k=J+1}^{n} e^{\mathrm{i} \gamma_k \phi_k} \right) e^{\mathrm{i} \gamma_{J} \mathbf{1}_{[t, t_{J+1})} \otimes \mathbf{1}_{A_J} }  - 1 \right) \\
     = & \ \E \big[ F \big| \mathbf{F}_{[S, s] \cup [t, T]} \big] \left( \left( \prod_{k=1}^{j-1} e^{\mathrm{i} \gamma_k \phi_k} \right) e^{\mathrm{i} \gamma_{j} \mathbf{1}_{[t_j, s)} \otimes \mathbf{1}_{A_j} } - 1 \right)  \left( \left( \prod_{k=J+1}^{n} e^{\mathrm{i} \gamma_k \phi_k} \right) e^{\mathrm{i} \gamma_{J} \mathbf{1}_{[t, t_{J+1})} \otimes \mathbf{1}_{A_J} }  - 1 \right).
    \end{split}
\end{equation*}
On the other hand, the derivative of $F$ is given by
\[ \D F = \prod_{k=0}^n e^{\mathrm{i} \gamma_k \tilde{N}\big( [t_k, t_{k+1}), A_k \big)} \left( \prod_{k=0}^n e^{\mathrm{i} \gamma_k \mathbf{1}_{[t_k, t_{k+1})} \otimes \mathbf{1}_{A_k}} -1 \right) = F \left( \prod_{k=0}^n e^{\mathrm{i} \gamma_k \mathbf{1}_{[t_k, t_{k+1})} \otimes \mathbf{1}_{A_k}} -1 \right). \]
Hence, since the $\gamma_j$'s are $\mathbf{F}_S$-measurable, so $\mathbf{F}_{[S, s] \cup [t, T]}$-measurable, we have 
\[ \E \big[ \D F \big| \mathbf{F}_{[S, s] \cup [t, T]} \big] = \E\big[ F \big| \mathbf{F}_{[S, s] \cup [t, T]} \big] \left( \prod_{k=0}^n e^{\mathrm{i} \gamma_k \mathbf{1}_{[t_k, t_{k+1})} \otimes \mathbf{1}_{A_k}} -1 \right). \]
This proves equality \eqref{eq:lemmeDeriveeNulle} for $F$ being a monomial.

It still holds when $F$ is linear random combination of those kind of monomials with $\mathbf{F}_S$-measurable coefficients. Hence, we conclude that if $F \in \mathcal{S} (\mathbf{F}_S, N)$, then $\E[F|\mathbf{F}_{[S, s] \cup [t, T]}] \in \mathbb{D}^{1,2}_{\mathrm{J}} (\mathbf{F}_S, N)$. For $F \in \mathbb{D}^{1,2}_{\mathrm{J}} (\mathbf{F}_S, N)$, $F$ is a limit for the norm graph $\norm{\cdot}{\mathbb{D}^{1,2}_{\mathrm{J}} (\mathbf{F}_S, W)}$ introduced in Proposition \ref{PropFermabilitéDérivée} of a sequence $(F_n)_{n \geqslant 1} \subset \mathcal{S} (\mathbf{F}_S, N)$. Hence, we have by \eqref{eq:lemmeDeriveeNulle} 
\[ \D \E \big[ F_n \left| \ \mathbf{F}_{[S, s] \cup [t, T]} \right. \big] = \E \big[ \D F_n \left| \ \mathbf{F}_{[S, s] \cup [t, T]} \right. \big] \mathbf{1}_{[S, r] \cup [t, T]} \xrightarrow[n \to \infty]{L^2 (\Omega \times E)} \E \big[ \D F \left| \ \mathbf{F}_{[S, s] \cup [t, T]} \right. \big] \mathbf{1}_{[S, s] \cup [t, T]}.  \]
Hence, the sequence $(\E [ F_n | \mathbf{F}_{[S, s] \cup [t, T]} ], n \geqslant 1)$ converges in $L^2 (\mathbf{F}_T)$ toward $(\E [ F | \mathbf{F}_{[S, s] \cup [t, T]}], n \geqslant 1)$ and the sequence $(\D \E [ F_n | \mathbf{F}_{[S, s] \cup [t, T]} ], n \geqslant 1)$ converges in $L^2 (\Omega \times E)$ by the previous computation. By the definition of the closability of $\D$, this means that $\E[F|\mathbf{F}_{[S, s] \cup [t, T]}] \in \mathbb{D}^{1,2}_{\mathrm{J}} (\mathbf{F}_S, N)$ and that the equality \eqref{eq:lemmeDeriveeNulle} is true for $F$.      
\end{Demo}

\begin{Lemma}\label{LemmeConvergenceDelta}
    Let $u \in L^2(\Omega \times E)$. We suppose that there exists a sequence of processes $(u_n)_{n \geqslant 1}$ in $L^2 (\Omega \times E)$ such that
    \begin{enumerate}[label=\Alph*.]
        \item For every $n \geqslant 1$, $u_n \in \mathrm{Dom}_{\mathrm{J}} \big( \delta(\mathbf{F}_{S}, N) \big)$; 
        \item The sequence $(u_n)_{n \geqslant 1}$ weakly converges to $u$ in $L^2 (\Omega \times E)$; 
        \item We have 
        \[ \sup_{n \geqslant 1} \E \left[ \left( \int_E u_n (r, z) \ \tilde{N}^{\mathbf{F}_S} (\delta r, \delta z) \right)^2 \right] < \infty. \]
    \end{enumerate}
    Then $u \in \mathrm{Dom}_{\mathrm{J}} \big( \delta(\mathbf{F}_S, N) \big)$ and $\int_E u_n (r, z) \ \tilde{N}^{\mathbf{F}_S} (\delta r, \delta z)$ converges weakly in $L^2 (\Omega)$ as $n$ goes to $\infty$ toward $\int_E u (r, z) \ \tilde{N}^{\mathbf{F}_S} (\delta r, \delta z)$.
\end{Lemma}

\begin{Demo}
    The following duality formula holds for every $n \geqslant 1$ and $F \in \mathbb{D}^{1, 2}_{\mathrm{J}} (\mathbf{F}_S, N)$
    \[ \E \left[ F \int_{E} u_n (r, z) \ \tilde{N}^{\mathbf{F}_S} (\delta r, \delta z) \right] = \E \left[ \int_{E} \D_{(r, z)} F \ u_n (r, z) \ \nu(\d z) \ \d r \right]. \]
    In particular, we have 
    \[ \left| \E \left[ \int_{E} \D_{(r, z)} F \ u_n (r, z) \ \nu(\d z) \ \d r  \right] \right| \leqslant \sup_{n \geqslant 1} \E \left[ \left( \int_E u_n (r, z) \ \tilde{N}^{\mathbf{F}_S} (\delta r, \delta z) \right)^2 \right]^{1/2} \E[F^2]^{1/2}. \]
    By $L^2$ weak convergence of $u_n$ toward $u$, we deduce by letting $n$ going to $\infty$ in the previous inequality that $u$ belongs to $\mathrm{Dom}_{\mathrm{J}} \big( \delta(\mathbf{F}_S, N) \big)$. Moreover, by weak convergence and duality formulas, we have in fact
    \[ \E \left[ F \int_{E} u_n (r, z) \ \tilde{N}^{\mathbf{F}_S} (\delta r, \delta z) \right] \xrightarrow[n \to \infty]{} \E \left[ F \int_{E} u (r, z) \ \tilde{N}^{\mathbf{F}_S} (\delta r, \delta z) \right], \]
    which holds for every $F \in \mathbb{D}^{1, 2}_{\mathrm{J}} (\mathbf{F}_S, N)$. It also holds for $F \in L^2(\mathbf{F}_T)$ by density in $L^2$ of $\mathbb{D}^{1, 2}_{\mathrm{J}} (\mathbf{F}_S, N)$. This proves the weak convergence in $L^2$ of the Skorohod integrals.
\end{Demo}

\begin{Theoreme}\label{thmMalliavin:Sko=Ito}
    The following inclusion holds: $\mathcal{P}_2^{\mathbf{F}_{[S, \cdot] \cup [t, T]}} (s, t) \subset \mathrm{Dom}_{\mathrm{J}} \big( \delta(\mathbf{F}_{[S, s] \cup [t, T]}, N_{|_{[s, t]}}) \big)$ (recall the first set is defined on \eqref{eqIto:defP2}). Moreover, for every $u \in \mathcal{P}_2^{\mathbf{F}_{[S, \cdot] \cup [t, T]}} (s, t)$, we have
    \[ \int_{s}^t \int_{\R \backslash \{ 0 \}} u(r, z) \ \tilde{N}^{\mathbf{F}_{[S, \cdot] \cup [t, T]}} (\d r, \d z) = \int_s^t \int_{\R \backslash \{ 0 \}} u(r, z) \ \tilde{N}^{\mathbf{F}_{[S, s] \cup [t, T]}} (\delta r, \delta z). \]
\end{Theoreme}

\begin{Demo}
Every process $u \in \mathcal{P}_2^{\mathbf{F}_{[S, \cdot] \cup [t, T]}} (s, t)$ is limit in $L^2 (\Omega \times E)$ of elementary processes which looks like $\sum_{i=1}^n F_i \mathbf{1}_{[s_i, t_i] \times A}$, where $s \leqslant s_1 < t_1 < s_2 < \cdots < t_n \leqslant t$, $A \in \mathcal{B}(\R \backslash \{ 0 \})$ and the $F_i$'s are $\mathbf{F}_{[S,s_i] \cup [t, T]}$-adapted. 

For $s \leqslant r < v \leqslant t$, let $F \in \mathbb{D}^{1,2}_{\mathrm{J}} (\mathbf{F}_{[S, s] \cup [t, T]}, N_{|_{[s, t]}})$ being $\mathbf{F}_{[S, r] \cup [t, T]}$-measurable. Then, for $A \in \mathcal{B}(\R \backslash \{ 0 \})$, we show that $F \mathbf{1}_{[r, v] \times A} \in \mathrm{Dom}_{\mathrm{J}} \big( \delta(\mathbf{F}_{[S, s] \cup [t, T]}, N_{|_{[s, t]}}) \big)$ and that the Skorohod integral coincides with the Itô's one. By the Integration by Parts formula (Proposition \ref{PropIPP}), for every $G \in \mathbb{D}^{1,2}_{\mathrm{J}} (\mathbf{F}_{[S, s] \cup [t, T]}, N_{|_{[s, t]}})$, we have
\begin{equation*}
    \begin{split}
        &  \E \left[ \left\langle \D G, F \mathbf{1}_{[r, v] \times A} \right\rangle_{L^2(E)} \right] \\
        = \ & \E \left[ G F \tilde{N}([r, v], A) \right] - \E \left[G \left\langle \D F,  \mathbf{1}_{[r, v] \times A} \right\rangle_{L^2(E)} \right] - \E \left[ \left\langle \D F \ \D G, \mathbf{1}_{[r, v] \times A} \right\rangle_{L^2(E)} \right].
    \end{split}
\end{equation*}
Since $F$ is $\mathbf{F}_{[S, r] \cup [t, T]}$-measurable, we conclude by Lemma \ref{LemmeDeriveeNulle} that $\D F = 0$ on $[r, t] \times \R \backslash \{ 0 \}$, so that
\[ \E \left[ \left\langle \D G, F \mathbf{1}_{[r, v] \times A} \right\rangle_{L^2(E)} \right] =  \E \left[ G F \tilde{N}([r, v], A) \right]. \]
This proves that $F \mathbf{1}_{[r, v] \times A}$ belongs to $\mathrm{Dom}_{\mathrm{J}} \big( \delta(\mathbf{F}_{[S, s] \cup [t, T]}, N_{|_{[s, t]}}) \big)$ and its Skorohod integral coincides with the Itô integral by unicity in the duality relation. 

By density of $\mathbb{D}^{1,2}_{\mathrm{J}} (\mathbf{F}_{[S, s] \cup [t, T]}, N_{|_{[s, t]}})$, we extend the previous result for $F \in L^2(\mathbf{F}_T)$ being $\mathbf{F}_{[S, r] \cup [t, T]}$-measurable. Indeed, since $\mathbb{D}^{1,2}_{\mathrm{J}} (\mathbf{F}_{[S, s] \cup [t, T]}, N_{|_{[s, t]}})$ is dense into $L^2 (\mathbf{F}_T)$, we can approximate $F \in L^2(\mathbf{F}_{[S, r] \cup [t, T]})$ (for the $L^2$-norm) by a sequence $(F_n)_n$ included in $\mathbb{D}^{1,2}_{\mathrm{J}}(\mathbf{F}_{[S, s] \cup [t, T]}, N_{|_{[s, t]}})$. We can even choose a sequence which is $\mathbf{F}_{[S, r] \cup [t, T]}$-measurable by considering the conditional expectation $\E[F_n | \mathbf{F}_{[S, r] \cup [t, T]}]$. Then, we can apply Lemma \ref{LemmeConvergenceDelta} to prove that $F \mathbf{1}_{[r, v] \times A}$ belongs to $\mathrm{Dom}_{\mathrm{J}} \big( \delta(\mathbf{F}_{[S, s] \cup [t, T]}, N_{|_{[s, t]}}) \big)$ and that once again, the Skorohod and Itô integrals coincides. By linearity, we proved that every elementary process belongs to the domain and that their Skorohod integrals coincide with the Itô's one. 

Now, for $u \in \mathcal{P}_2^{\mathbf{F}_{[S, \cdot] \cup [t, T]}} (s, t)$, it is once again the use of Lemma \ref{LemmeConvergenceDelta} that allows us to conclude, added with the fact that when we approximate $u$ by a sequence of elementary process, we approximate their Itô's integrals too, in $L^2$. 
\end{Demo}

\begin{Theoreme}\label{thmMalliavin:Chasles}
    Let $u: \Omega \times E \longrightarrow \R$ a stochastic process and $S \leqslant s < t \leqslant T$. The following are equivalent.

\begin{enumerate}
    \item The restricted process $u: \Omega \times (s, t) \times \R \backslash \{ 0 \} \longrightarrow \R$ belongs to $\mathrm{Dom}_{\mathrm{J}} \big( \delta(\mathbf{F}_{[S, s] \cup [t, T]}, N_{|_{[s, t]}}) \big)$; 

    \item The process $u \mathbf{1}_{[s, t]}$ belongs to $\mathrm{Dom}_{\mathrm{J}} \big( \delta(\mathbf{F}_S, N) \big)$.
\end{enumerate}
In this case, we have 
\[ \int_s^t \int_{\R \backslash \{ 0 \}} u(z, r) \ \tilde{N}^{\mathbf{F}_{[S, s] \cup [t, T]}} (\delta t, \delta z) = \int_S^T \int_{\R \backslash \{ 0 \}} u(z, r) \mathbf{1}_{[s, t]} (r) \ \tilde{N}^{\mathbf{F}_S} (\delta t, \delta z).  \]
\end{Theoreme}

The proof of this Theorem is exactly the same as \cite{ItoAG}, Proposition A.7.

\section{Proof of Theorem \ref{maintheo}} \label{sec:proof_maintheo}

\subsection{Time discretisation setting}

For $x \in \R$ and $0 \leqslant s \leqslant t \leqslant T$, we recall the definition of the flow process 
\begin{equation*}
    X_{s,t}^{x} = x + \int_s^t \mu(r,X_{s,r}^{x}) \ \d r + \int_s^t \int_{\R \backslash \{ 0 \}} \sigma ( r,X_{s,r^-}^{x},z ) \ \tilde{N}(\d r, \d z),
\end{equation*}
along with the stochastic process
\begin{equation*}
    Y_t = Y_0 + \int_0^t A_s \ \d s + \int_0^t \int_{\R \backslash \{ 0 \}} B_{s,z} \ \tilde{N}(\d s, \d z), \quad t \in [0,T].
\end{equation*}
Let $n \geqslant 1$ and $h = \frac{T}{n}$. We have $\P$-a.s., 
\begin{equation}\label{eq:decoupage}
    \begin{split}
        f \left( X_{0,T}^{Y_0} \right) - f(Y_T) & = \sum_{i=0}^{n-1} \left( f \left( X_{ih,T}^{Y_{ih}} \right) - f \left( X_{(i+1)h,T}^{Y_{(i+1)h}} \right) \right) \\ & = \sum_{i=0}^{n-1} \left(  f \left( X_{ih,T}^{Y_{ih}} \right) -  f\left(X_{(i+1)h,T}^{Y_{ih}} \right) \right) - \sum_{i=0}^{n-1} \left(  f \left( X_{(i+1)h,T}^{Y_{(i+1)h}} \right) -  f \left( X_{(i+1)h,T}^{Y_{ih}} \right) \right).
    \end{split}
\end{equation}

\subsection{From Itô to Skorohod Poisson random integrals}
\textbf{We start by analysing the second term} in the right-hand side of \eqref{eq:decoupage}, which is slightly simpler than the first one. 
We want apply the Itô formula for independent random fields of Poisson random integrals (Proposition \ref{prop:ito_formula_random_functions_poisson}). For this purpose we fix $i \in \{0, \ldots, n-1\}$ and consider the stochastic process $(Y_r)_{r \in [ih, (i+1)h]}$ which is $(\mathbf{F}_{[ih,t]})_{t \in [ih, (i+1)h]}$-adapted, where we recall that 
\begin{equation*}
    \mathbf{F}_{[ih,t]} = \sigma( N([ih,s],A), \ ih \leqslant s \leqslant t, \ A \in \mathcal{B}(\R \backslash \{0 \})),
\end{equation*}
along with the application $(\omega,x) \ni \Omega \times \R \mapsto f(X_{(i+1)h,T}^x(\omega))$. Define 
\begin{equation*}
    \mathbf{G}_{T}^{i,h} \defeq \sigma(N([(i+1)h,T],A), \ A \in \mathcal{B}(\R \backslash \{0 \})),
\end{equation*}
so that for fixed $h$ and $i$, $\mathbf{G}_{T}^{i,h}$ is independent from $\mathbf{F}_{[ih,t]}$ for every $t \in [ih, (i+1)h]$, and  $\omega \mapsto f(X_{(i+1)h,T}^x(\omega))$ is $\mathbf{G}_{T}^{i,h}$-measurable. Finally, we define
\begin{equation*}
    \mathbf{H}^{i,h}_t \defeq \mathbf{F}_{[ih,t]} \vee \mathbf{G}_{T}^{i,h}.
\end{equation*}
By Proposition \ref{prop:ito_formula_random_functions_poisson}, we obtain
\begin{align*}
    & f \left( X_{(i+1)h,T}^{Y_{(i+1)h}} \right) -  f \left( X_{(i+1)h,T}^{Y_{ih}} \right) = \int_{ih}^{(i+1)h} f' \left( X_{(i+1)h,T}^{Y_r} \right) \partial_x X_{(i+1)h,T}^{Y_r} A_r \ \d r \\
     &  + \int_{ih}^{(i+1)h} \int_{\R \backslash \{ 0 \}} \left( f \Big( X_{(i+1)h,T}^{Y_{r^-}+B_{r,z}} \Big) - f \Big( X_{(i+1)h,T}^{Y_{r^-}} \Big) \right) \tilde{N}^{\mathbf{H}^{i, h}}(\d r, \d z) \\ & + \int_{ih}^{(i+1)h} \int_{\R \backslash \{ 0 \}} \left( f \left( X_{(i+1)h,T}^{Y_{r}+B_{r,z}} \right) - f \left( X_{(i+1)h,T}^{Y_{r}} \right) - B_{r,z} f' \left( X_{(i+1)h,T}^{Y_r} \right) \partial_x X_{(i+1)h,T}^{Y_r} \right) \ \nu( \d z) \ \d r. 
\end{align*}
We now want to transform the Itô integral into a Skorohod integral, which requires to apply Theorem \ref{thmMalliavin:Sko=Ito}. Toward this goal we prove in section \ref{proof:controle_l2_second_terme} the following Lemma.

\begin{Lemma} \label{lm:controle_l2_second_terme}
There exists a positive function $\phi_1 \in L^1(\nu)$ such that for $\nu$-almost every $z \in \R \backslash \{ 0 \}$,
\begin{equation*}
   \sup_{h \in (0,T)} \sup_{0\leqslant i \leqslant \lfloor T/h \rfloor - 1} \sup_{r \in [ih, (i+1)h)} \E \left[  \left( f \Big( X_{(i+1)h,T}^{Y_{r}+B_{r,z}} \Big) - f \Big( X_{(i+1)h,T}^{Y_{r}} \Big) \right)^2 \right] \leq \phi_1(z).
\end{equation*}
\end{Lemma}

In particular, we deduce from Lemma \ref{lm:controle_l2_second_terme}  that 
\begin{equation}\label{eq:MomentL2DeuxiemeTerme}
   \E \left[ \int_{ih}^{(i+1)h} \int_{\R \backslash \{ 0 \}} \left( f \Big( X_{(i+1)h,T}^{Y_{r}+B_{r,z}} \Big) - f \Big( X_{(i+1)h,T}^{Y_{r}} \Big) \right)^2 \ \nu(\d z) \ \d r  \right] \leqslant T \|\phi_1 \|_{L^1(\nu)} < + \infty.
\end{equation}
Then due to \eqref{eq:MomentL2DeuxiemeTerme} and Theorem \ref{thmMalliavin:Sko=Ito}, the Itô and Skorohod integrals coincide:
\[ \begin{split}
     & \int_{ih}^{(i+1)h} \int_{\R \backslash \{ 0 \}} \left( f \Big( X_{(i+1)h,T}^{Y_{r^-}+B_{r,z}} \Big) - f \Big( X_{(i+1)h,T}^{Y_{r^-}} \Big) \right) \tilde{N}^{\mathbf{H}^{i, h}}(\d r, \d z) \\
     = & \ \int_{ih}^{(i+1)h} \int_{\R \backslash \{ 0 \}} \left( f \Big( X_{(i+1)h,T}^{Y_{r}+B_{r,z}} \Big) - f \Big( X_{(i+1)h,T}^{Y_{r}} \Big) \right) \tilde{N}^{\mathbf{F}_{[0, ih] \cup [(i+1)h, T]}}(\delta r, \delta z) .
\end{split}  \]
Note that we may replace $Y_{r^-}$ by $Y_r$ in the Skorohod integral due to Remark \ref{rq:skorohod_without_minus}. 
By using Theorem \ref{thmMalliavin:Chasles} (the Chasles relation) afterwards , we have 
\[ \begin{split}
     & \int_{ih}^{(i+1)h} \int_{\R \backslash \{ 0 \}} \left( f \Big( X_{(i+1)h,T}^{Y_{r}+B_{r,z}} \Big) - f \Big( X_{(i+1)h,T}^{Y_{r}} \Big) \right) \tilde{N}^{\mathbf{H}^{i,h}}(\d r, \d z) \\
     = & \ \int_{0}^{T}  \int_{\R \backslash \{ 0 \}} \mathbf{1}_{[ih, (i+1)h]} (r) \left( f \Big( X_{(i+1)h,T}^{Y_{r}+B_{r,z}} \Big) - f \Big( X_{(i+1)h,T}^{Y_{r}} \Big) \right) \tilde{N}^{\mathbf{F}_{0}}(\delta r, \delta z),
\end{split}  \]
where we recall that $\mathbf{F}_0$ is the $\sigma$-algebra generated by $N(\{ 0 \},A)$, for $A \in \mathcal{B}(\R \backslash \{ 0 \})$. In other words, $\mathbf{F}_0$ is exactly the $\sigma$-algebra generated by $\mathcal{N}$, the negligible sets, so that 
\[ \begin{split}
     & \int_{ih}^{(i+1)h} \int_{\R \backslash \{ 0 \}} \left( f \Big( X_{(i+1)h,T}^{Y_{r^-}+B_{r,z}} \Big) - f \Big( X_{(i+1)h,T}^{Y_{r^-}} \Big) \right) \tilde{N}^{\mathbf{H}^{i,h}}(\d r, \d z) \\
     = & \ \int_{0}^{T}  \int_{\R \backslash \{ 0 \}} \mathbf{1}_{[ih, (i+1)h]} (r) \left( f \Big( X_{(i+1)h,T}^{Y_{r}+B_{r,z}} \Big) - f \Big( X_{(i+1)h,T}^{Y_{r}} \Big) \right) \tilde{N}(\delta r, \delta z) \\
     = & \int_{0}^{T}  \int_{\R \backslash \{ 0 \}} \mathbf{1}_{[ih, (i+1)h]} (r) \left( f \Big( X_{\lceil r \rceil_h,T}^{Y_{r}+B_{r,z}} \Big) - f \Big( X_{\lceil r \rceil_h,T}^{Y_{r}} \Big) \right) \tilde{N}(\delta r, \delta z),
\end{split}  \]
where the latter integral is the classical Skorohod integral and $\lceil r \rceil_h$ denotes the lowest $ih$, with $i \geqslant 1$ such that $r \leqslant ih$. We finally conclude by summing this equalities for $0 \leqslant i \leqslant n-1$ that
\[ \begin{split}
     & \sum_{i=0}^{n-1} \int_{ih}^{(i+1)h} \int_{\R \backslash \{ 0 \}} \left( f \Big( X_{(i+1)h,T}^{Y_{r^-}+B_{r,z}} \Big) - f \Big( X_{(i+1)h,T}^{Y_{r^-}} \Big) \right) \tilde{N}^{\mathbf{H}^{i,h}}(\d r, \d z) \\
     = & \ \int_{0}^{T}  \int_{\R \backslash \{ 0 \}}  \left( f \Big( X_{\lceil r \rceil_h,T}^{Y_{r}+B_{r,z}} \Big) - f \Big( X_{\lceil r \rceil_h,T}^{Y_{r}} \Big) \right) \tilde{N}(\delta r, \delta z).
\end{split}  \]
Consequently, the second term of \eqref{eq:decoupage} finally writes
\begin{equation}\label{eq:Decoupage2emeterme}
    \begin{split}
        & \sum_{i=0}^{n-1}  f \left( X_{(i+1)h,T}^{Y_{(i+1)h}} \right) -  f \left( X_{(i+1)h,T}^{Y_{ih}} \right) = \int_{0}^{T} f' \left( X_{\lceil r \rceil_h,T}^{Y_r} \right) \partial_x X_{\lceil r \rceil_h,T}^{Y_r} A_r \ \d r \\
     + &  \int_{0}^{T}  \int_{\R \backslash \{ 0 \}}  \left( f \Big( X_{\lceil r \rceil_h,T}^{Y_{r}+B_{r,z}} \Big) - f \Big( X_{\lceil r \rceil_h,T}^{Y_{r}} \Big) \right) \tilde{N}(\delta r, \delta z) \\ 
     + & \int_{0}^{T} \int_{\R \backslash \{ 0 \}} \left( f \left( X_{\lceil r \rceil_h,T}^{Y_{r}+B_{r,z}} \right) - f \left( X_{\lceil r \rceil_h,T}^{Y_{r}} \right) - B_{r,z} f' \left( X_{\lceil r \rceil_h,T}^{Y_r} \right) \partial_x X_{\lceil r \rceil_h,T}^{Y_r} \right) \ \nu( \d z) \ \d r 
    \end{split}
\end{equation}

\textbf{We turn to the analysis of the first term in} the right hand side of \eqref{eq:decoupage}. For the latter, the issue is that the flow process $X^{Y_{ih}}_{\cdot,T}$ is not taken at the same time in the two terms $f \left( X_{ih,T}^{Y_{ih}} \right) $ and $ f \left( X_{(i+1)h,T}^{Y_{ih}} \right)$. We follow \cite{ItoAG} and push the flow. For $n \in \N$, $h = T/n$, and $i \in \{0, \ldots, n-1 \}$, using \ref{hyp:PushTheFlow} and the continuity of $x \mapsto X_{ih,T}^x$ and $x \mapsto X_{(i+1)h,T}^x$, we have that $\P$-almost surely,
\begin{equation*}
    X_{ih,T}^{Y_{ih}} = X_{ih,T}^{X_{ih,(i+1)h}^{Y_{ih}}}
\end{equation*}
and 
\begin{equation*}
    X_{(i+1)h,T}^{Y_{ih}} = X_{(i+1)h,T}^{X_{ih,ih}^{Y_{ih}}}.
\end{equation*}
It follows that $\P$-almost surely,
\begin{equation*}
     f \left( X_{ih,T}^{Y_{ih}} \right) -  f \left( X_{(i+1)h,T}^{Y_{ih}} \right) = f \left( X_{(i+1)h,T}^{X^{Y_{ih}}_{ih,(i+1)h}} \right) -  f \left( X_{(i+1)h,T}^{X^{Y_{ih}}_{ih,ih}} \right).
\end{equation*}
which allows us to apply the Itô formula for independent random fields (Proposition \ref{prop:ito_formula_random_functions_poisson}) to the map $(\omega, x) \mapsto f(X_{(i+1)h,T}^x(\omega))$ satisfying that $\omega \mapsto f(X_{(i+1)h,T}^x(\omega))$ is $\mathbf{G}_{T}^{i,h}$-measurable for every $x \in \R$, and to the $(\mathbf{F}_{[ih,t]})_{t \in [ih, (i+1)h]}$-adapted stochastic process $(X^{Y_{ih}}_{ih,r}, r \in [ih, (i+1)h])$ between times $r=ih$ and $r=(i+1)h$, leading to:
\[ \begin{split}
      & f \left( X_{ih,T}^{Y_{ih}} \right) -  f \left( X_{(i+1)h,T}^{Y_{ih}} \right) = \int_{ih}^{(i+1)h} f' \left( X_{(i+1)h,T}^{X^{Y_{ih}}_{ih,r}} \right) \partial_x X_{(i+1)h,T}^{X^{Y_{ih}}_{ih,r}} \ \mu \left( r, X^{Y_{ih}}_{ih,r} \right) \ \d r \\
      + & \int_{ih}^{(i+1)h} \int_{\R \backslash \{ 0 \}} \left\{ f \left( X_{(i+1)h,T}^{X^{Y_{ih}}_{ih,{r^-}} + \sigma \left( r,X^{Y_{ih}}_{ih,{r^-}},z \right)} \right) - f \left( X_{(i+1)h,T}^{X^{Y_{ih}}_{ih,{r^-}}} \right) \right\} \ \tilde{N}^{\mathbf{H}^{i, h}}(\d r, \d z) \\ 
      + & \int_{ih}^{(i+1)h} \int_{\R \backslash \{ 0 \}} \left\{ f \left( X_{(i+1)h,T}^{X^{Y_{ih}}_{ih,{r}} + \sigma \left( r,X^{Y_{ih}}_{ih,{r}},z \right)} \right) - f \left( X_{(i+1)h,T}^{X^{Y_{ih}}_{ih,{r}}} \right) \right. \\
      & \left. - \sigma \left( r,X^{Y_{ih}}_{ih,{r}},z \right) f'\left( X_{(i+1)h,T}^{X^{Y_{ih}}_{ih,r}} \right) \partial_x X_{(i+1)h,T}^{X^{Y_{ih}}_{ih,r}} \right\} \ \nu( \d z) \d r.
\end{split} \]
As for the second term, we need the finiteness of the $L^2$-norm of the stochastic integral. This is provided by the following Lemma, proved in Section \ref{proof:controle_l2_premier_terme}. 

\begin{Lemma} \label{lm:controle_l2_premier_terme}
There exists a positive function $\phi_2 \in L^1(\nu)$ such that for almost every $r \in [0,T]$ and $z \in \R \backslash \{ 0 \}$,
\begin{equation*}
    \sup_{h \in (0,T)} \sup_{0 \leqslant i \leqslant \lfloor T/h \rfloor-1} \sup_{r \in [ih, (i+1)h)} \E \left[ \left( f \Big( X_{(i+1)h,T}^{X_{ih,r}^{Y_{ih}}+\sigma(r,X_{ih,r}^{Y_{ih}},z)} \Big) -  f \Big( X_{(i+1)h,T}^{X_{ih,r}^{Y_{ih}}} \Big) \right)^2 \right] \leqslant \phi_2(z).
\end{equation*}
\end{Lemma}
In particular, due to Lemma \ref{lm:controle_l2_premier_terme} we get
 \begin{equation}
        \E \left[ \int_{ih}^{(i+1)h} \int_{\R \backslash \{0 \}} \left\{ f \Big( X_{(i+1)h,T}^{X_{ih,r}^{Y_{ih}}+\sigma(r,X_{ih,r}^{Y_{ih}},z)} \Big) -  f \Big( X_{(i+1)h,T}^{X_{ih,r}^{Y_{ih}}} \Big) \right\}^2 \nu(\d z) \d r \right] \leqslant T \| \phi_2 \|_{L^1(\nu)} < + \infty.
    \end{equation}
We can apply Theorem \ref{thmMalliavin:Sko=Ito}, giving for the stochastic integral:
\[ \begin{split}
    & \int_{ih}^{(i+1)h} \int_{\R \backslash \{ 0 \}} \left\{ f \left( X_{(i+1)h,T}^{X^{Y_{ih}}_{ih,{r^-}} + \sigma \left( r,X^{Y_{ih}}_{ih,{r^-}},z \right)} \right) - f \left( X_{(i+1)h,T}^{X^{Y_{ih}}_{ih,{r^-}}} \right) \right\} \ \tilde{N}^{\mathbf{H}^{i, h}}(\d r, \d z) \\
    = \ & \int_{ih}^{(i+1)h} \int_{\R \backslash \{ 0 \}} \left\{ f \left( X_{(i+1)h,T}^{X^{Y_{ih}}_{ih,{r}} + \sigma \left( r,X^{Y_{ih}}_{ih,{r}},z \right)} \right) - f \left( X_{(i+1)h,T}^{X^{Y_{ih}}_{ih,{r}}} \right) \right\} \ \tilde{N}^{\mathbf{F}_{[0, ih] \cup [(i+1)h, T]}}(\delta r, \delta z).
\end{split} \]
By Theorem \ref{thmMalliavin:Chasles}, we have
\[ \begin{split}
    & \int_{ih}^{(i+1)h} \int_{\R \backslash \{ 0 \}} \left\{ f \left( X_{(i+1)h,T}^{X^{Y_{ih}}_{ih,{r^-}} + \sigma \left( r,X^{Y_{ih}}_{ih,{r^-}},z \right)} \right) - f \left( X_{(i+1)h,T}^{X^{Y_{ih}}_{ih,{r^-}}} \right) \right\} \ \tilde{N}^{\mathbf{H}^{i, h}}(\d r, \d z) \\
    = \ & \int_{0}^{T} \int_{\R \backslash \{ 0 \}} \mathbf{1}_{[ih, (i+1)h]} (r) \left\{ f \left( X_{(i+1)h,T}^{X^{Y_{ih}}_{ih,{r}} + \sigma \left( r,X^{Y_{ih}}_{ih,{r}},z \right)} \right) - f \left( X_{(i+1)h,T}^{X^{Y_{ih}}_{ih,{r}}} \right) \right\} \ \tilde{N}^{\mathbf{F}_{0}}(\delta r, \delta z) \\
    = \ & \int_{0}^{T} \int_{\R \backslash \{ 0 \}} \mathbf{1}_{[ih, (i+1)h]} (r) \left\{ f \left( X_{\lceil r \rceil_h,T}^{X^{Y_{\lfloor r \rfloor_h}}_{\lfloor r \rfloor_h,{r}} + \sigma \left( r,X^{Y_{\lfloor r \rfloor_h}}_{\lfloor r \rfloor_h,{r}},z \right)} \right) - f \left( X_{\lceil r \rceil_h,T}^{X^{Y_{\lfloor r \rfloor_h}}_{\lfloor r \rfloor_r,{r}}} \right) \right\} \ \tilde{N} (\delta r, \delta z),
\end{split} \]
where $\lfloor r \rfloor_h$ is the biggest $ih$, with $i \geqslant 0$ such that $ih \leqslant r$. By summing those equalities for $0 \leqslant i \leqslant n-1$, we get
\[ \begin{split}
    & \sum_{i=0}^{n-1} \int_{ih}^{(i+1)h} \int_{\R \backslash \{ 0 \}} \left\{ f \left( X_{(i+1)h,T}^{X^{Y_{ih}}_{ih,{r^-}} + \sigma \left( r,X^{Y_{ih}}_{ih,{r^-}},z \right)} \right) - f \left( X_{(i+1)h,T}^{X^{Y_{ih}}_{ih,{r^-}}} \right) \right\} \ \tilde{N}^{\mathbf{H}^{i, h}}(\d r, \d z) \\
    = \ & \int_0^T \int_{\R \backslash \{ 0 \}}  \left\{ f \left( X_{\lceil r \rceil_h,T}^{X^{Y_{\lfloor r \rfloor_h}}_{\lfloor r \rfloor_h,{r}} + \sigma \left( r,X^{Y_{\lfloor r \rfloor_h}}_{\lfloor r \rfloor_h,{r}},z \right)} \right) - f \left( X_{\lceil r \rceil_h,T}^{X^{Y_{\lfloor r \rfloor_h}}_{\lfloor r \rfloor_h,{r}}} \right) \right\} \ \tilde{N} (\delta r, \delta z).
\end{split} \]
Finally, the first term writes
\begin{equation}\label{eq:Decoupage1erTerme}
        \begin{split}
      & \sum_{i=0}^{n-1} \left( f \left( X_{ih,T}^{Y_{ih}} \right) -  f \left( X_{(i+1)h,T}^{Y_{ih}} \right) \right) = \int_{0}^{T} f' \left( X_{\lceil r \rceil_h,T}^{X^{Y_{\lfloor r \rfloor_h}}_{\lfloor r \rfloor_h,r}} \right) \partial_x X_{\lceil r \rceil_h,T}^{X^{Y_{\lfloor r \rfloor_h}}_{\lfloor r \rfloor_h,r}} \ \mu \left( r, X^{Y_{\lfloor r \rfloor_h}}_{\lfloor r \rfloor_h,T} \right) \ \d r \\
      + & \int_0^T \int_{\R \backslash \{ 0 \}}  \left\{ f \left( X_{\lceil r \rceil_h,T}^{X^{Y_{\lfloor r \rfloor_h}}_{\lfloor r \rfloor_h,{r}} + \sigma \left( r,X^{Y_{\lfloor r \rfloor_h}}_{\lfloor r \rfloor_h,{r}},z \right)} \right) - f \left( X_{\lceil r \rceil_h,T}^{X^{Y_{\lfloor r \rfloor_h}}_{\lfloor r \rfloor_h,{r}}} \right) \right\} \ \tilde{N} (\delta r, \delta z) \\ 
      + & \int_{0}^{T} \int_{\R \backslash \{ 0 \}} \left\{ f \left( X_{\lceil r \rceil_h,T}^{X^{Y_{\lfloor r \rfloor_h}}_{\lfloor r \rfloor_h,{r}} + \sigma \left( r,X^{Y_{\lfloor r \rfloor_h}}_{\lfloor r \rfloor_h,{r}},z \right)} \right) - f \left( X_{\lceil r \rceil_h,T}^{X^{Y_{\lfloor r \rfloor_h}}_{\lfloor r \rfloor_h,{r}}} \right) \right. \\
      & \left. - \sigma \left( r,X^{Y_{\lfloor r \rfloor_h}}_{\lfloor r \rfloor_h,{r}},z \right) f'\left( X_{\lceil r \rceil_h,T}^{X^{Y_{\lfloor r \rfloor_h}}_{\lfloor r \rfloor_h,r}} \right) \partial_x X_{\lceil r \rceil_h,T}^{X^{Y_{\lfloor r \rfloor_h}}_{\lfloor r \rfloor_h,r}} \right\} \ \nu( \d z) \d r.
\end{split}
\end{equation}
By combining \eqref{eq:decoupage}, \eqref{eq:Decoupage2emeterme} and \eqref{eq:Decoupage1erTerme}, we conclude that for every $n \geqslant 1$ and $h=T/n$:
\begin{equation}\label{eq:FormuleIntermediaireH}
    \begin{split}
         & f \left( X_{0,T}^{Y_0} \right) - f(Y_T) \\ 
         = \quad & \int_{0}^{T} \left\{  f' \left( X_{\lceil r \rceil_h,T}^{X^{Y_{\lfloor r \rfloor_h}}_{\lfloor r \rfloor_h,r}} \right) \partial_x X_{\lceil r \rceil_h,T}^{X^{Y_{\lfloor r \rfloor_h}}_{\lfloor r \rfloor_h,r}} \ \mu \left( r, X^{Y_{\lfloor r \rfloor_h}}_{\lfloor r \rfloor_h,r} \right) - f' \left( X_{\lceil r \rceil_h,T}^{Y_r} \right) \partial_x X_{\lceil r \rceil_h,T}^{Y_r} A_r \right\} \ \d r \\
      + & \int_{0}^{T} \int_{\R \backslash \{ 0 \}} \left\{ f \left( X_{\lceil r \rceil_h,T}^{X^{Y_{\lfloor r \rfloor_h}}_{\lfloor r \rfloor_h,{r}} + \sigma \left( r,X^{Y_{\lfloor r \rfloor_h}}_{\lfloor r \rfloor_h,{r}},z \right)} \right) - f \left( X_{\lceil r \rceil_h,T}^{Y_{r}+B_{r,z}} \right) \right. \\
      & + f \left( X_{\lceil r \rceil_h,T}^{Y_{r}} \right) - f \left( X_{\lceil r \rceil_h,T}^{X^{Y_{\lfloor r \rfloor_h}}_{\lfloor r \rfloor_h,{r}}} \right)  \\
      & \left. +  B_{r,z} f' \left( X_{\lceil r \rceil_h,T}^{Y_r} \right) \partial_x X_{\lceil r \rceil_h,T}^{Y_r} - \sigma \left( r,X^{Y_{\lfloor r \rfloor_h}}_{\lfloor r \rfloor_h,{r}},z \right) f'\left( X_{\lceil r \rceil_h,T}^{X^{Y_{\lfloor r \rfloor_h}}_{\lfloor r \rfloor_h,r}} \right) \partial_x X_{\lceil r \rceil_h,T}^{X^{Y_{\lfloor r \rfloor_h}}_{\lfloor r \rfloor_h,r}} \right\} \ \nu( \d z) \d r \\ 
      + & \int_0^T \int_{\R \backslash \{ 0 \}}  \left\{ f \left( X_{\lceil r \rceil_h,T}^{X^{Y_{\lfloor r \rfloor_h}}_{\lfloor r \rfloor_h,{r}} + \sigma \left( r,X^{Y_{\lfloor r \rfloor_h}}_{\lfloor r \rfloor_h,{r}},z \right)} \right) - f \Big( X_{\lceil r \rceil_h,T}^{Y_{r}+B_{r,z}} \Big)  \right. \\
      & \left. + f \Big( X_{\lceil r \rceil_h,T}^{Y_{r}} \Big) - f \left( X_{\lceil r \rceil_h,T}^{X^{Y_{\lfloor r \rfloor_h}}_{\lfloor r \rfloor_h,{r}}} \right)  \right\} \ \tilde{N} (\delta r, \delta z).
    \end{split}
\end{equation}

\subsection{Passage to the limit $h \to 0$} \label{subsec:cvg_h_zero}

We want to let $h$ going to zero in \eqref{eq:FormuleIntermediaireH}. To that extend, we will analyse the limit when $h \to 0$ of the three terms appearing in the right-hand-side of \eqref{eq:FormuleIntermediaireH}, namely: the Lebesgue integral with respect to $\d r$, the Lebesgue integral with respect to $\nu(\d z) \ \d r$ and the Poisson--Skorohod integral with respect to $\tilde{N}(\delta r, \delta z)$. 

We proceed first by \textbf{the Lebesgue integral with respect to} $\d r$. We apply Vitali's Theorem for finite measures, to prove a convergence in $L^1 ( \P \otimes \lambda_{[0,T]})$, by showing that 
\begin{enumerate}
    \item we have the convergence in measure
    \[ \begin{split}
        &  f' \left( X_{\lceil r \rceil_h,T}^{X^{Y_{\lfloor r \rfloor_h}}_{\lfloor r \rfloor_h,r}} \right) \partial_x X_{\lceil r \rceil_h,T}^{X^{Y_{\lfloor r \rfloor_h}}_{\lfloor r \rfloor_h,r}} \ \mu \left( r, X^{Y_{\lfloor r \rfloor_h}}_{\lfloor r \rfloor_h,r} \right) - f' \left( X_{\lceil r \rceil_h,T}^{Y_r} \right) \partial_x X_{\lceil r \rceil_h,T}^{Y_r} A_r \\
        \xrightarrow[h \to 0^+]{\P \otimes \lambda_{[0,T]}} {} & f' \left( X_{r,T}^{Y_{r}} \right) \partial_x X_{r,T}^{Y_{r}} \ \mu \left( r, Y_r \right) - f' \left( X_{r,T}^{Y_r} \right) \partial_x X_{r,T}^{Y_r} A_r,
    \end{split}  \]

    \item and we have the following uniform $L^2$-bound:
    \begin{multline} \label{eq:borne_integrale_dr}
     \sup_{h \in (0, T)} \E \left[  \int_{0}^{T} \left|  f' \left( X_{\lceil r \rceil_h,T}^{X^{Y_{\lfloor r \rfloor_h}}_{\lfloor r \rfloor_h,r}} \right) \partial_x X_{\lceil r \rceil_h,T}^{X^{Y_{\lfloor r \rfloor_h}}_{\lfloor r \rfloor_h,r}} \ \mu \left( r, X^{Y_{\lfloor r \rfloor_h}}_{\lfloor r \rfloor_h,r} \right) - f' \left( X_{\lceil r \rceil_h,T}^{Y_r} \right) \partial_x X_{\lceil r \rceil_h,T}^{Y_r} A_r \right|^{2} \ \d r  \right] \\ < \infty.    
    \end{multline}
\end{enumerate}
First, we check the convergence in measure. For every $r \in [0,T]$, by the stochastic continuity of $Y$, $X_{\cdot,T}^{\cdot}$, $X_{\cdot,r}^{\cdot}$ and $\partial_xX_{\cdot,T}^{\cdot}$, the continuity of $f'$, and the fact that $X_{r,T}^{X_{r,r}^{Y_r}} = X_{r,T}^{Y_r}$ $\P$-a.s imply that for almost every $r \in [0, T]$
\[ \begin{split}
     & f' \left( X_{\lceil r \rceil_h,T}^{X^{Y_{\lfloor r \rfloor_h}}_{\lfloor r \rfloor_h,r}} \right) \partial_x X_{\lceil r \rceil_h,T}^{X^{Y_{\lfloor r \rfloor_h}}_{\lfloor r \rfloor_h,r}} \ \mu \left( r, X^{Y_{\lfloor r \rfloor_h}}_{\lfloor r \rfloor_h,r} \right) - f' \left( X_{\lceil r \rceil_h,T}^{Y_r} \right) \partial_x X_{\lceil r \rceil_h,T}^{Y_r} A_r \\
     \xrightarrow[h \to 0^+]{\P} {} & f' \left( X_{r,T}^{Y_{r}} \right) \partial_x X_{r,T}^{Y_{r}} \ \mu \left( r, Y_r \right) - f' \left( X_{r,T}^{Y_r} \right) \partial_x X_{r,T}^{Y_r} A_r.
\end{split} \]
Since $\lambda_{[0,T]}$ is a \textbf{finite measure}, it follows directly that 
\[ \begin{split}
        &  f' \left( X_{\lceil r \rceil_h,T}^{X^{Y_{\lfloor r \rfloor_h}}_{\lfloor r \rfloor_h,r}} \right) \partial_x X_{\lceil r \rceil_h,T}^{X^{Y_{\lfloor r \rfloor_h}}_{\lfloor r \rfloor_h,r}} \ \mu \left( r, X^{Y_{\lfloor r \rfloor_h}}_{\lfloor r \rfloor_h,r} \right) - f' \left( X_{\lceil r \rceil_h,T}^{Y_r} \right) \partial_x X_{\lceil r \rceil_h,T}^{Y_r} A_r \\
        \xrightarrow[h \to 0^+]{\P \otimes \lambda_{[0,T]}} {} & f' \left( X_{r,T}^{Y_{r}} \right) \partial_x X_{r,T}^{Y_{r}} \ \mu \left( r, Y_r \right) - f' \left( X_{r,T}^{Y_r} \right) \partial_x X_{r,T}^{Y_r} A_r.
    \end{split}  \]

The $L^2$-uniform bound \eqref{eq:borne_integrale_dr}, whose proof is postponed to Section \ref{subsec:Borne_integrale_dr}, implies in particular that
\begin{equation*}
    \sup_{h \in (0,T)} \Bigg\|  f' \left( X_{\lceil r \rceil_h,T}^{X^{Y_{\lfloor r \rfloor_h}}_{\lfloor r \rfloor_h,r}} \right) \partial_x X_{\lceil r \rceil_h,T}^{X^{Y_{\lfloor r \rfloor_h}}_{\lfloor r \rfloor_h,r}} \ \mu \left( r, X^{Y_{\lfloor r \rfloor_h}}_{\lfloor r \rfloor_h,r} \right) - f' \left( X_{\lceil r \rceil_h,T}^{Y_r} \right) \partial_x X_{\lceil r \rceil_h,T}^{Y_r} A_r \Bigg\|_{L^2(\P \otimes \lambda_{[0,T]})}  \\ < \infty. 
\end{equation*}
In particular, due to La Vallée Poussin's criterion, the family
\begin{equation*}
    \left( f' \left( X_{\lceil r \rceil_h,T}^{X^{Y_{\lfloor r \rfloor_h}}_{\lfloor r \rfloor_h,r}} \right) \partial_x X_{\lceil r \rceil_h,T}^{X^{Y_{\lfloor r \rfloor_h}}_{\lfloor r \rfloor_h,r}} \ \mu \left( r, X^{Y_{\lfloor r \rfloor_h}}_{\lfloor r \rfloor_h,r} \right) - f' \left( X_{\lceil r \rceil_h,T}^{Y_r} \right) \partial_x X_{\lceil r \rceil_h,T}^{Y_r} A_r \right)_{h \in (0,T)}
\end{equation*}
is $\P \otimes \lambda_{[0,T]}$-uniformly integrable.
Then by Vitali's Theorem for finite measures we obtain the convergence
\begin{multline*}
    \lim_{h \to 0} \Bigg\| f' \left( X_{\lceil r \rceil_h,T}^{X^{Y_{\lfloor r \rfloor_h}}_{\lfloor r \rfloor_h,r}} \right) \partial_x X_{\lceil r \rceil_h,T}^{X^{Y_{\lfloor r \rfloor_h}}_{\lfloor r \rfloor_h,r}} \ \mu \left( r, X^{Y_{\lfloor r \rfloor_h}}_{\lfloor r \rfloor_h,r} \right) - f' \left( X_{\lceil r \rceil_h,T}^{Y_r} \right) \partial_x X_{\lceil r \rceil_h,T}^{Y_r} A_r \\
    - \left(f' \left( X_{r,T}^{Y_{r}} \right) \partial_x X_{r,T}^{Y_{r}} \ \mu \left( r, Y_r \right) - f' \left( X_{r,T}^{Y_r} \right) \partial_x X_{r,T}^{Y_r} A_r\right)  \Bigg\|_{L^1(\P \otimes \lambda_{[0,T]})} = 0.
\end{multline*}
This does imply the convergence 
\begin{multline*}
    \lim_{h \to 0} \E \left[ \left| \int_0^T \Bigg\{ f' \left( X_{\lceil r \rceil_h,T}^{X^{Y_{\lfloor r \rfloor_h}}_{\lfloor r \rfloor_h,r}} \right) \partial_x X_{\lceil r \rceil_h,T}^{X^{Y_{\lfloor r \rfloor_h}}_{\lfloor r \rfloor_h,r}} \ \mu \left( r, X^{Y_{\lfloor r \rfloor_h}}_{\lfloor r \rfloor_h,r} \right) - f' \left( X_{\lceil r \rceil_h,T}^{Y_r} \right) \partial_x X_{\lceil r \rceil_h,T}^{Y_r} A_r \right. \right. \\
    - \left. \left. \left( f' \left( X_{r,T}^{Y_{r}} \right) \partial_x X_{r,T}^{Y_{r}} \ \mu \left( r, Y_r \right) - f' \left( X_{r,T}^{Y_r} \right) \partial_x X_{r,T}^{Y_r} A_r\right) \Bigg\} \ \d r \right| \right] = 0.
\end{multline*}

We consider now the \textbf{second Lebesgue integral with respect to} $\nu(\d z) \ \d r$. The situation is more intricate, because $\nu$ is not a finite measure. To simplify the presentation, we will use the notations
\begin{align*}
    F_h(r,z) \defeq {} &  f \left( X_{\lceil r \rceil_h,T}^{X^{Y_{\lfloor r \rfloor_h}}_{\lfloor r \rfloor_h,{r}} + \sigma \left( r,X^{Y_{\lfloor r \rfloor_h}}_{\lfloor r \rfloor_h,{r}},z \right)} \right) - f \left( X_{\lceil r \rceil_h,T}^{Y_{r}+B_{r,z}} \right) + f \left( X_{\lceil r \rceil_h,T}^{Y_{r}} \right) - f \left( X_{\lceil r \rceil_h,T}^{X^{Y_{\lfloor r \rfloor_h}}_{\lfloor r \rfloor_h,{r}}} \right)  \\
    & {} +  B_{r,z} f' \left( X_{\lceil r \rceil_h,T}^{Y_r} \right) \partial_x X_{\lceil r \rceil_h,T}^{Y_r} - \sigma \left( r,X^{Y_{\lfloor r \rfloor_h}}_{\lfloor r \rfloor_h,{r}},z \right) f'\left( X_{\lceil r \rceil_h,T}^{X^{Y_{\lfloor r \rfloor_h}}_{\lfloor r \rfloor_h,r}} \right) \partial_x X_{\lceil r \rceil_h,T}^{X^{Y_{\lfloor r \rfloor_h}}_{\lfloor r \rfloor_h,r}},
\end{align*}
and
\begin{multline*}
    F(r,z) \defeq f \left( X_{r,T}^{Y_r + \sigma \left( r,Y_r,z \right)} \right) - f \left( X_{r,T}^{Y_{r}+B_{r,z}} \right) \\
      +  B_{r,z} f' \left( X_{r,T}^{Y_r} \right) \partial_x X_{r,T}^{Y_r} - \sigma \left( r,Y_r,z \right) f'\left( X_{r,T}^{Y_r} \right) \partial_x X_{r,T}^{Y_r}.
\end{multline*}

The objective is to apply Vitali's $L^1$-convergence Theorem \ref{thm:vitali_general} for $\sigma$-finite measure spaces. We will rely on the following Lemma, proven in Section \ref{subsec:Borne_integrale_dNu(z)dr}.

\begin{Lemma} \label{lm:CroissancePolynomialeDeFhEnZ}
    There exists a positive function $\psi \in L^1(\nu) \cap L^2(\nu)$ such that for every $z \in \R \backslash \{ 0 \}$,
    \begin{equation*}
        \sup_{h > 0} \| F_h(\cdot,z) \|_{L^2(\P \otimes \lambda_{[0,T]})} \leq \psi(z) \ \text{ and } \ \| F(\cdot,z) \|_{L^2(\P \otimes \lambda_{[0,T]})} \leq \psi(z).
     \end{equation*}
\end{Lemma}

Firstly, using similar arguments than for the first Lebesgue integral, we have that for almost every $r \in [0,T]$ and $z \in \R \backslash \{0 \}$, 
\begin{equation*}
    F_h(r,z) \xrightarrow[h \to 0^+]{\P} F(r,z). 
\end{equation*}
This is, however, not enough to conclude the convergence in measure for $\P \otimes \lambda_{[0,T]} \otimes \nu$ because $\nu$ is only $\sigma$-finite. Instead, fix $\eps > 0$. We apply Markov inequality, Cauchy--Schwarz inequality and Lemma \ref{lm:CroissancePolynomialeDeFhEnZ}:
\begin{align*}
    \E \left[ \int_0^T \mathbf{1}_{ \big\{ |F_h(r,z) - F(r,z)| > \eps \big\} } \ \d r \right] & \leqslant \frac{1}{\eps} \left\| F_h(\cdot,z) - F(\cdot,z) \right\|_{L^1(\P \otimes \lambda_{[0,T]})} \\ 
    & \leqslant  \frac{\sqrt{T}}{\eps} (\sup_{h > 0} \| F_h(\cdot,z) \|_{L^2(\P \otimes \lambda_{[0,T]})} + \| F(\cdot,z) \|_{L^2(\P \otimes \lambda_{[0,T]})}) \\
     & \leqslant \frac{\sqrt{T}}{\eps} \psi(z).
\end{align*}
By the dominated convergence Theorem it follows that for every $\eps > 0$,
\begin{equation*}
    \int_{\R \backslash \{ 0 \}} \E \left[ \int_0^T \mathbf{1}_{ \{ |F_h(r,z) - F(r,z)| > \eps \} } \ \d r \right] \ \nu (\d z) \xrightarrow[h \to 0^+]{} 0.
\end{equation*}
In other words, we have the limit in measure $\lim_{h \to 0^+} F_h = F$ for the product measure $\P \otimes \lambda_{[0,T]} \otimes \nu$. 

Secondly, we will show that the family $(F_h)_{h \in (0,1]}$ has uniformly absolutely continuous integrals. Observe that due to Lemma \ref{lm:CroissancePolynomialeDeFhEnZ} we have for every $h \in (0,1]$ and $z \in \mathbb{R} \backslash \{ 0 \}$,
\begin{equation*}
    \| F_h(\cdot,z) \|_{L^1(\P \otimes \lambda_{[0,T]})} \leq T \| F_h(\cdot,z) \|_{L^2(\P \otimes \lambda_{[0,T]})} \leq T \psi(z),
\end{equation*}
where $\psi \in L^1(\nu)$. It follows that $F_h \in L^1(\P \otimes \lambda_{[0,T]} \otimes \nu)$. 
Then, using again Lemma \ref{lm:CroissancePolynomialeDeFhEnZ}, we have that for any $h\in (0,T)$
\[ \| F_h \|_{L^2(\P \otimes \lambda_{[0,T]} \otimes \nu)} = \int_{\R \backslash \{ 0 \}} \| F_h (\cdot, z) \|_{L^2(\P \otimes \lambda_{[0,T]})}^2 \ \nu(\d z) \leqslant \int_{\R \backslash \{ 0 \}} \psi(z)^2 \ \nu(\d z), \]
hence since $\psi \in L^2(\nu)$, we obtain the uniform upper-bound
\begin{equation} \label{eq:BorneMomentIntegraleNudZdR}
         \sup_{h > 0} \| F_h \|_{L^2(\P \otimes \lambda_{[0,T]} \otimes \nu)} < \infty.
\end{equation}
This allows us to apply Lemma \ref{lm:critere_vallee_poussin} with $G(t) = t^2$, which shows that $(F_h)_{h \in (0,1]}$ has uniformly absolutely continuous integrals.

Finally, we have to show the tension criterion. As $\int_{\R \backslash \{ 0 \}} \psi(z) \nu(\d z) < + \infty$, the dominated convergence Theorem implies that for any $\eps > 0$ there exists $\delta_{\varepsilon} > 0$ such that
\begin{equation*}
    \int_{-\delta_{\varepsilon}}^{\delta_{\varepsilon}} \psi(z) \nu(\d z) \leq \frac{\eps}{\sqrt{T}}. 
\end{equation*}
Hence, by letting $K_{\eps} = \Omega \times [0,T] \times \R \backslash [-\delta_{\varepsilon},\delta_{\eps}]$, using Fubini--Tonelli's Theorem, Cauchy--Schwarz inequality and Lemma \ref{lm:CroissancePolynomialeDeFhEnZ}, we obtain the inequality
\begin{align*}
    \sup_{h \in (0, T)} \int_{ (\Omega \times [0,T] \times \R) \backslash K_{\eps} } |F_h(r,z,\omega)| \d (\P \otimes \lambda_{[0,T]} \otimes \nu) (r, z, \omega)  & = \sup_{h \in (0, T)}  \E \left[ \int_0^T \int_{-\delta}^\delta |F_h(r,z)| \ \d r \ \nu(\d z) \right] \\ & \leqslant \sqrt{T} \int_{-\delta}^\delta \psi(z) \nu(\d z) \leqslant \eps.
\end{align*}
This shows that the family $(F_h)_{h \in (0, T)}$ is tight. Finally, we conclude using Vitali's $L^1$-convergence Theorem \ref{thm:vitali_general} for $\sigma$-finite measure spaces that 
\begin{equation*}
    \lim_{h \to 0^+} \| F_h - F \|_{L^1(\P \otimes \lambda_{[0,T]} \otimes \nu)} = 0.
\end{equation*}

We deal now with the \textbf{third integral}, the Skorohod one. We use Lemma \ref{LemmeConvergenceDelta} to prove a weak convergence in $L^2$ of this integral. To that extend, we need to prove the two following statements. 

\begin{enumerate}
    \item The family of processes
    \begin{multline*}
        \left\{  f \left( X_{\lceil r \rceil_h,T}^{X^{Y_{\lfloor r \rfloor_h}}_{\lfloor r \rfloor_h,{r}} + \sigma \left( r,X^{Y_{\lfloor r \rfloor_h}}_{\lfloor r \rfloor_h,{r}},z \right)} \right) - f \Big( X_{\lceil r \rceil_h,T}^{Y_{r}+B_{r,z}} \Big) \right.  \\ \left. +  f \Big( X_{\lceil r \rceil_h,T}^{Y_{r}} \Big)  -  f \left( X_{\lceil r \rceil_h,T}^{X^{Y_{\lfloor r \rfloor_h}}_{\lfloor r \rfloor_h,{r}}} \right), r \in [0, T], z \in \R \backslash \{ 0 \} \right\}
    \end{multline*}
     weakly converges in $L^2(\Omega \times E)$ when $h$ goes to $0^+$ toward
     \[ \left\{  f \left( X_{r,T}^{Y_{r} + \sigma \left( r,Y_{r},z \right)} \right) - f \Big( X_{r,T}^{Y_{r}+B_{r,z}} \Big)  , r \in [0, T], z \in \R \backslash \{ 0 \} \right\}.
     \]

     \item We have the bound
     \begin{multline}\label{eq:Borne_Moment_Int_3eme_Terme}
         \sup_{h > 0}  \E \left[ \left( \int_0^T \int_{\R \backslash \{ 0 \}}  \left\{ f \left( X_{\lceil r \rceil_h,T}^{X^{Y_{\lfloor r \rfloor_h}}_{\lfloor r \rfloor_h,{r}} + \sigma \left( r,X^{Y_{\lfloor r \rfloor_h}}_{\lfloor r \rfloor_h,{r}},z \right)} \right) - f \Big( X_{\lceil r \rceil_h,T}^{Y_{r}+B_{r,z}} \Big)  \right. \right. \right. \\
       \left. \left. \left. +  f \left( X_{\lceil r \rceil_h,T}^{X^{Y_{\lfloor r \rfloor_h}}_{\lfloor r \rfloor_h,{r}}} \right)  - f \Big( X_{\lceil r \rceil_h,T}^{Y_{r}} \Big)  \right\} \ \tilde{N} (\delta r, \delta z) \right)^2 \right] < \infty. 
     \end{multline}
\end{enumerate}

We may use the notations
\begin{align*}
    G_h(r,z) & \defeq f \left( X_{\lceil r \rceil_h,T}^{X^{Y_{\lfloor r \rfloor_h}}_{\lfloor r \rfloor_h,{r}} + \sigma \left( r,X^{Y_{\lfloor r \rfloor_h}}_{\lfloor r \rfloor_h,{r}},z \right)} \right) - f \Big( X_{\lceil r \rceil_h,T}^{Y_{r}+B_{r,z}} \Big)  +  f \left( X_{\lceil r \rceil_h,T}^{X^{Y_{\lfloor r \rfloor_h}}_{\lfloor r \rfloor_h,{r}}} \right)  - f \Big( X_{\lceil r \rceil_h,T}^{Y_{r}} \Big), \\
    G(r,z) &\defeq f \left( X_{r,T}^{Y_{r} + \sigma \left( r,Y_{r},z \right)} \right) - f \Big( X_{r,T}^{Y_{r}+B_{r,z}} \Big).
\end{align*}
Firsty, we will prove point 1. The goal is to apply Proposition \ref{prop:measure_CVG_implies_weakL2_CVG}, that is the fact that $L^2$-bounded convergence in measure implies
$L^2$-weak convergence. For that purpose, we need two things: the convergence of $G_h$ to $G$ in measure for the product measure $\P \otimes \lambda_{[0,T]} \otimes \nu$, and the uniform $L^2$-bound
\begin{equation} \label{eq:uniform_L2_bound_troisieme_integrande}
    \sup_{h > 0} \| G_h \|_{L^2(\P \otimes \lambda_{[0,T]} \otimes \nu)} < + \infty.
\end{equation}
This bound is in fact straightforward to obtain by combining the results of Lemmas \ref{lm:controle_l2_premier_terme} and \ref{lm:controle_l2_second_terme} with the standard inequality $(x+y)^2 \leqslant 2 (x^2 + y^2)$, leading for every $h \in (0,T)$ to 
\begin{equation*}
      \| G_h \|_{L^2(\P \otimes \lambda_{[0,T]} \otimes \nu)} \leqslant T( \|\phi_1 \|_{L^1(\nu)} + \| \phi_2 \|_{L^1(\nu)}) < + \infty.
\end{equation*}

For the convergence in probability, we proceed as for the second integral. Firstly, observe that for any $r\in[0,T]$ and $z \in \R \backslash \{0 \}$, $G_h(r,z)$ converges in probability when $h \to 0^+$ toward $ G(r,z)$. Then, applying Markov inequality, for any $\eps > 0$,
\begin{align*}
    \E \left[ \mathbf{1}_{ \left\{ | G_h(r,z) - G(r,z)| < \eps \right\} } \right]   & = \E \left[ \mathbf{1}_{ \left\{ | G_h(r,z) - G(r,z) |^2 < \eps^2 \right\} } \right] \\
    & \leqslant \frac{1}{\eps^{2}} \E \left[ \big| G_h(r,z) - G(r,z) \big|^2 \right].
\end{align*}
By Lemmas \ref{lm:controle_l2_second_terme} and \ref{lm:controle_l2_premier_terme} we have that for almost every $r \in [0,T]$ and $z \in \R \backslash \{ 0 \}$, 
\begin{equation*}
    \sup_{h \in (0,T)} \E \left[ G_h(r,z)^2 \right] \leqslant \phi_1(z) + \phi_2(z).
\end{equation*}
Then, using similar arguments as in the proof of Lemmas \ref{lm:controle_l2_second_terme} and \ref{lm:controle_l2_premier_terme}, it is straightforward to prove the existence of a positive function $\phi_3 \in L^1(\nu)$ such that for almost every $r \in [0,T]$ and $z \in \R \backslash \{ 0 \}$,
\begin{equation*}
     \E \left[ G(r,z)^2 \right] \leqslant \phi_3(z).
\end{equation*}

Using the dominated convergence Theorem and the fact that $\lambda_{[0,T]}$ is a finite measure, we conclude that the convergence in measure holds for the product measure $\P \otimes \lambda_{[0,T]} \otimes \nu$. Together with the $L^2(\P \otimes \lambda_{[0,T]} \otimes \nu)$-uniform bound \eqref{eq:uniform_L2_bound_troisieme_integrande}, this allows to apply Proposition \ref{prop:measure_CVG_implies_weakL2_CVG}, demonstrating the expected weak $L^2$-convergence result.
Finally, the inequality \eqref{eq:Borne_Moment_Int_3eme_Terme} announced in Point 2 is postponed to Section \ref{subsec:Borne_integrale_Sto}.

\subsection{Conclusion}

In view of the three previous convergence ($L^1$, $L^1$ and weak $L^2$), we conclude that for every $Z \in L^{\infty} (\mathbf{F}_T)$, 
\[ \begin{split}
    & \E \left[ Z \left( f \left( X_{0,T}^{Y_0} \right) - f(Y_T) - \int_{0}^{T} \left\{  f' \left( X_{r,T}^{Y_{r}} \right) \partial_x X_{r,T}^{Y_{r}} \ \mu \left( r, Y_r \right) - f' \left( X_{r,T}^{Y_r} \right) \partial_x X_{r,T}^{Y_r} A_r \right\} \ \d r \right. \right. \\
    & - \int_{0}^{T} \int_{\R \backslash \{ 0 \}} \left\{ f \left( X_{r,T}^{Y_r + \sigma \left( r,Y_r,z \right)} \right) - f \left( X_{r,T}^{Y_{r}+B_{r,z}} \right) \right. \\
      & \left. +  B_{r,z} f' \left( X_{r,T}^{Y_r} \right) \partial_x X_{r,T}^{Y_r} - \sigma \left( r,Y_r,z \right) f'\left( X_{r,T}^{Y_r} \right) \partial_x X_{r,T}^{Y_r} \right\} \ \nu( \d z) \d r \\
    & \left. \left. - \int_0^T \int_{\R \backslash \{ 0 \}} \left( f \left( X_{r,T}^{Y_{r} + \sigma \left( r,Y_{r},z \right)} \right) - f \Big( X_{r,T}^{Y_{r}+B_{r,z}} \Big) \right) \ \tilde{N}(\delta r, \delta z) \right) \right] = 0.
\end{split} \]
This concludes the proof. 

\section{Proofs of technical lemmas and estimates} \label{sec:Technical_Lemmas}

\subsection{Proof of Lemma \ref{lm:controle_l2_second_terme} } \label{proof:controle_l2_second_terme}

Let $ \mathbb{N} \backslash\{0 \}, h = T/n, i \in \{ 0, \ldots, n-1\}$, $r \in [ih,(i+1)h)$ and $z \in \R \backslash \{ 0 \}$. Using Taylor's formula with integral leftover, we may factorize the integrand by the jump coefficient $B_{r,z}$:
\begin{equation*}
     f \Big( X_{(i+1)h,T}^{Y_{r}+B_{r,z}} \Big) - f \Big( X_{(i+1)h,T}^{Y_{r}} \Big) = B_{r,z} \ \int_0^1 f'\left( X_{(i+1)h,T}^{\lambda B_{r,z} + Y_r} \right) \ \partial_x X_{(i+1)h,T}^{\lambda B_{r,z} + Y_{r}} \ \d \lambda.
\end{equation*}
It follows, after the use of Jensen's inequality that
\begin{equation*}
     \E \left[   \left( f \Big( X_{(i+1)h,T}^{Y_{r}+B_{r,z}} \Big) - f \Big( X_{(i+1)h,T}^{Y_{r}} \Big) \right)^2  \right]  \leq  \int_0^1 \E \left[ \left( B_{r,z} \ f'(X_{(i+1)h,T}^{\lambda B_{r,z} + Y_r}) \ \partial_x X_{(i+1)h,T}^{\lambda B_{r,z} + Y_{r}} \right)^2\right] \d \lambda.  
\end{equation*}
Let $\lambda \in [0, 1]$. We apply Hölder's inequality with $a = \frac{p}{2}$, $b = \frac{p}{2q}$ and $c = \frac{p}{p-2(q+1)}$, so that we have $\frac{1}{a} + \frac{1}{b} + \frac{1}{c} = 1$. We deduce that
\begin{multline*}
    \E \left[ \left( B_{r,z} \ f'(X_{(i+1)h,T}^{\lambda B_{r,z} + Y_r}) \ \partial_x X_{(i+1)h,T}^{\lambda B_{r,z} + Y_{r}} \right)^2\right] \\ \leqslant \E \left[ \big| B_{r,z} \big| ^p \right]^{\frac{2}{p}} \ \E \left[ \left| f' \left( X_{(i+1)h,T}^{\lambda B_{r,z} + Y_r} \right) \right|^{\frac{p}{q}} \right]^{\frac{2q}{p}} \ \E \left[ \left| \partial_x X_{(i+1)h,T}^{\lambda B_{r,z} + Y_{r}} \right|^\frac{2p}{p-2(1+q)} \right]^{\frac{p-2(1+q)}{p}}.
\end{multline*}
We will analyse the three terms of the right-hand side of the latter inequality one by one. 

First, using \eqref{hyp:coeffs} we have
\[ \E \left[ \big| B_{r,z} \big|^p \right] \leqslant |z|^{p k} \ \E \left[ \left( \overline{B}_r \right)^p \right] \leqslant |z|^{p k} \ \sup_{r \in [0,T]}  \E \left[ \left( \overline{B}_r \right)^p \right]. \]
It follows that
\begin{equation*}
    \sup_{r \in [ih,(i+1)h)} \E \left[ \big| B_{r,z} \big|^p \right]^\frac{2}{p} \leqslant C_0 \  |z|^{2p},
\end{equation*}
where $C_0 \coloneqq \sup_{r \in [0,T]} \E[ (\overline{B}_r)^p ]^{2/p} \in (0,+\infty)$. 

Then, using \eqref{hyp:testfunction} we have by using convexity of $x \mapsto x^{p/q}$ (since $p/q > 2 + 2/q > 1$)
\begin{align*}
    \E \left[ \left| f'\left( X_{(i+1)h,T}^{\lambda B_{r,z} + Y_r} \right) \right|^{\frac{p}{q}} \right] & \leqslant C_f^\frac{p}{q} \ \E \left[ \left( 1 + \left| X_{(i+1) h, T}^{\lambda B_{r,z} + Y_r} \right|^q \right)^\frac{p}{q} \right] \\
    & \leqslant 2^{p/q-1} C_f^{p/q} \left( 1 + \E \left[ \left| X_{(i+1) h, T}^{\lambda B_{r,z} + Y_r} \right|^p \right] \right).
\end{align*}
Let us write 
\begin{equation*}
    \E \left[ \left| X_{(i+1) h, T}^{\lambda B_{r,z} + Y_r} \right|^p \right] = \E \left[ \E \left[ \left| X_{(i+1) h, T}^{\lambda B_{r,z} + Y_r} \right|^p \ \bigg| \ \mathbf{F}_{[ih,(i+1)h]}  \right] \right].
\end{equation*}
Then $\lambda B_{r,z} + Y_r$ is $\mathbf{F}_{[ih,(i+1)h]}$-measurable for every $r \in [ih,(i+1)h]$, while the family of random variable $ (X_{(i+1) h, T}^x, \ x \in \R)$ is independent of $\mathbf{F}_{[ih,(i+1)h]}$. Thus from the properties of conditional expectation we deduce that 
\begin{equation*}
    \E \left[ \left| X_{(i+1) h, T}^{\lambda B_{r,z} + Y_r} \right|^p \right] = \E \left[ \phi \big( \lambda B_{r,z} + Y_r \big) \right],
\end{equation*}
where for every $x \in \R$, $\phi(x) \coloneqq \E[ |X_{(i+1) h, T}^x|^p]$. Now due to \eqref{hyp:MomentsX} we have 
\begin{equation*}
    \phi(x) \leqslant C_{X,p} \ (1 + |x|^{\mx}),
\end{equation*}
leading to 
\begin{align*}
    \E \left[ \left| X_{(i+1) h, T}^{\lambda B_{r,z} + Y_r} \right|^p \right]  & \leqslant C_{X,p} \ \left( 1 + \E \left[ \big| \lambda B_{r,z} + Y_r \big|^{\mx} \right] \right) \\
   &  \leqslant C_{X,p} \left( 1 + 2^{\mx-1}  \E \left[ \big| B_{r,z} \big|^{\mx} \right] + 2^{\mx-1} \E \left[ \big| Y_r \big|^{\mx} \right] \right),
\end{align*}
where we used that $\lambda^{\mx} \leqslant 1$, and the convexity of $x \mapsto x^{\mx}$. Applying \eqref{hyp:coeffs} and \eqref{hyp:MomentsY}, we get 
\[ \E \left[ \left| X_{(i+1) h, T}^{\lambda B_{r,z} + Y_r} \right|^p \right] \leqslant C_{X,p} \left( 1 + 2^{\mx-1} |z|^{k \mx} \sup_{r \in [0,T]}  \E \left[ \left( \overline{B}_r \right)^{\mx} \right] + 2^{\mx-1} \sup_{r \in [0,T]} \E \left[ \big| Y_r \big|^{\mx} \right] \right). \]
Finally, by using $(x+y)^{q/2p} \leqslant \max\{ 1, 2^{q/2p-1} \} (x^{q/2p} + y^{q/2p})$, for $x,y > 0$, we obtain the existence of a finite constant $C_1 > 0$ depending on $X,Y,\overline{B},p,q,\mx$ such that 
\begin{equation*}
    \sup_{r \in [ih,(i+1)h)} \sup_{\lambda \in [0,1]} \E \left[ \left| f' \left( X_{(i+1)h,T}^{\lambda B_{r,z} + Y_r} \right) \right|^{\frac{p}{q}} \right]^{\frac{2q}{p}} \leqslant C_1 \left( 1 +  |z|^{\frac{2 k \mx q}{p}} \right).
\end{equation*}

To obtain an upper-bound for the last term, the arguments are quite similar except that we use \eqref{hyp:MomentsX}
instead of \eqref{hyp:MomentsX}. We obtain the existence of $C_2 > 0$ depending on $\partial_xX,Y,\overline{B},p,q,\mdx$ such that 
\begin{equation*}
    \sup_{r \in [ih,(i+1)h)} \sup_{\lambda \in [0,1]} \E \left[ \left| \partial_x X_{(i+1)h,T}^{\lambda B_{r,z} + Y_{r}} \right|^\frac{2p}{p-2(1+q)} \right]^{\frac{p-2(1+q)}{p}} \leq C_2 \left( 1 +  |z|^{\frac{k \mdx (p - 2(1+q))}{p}} \right).
\end{equation*}

We conclude that for every $r \in [ih, (i+1)h)$,
\begin{equation*}
\begin{split}
    \E \left[  \left( f \Big( X_{(i+1)h,T}^{Y_{r}+B_{r,z}} \Big) - f \Big( X_{(i+1)h,T}^{Y_{r}} \Big) \right)^2   \right] & \leqslant  C_0 C_1 C_2 |z|^{2k}  \left( 1 + |z|^\frac{2 k \mx q}{p} \right) \left( 1 + |z|^\frac{k \mdx (p - 2(q+1))}{p} \right) \\
    &:= \phi_1 (z).
\end{split}
\end{equation*}
Moreover, using the moment hypothesis for $\nu$ given by \eqref{hyp:mesure_levy}, we conclude that $\phi_1 \in L^1 (\nu)$. 

\subsection{Proof of Lemma \ref{lm:controle_l2_premier_terme}} \label{proof:controle_l2_premier_terme}

The proof is close to that of Lemma \ref{lm:controle_l2_second_terme}. Let $h \in (0, T)$, $0 \leqslant i \leqslant \lfloor T/h \rfloor$, $r \in [ih, (i+1)h)$ and $z \in \R \backslash \{ 0 \}$. First, we apply Taylor's formula with integral leftover, leading to
\begin{multline*}
    f \left( X_{(i+1)h,T}^{X_{ih,r}^{Y_{ih}} + \sigma(r,X_{ih,r}^{Y_{ih}},z)} \right) -  f \left( X_{(i+1)h,T}^{X_{ih,r}^{Y_{ih}}} \right) \\ = \sigma \left( r,X_{ih,r}^{Y_{ih}},z \right) \int_0^1 f' \left( X_{(i+1)h,T}^{X_{ih,r}^{Y_{ih}} + \lambda \sigma \left( r,X_{ih,r}^{Y_{ih}},z \right)} \right) \partial_x X_{(i+1)h,T}^{X_{ih,r}^{Y_{ih}} + \lambda \sigma \left( r,X_{ih,r}^{Y_{ih}},z \right)} \ \d \lambda.
\end{multline*}
Then, we apply Jensen's inequality, Fubini-Tonelli's Theorem and Hölder's inequality to get
\begin{multline*}
    \E \left[  \left( f \left( X_{(i+1)h,T}^{X_{ih,r}^{Y_{ih}} + \sigma(r,X_{ih,r}^{Y_{ih}},z)} \right) -  f \left( X_{(i+1)h,T}^{X_{ih,r}^{Y_{ih}}} \right) \right)^2  \right] \\ 
    \leqslant \int_0^1 \left\{ \E \left[ \left| \sigma \left( r,X_{ih,r}^{Y_{ih}},z \right) \right|^p \right]^{2/p} \ \E \left[ \left| f'\left( X_{(i+1)h,T}^{X_{ih,r}^{Y_{ih}} + \lambda \sigma(r,X_{ih,r}^{Y_{ih}},z)} \right) \right|^{\frac{p}{q}} \right]^\frac{2q}{p} \right. \\ \left. \times \E \left[ \left| \partial_x X_{(i+1)h,T}^{X_{ih,r}^{Y_{ih}}+ \lambda \sigma(r,X_{ih,r}^{Y_{ih}},z)} \right|^\frac{2p}{p-2(q+1)} \right]^\frac{p-2(q+1)}{p} \right\} \ \d \lambda .
\end{multline*}
Making use of \eqref{hyp:coeffs} we have
\[ \E \left[ \left| \sigma \left( r,X_{ih,r}^{Y_{ih}},z \right) \right|^p \right] \leqslant C_\sigma^p \ |z|^{pk} \ \E \left[ \left( 1 + \left| X_{ih,r}^{Y_{ih}} \right|^{\ms} \right)^p \right] \leqslant 2^{p-1} C_\sigma^p |z|^{pk} \left( 1 + \E \left[ \left| X_{ih,r}^{Y_{ih}} \right|^{\ms p} \right] \right). \]
Then conditioning by $\mathbf{F}_{[0,ih]}$, using that $Y_{ih}$ is $\mathbf{F}_{[0,ih]}$-measurable and that $(X_{ih,r}^x , \ x \in \R)$ is independent from $\mathbf{F}_{[0,ih]}$, we obtain 
\begin{equation*}
    \E \left[ \left| X_{ih,r}^{Y_{ih}} \right|^{\ms p} \right] = \E \big[ \phi(Y_{ih}) \big]
\end{equation*}
where for every $x \in \R$, by using \eqref{hyp:MomentsX}
\[  \phi(x):= \E \left[ \left| X_{ih,r}^x \right|^{ \ms p} \right] \leqslant C_{X,p} (1 + |x|^{\msx}). \] 
Hence
\begin{equation*}
    \E \left[ \left| \sigma \left( r,X_{ih,r}^{Y_{ih}},z \right) \right|^p \right] \leqslant 2^{p-1} C_\sigma^p |z|^{pk} \left( 1 + C_{X,p} \left(1 + \E \left[ \big| Y_{ih} \big|^{\msx} \right] \right) \right).
\end{equation*}
It follows from \eqref{hyp:MomentsY} that there exists a finite
constant $C_0>0$ such that 
\begin{equation*}
    \E \left[ \left| \sigma \left( r,X_{ih,r}^{Y_{ih}},z \right) \right|^p \right]^\frac{2}{p} \leqslant C_0 |z|^{2k}.
\end{equation*}

To estimate the two other terms, we condition successively by $\mathbf{F}_{[ih,(i+1)h]}$ and by $\mathbf{F}_{[0,ih]}$, we use the moment bounds from \eqref{hyp:testfunction}, from \eqref{hyp:MomentsX} and from \eqref{hyp:MomentsX}. After some tedious computations, but without any additional difficulties, we obtain the existence of finite constants $C_1,C_2 > 0$ such that
\begin{align*}
     \E \left[ \left| f' \left( X_{(i+1)h,T}^{X_{ih,r}^{Y_{ih}}+ \lambda \sigma(r,X_{ih,r}^{Y_{ih}},z)} \right) \right|^{\frac{p}{q}} \right]^\frac{2q}{p} & \leqslant C_1 \left( 1 + |z|^{\frac{2 k \mx q}{p}} \right), \\
     \E \left[ \left| \partial_x X_{(i+1)h,T}^{X_{ih,r}^{Y_{ih}}+ \lambda \sigma(r,X_{ih,r}^{Y_{ih}},z)} \right|^\frac{2p}{p-2(q+1)} \right]^\frac{p-2(q+1)}{p} & \leqslant C_2 \left( 1 + |z|^{\frac{k \mdx (p-2(q+1))}{p}} \right).
\end{align*}
The result follows by letting 
\begin{equation*}
    \phi_2 (z) \defeq C_0 C_1 C_2 |z|^{2k} \left( 1 +  |z|^{\frac{2 k \mx q}{p}} \right) \left( 1 +  |z|^{\frac{k \mdx (p-2(q+1))}{p}} \right),
\end{equation*}
which defines a map in $L^1 (\nu)$ due to \eqref{hyp:mesure_levy}.

\subsection{Proof of the estimate \eqref{eq:borne_integrale_dr}}\label{subsec:Borne_integrale_dr}

Let $n \in \N$ with $n \geqslant 1$ fixed and $h = \frac{T}{n}$. We apply Fubini-Tonelli's Theorem and make use of the standard inequality $(X+Y)^2 \leq 2(X^2 + Y^2)$:
\begin{multline*}
    \E \left[ \int_{0}^{T}\left|  f' \left( X_{\lceil r \rceil_h,T}^{X^{Y_{\lfloor r \rfloor_h}}_{\lfloor r \rfloor_h,r}} \right) \partial_x X_{\lceil r \rceil_h,T}^{X^{Y_{\lfloor r \rfloor_h}}_{\lfloor r \rfloor_h,r}} \ \mu \left( r, X^{Y_{\lfloor r \rfloor_h}}_{\lfloor r \rfloor_h,r} \right) - f' \left( X_{\lceil r \rceil_h,T}^{Y_r} \right) \partial_x X_{\lceil r \rceil_h,T}^{Y_r} A_r \right|^{2} \ \d r  \right] \\
    \leqslant 2 \int_0^T \E \left[ f' \left( X_{\lceil r \rceil_h,T}^{X^{Y_{\lfloor r \rfloor_h}}_{\lfloor r \rfloor_h,r}} \right)^2 \left( \partial_x X_{\lceil r \rceil_h,T}^{X^{Y_{\lfloor r \rfloor_h}}_{\lfloor r \rfloor_h,r}} \right)^2 \ \mu \left( r, X^{Y_{\lfloor r \rfloor_h}}_{\lfloor r \rfloor_h,r} \right)^2 \right] \d r \\ 
    + \int_0^T \E\left[ f' \left( X_{\lceil r \rceil_h,T}^{Y_r} \right)^2 \left( \partial_x X_{\lceil r \rceil_h,T}^{Y_r} \right)^2 \ A_r^2   \right] \d r.
\end{multline*}
Let $r \in [0, T]$. In both terms appearing in the right-hand side of the latter inequality, we apply Hölder inequality with $a = \frac{p}{2q}$, $b = \frac{p}{p-2(1+q)}$ and $c=\frac{p}{2}$. 
It follows that 
\begin{multline} \label{eq:drift_estimate_holder_X}
    \E \left[ f' \left( X_{\lceil r \rceil_h,T}^{X^{Y_{\lfloor r \rfloor_h}}_{\lfloor r \rfloor_h,r}} \right)^2 \left( \partial_x X_{\lceil r \rceil_h,T}^{X^{Y_{\lfloor r \rfloor_h}}_{\lfloor r \rfloor_h,r}} \right)^2 \ \mu \left( r, X^{Y_{\lfloor r \rfloor_h}}_{\lfloor r \rfloor_h,r} \right)^2 \right] \\
    \leqslant \E \left[ f' \left( X_{\lceil r \rceil_h,T}^{X^{Y_{\lfloor r \rfloor_h}}_{\lfloor r \rfloor_h,r}} \right)^\frac{p}{q} \right]^\frac{2q}{p} \E \left[ \left( \partial_x X_{\lceil r \rceil_h,T}^{X^{Y_{\lfloor r \rfloor_h}}_{\lfloor r \rfloor_h,r}} \right)^\frac{2p}{p-2(q+1)} \right]^\frac{p-2(q+1)}{p} \E\left[  \mu \left( r, X^{Y_{\lfloor r \rfloor_h}}_{\lfloor r \rfloor_h,r} \right)^p \right]^\frac{2}{p}
\end{multline}
and 
\begin{equation} \label{eq:drift_estimate_holder_Y}
    \E\left[ f' \left( X_{\lceil r \rceil_h,T}^{Y_r} \right)^2 \left( \partial_x X_{\lceil r \rceil_h,T}^{Y_r} \right)^2 \ A_r^2   \right] \leq \E\left[  f' \left( X_{\lceil r \rceil_h,T}^{Y_r} \right)^\frac{p}{q} \right]^\frac{q}{2p} \E \left[  \left( \partial_x X_{\lceil r \rceil_h,T}^{Y_r} \right)^{\frac{p}{p-2(q+1)}} \right]^\frac{p-2(q+1)}{2p} \E \left[ A_r^p \right]^\frac{2}{p}. 
\end{equation}
The goal is now to find an upper-bound of each of the six expectations in the right hand side of \eqref{eq:drift_estimate_holder_X} and \eqref{eq:drift_estimate_holder_Y}, uniformly in $h$ and in $r$. We detail only the first one, since the others works similarly. Using \eqref{hyp:testfunction} we have 
\begin{equation*}
    \E\left[ f' \left( X_{\lceil r \rceil_h,T}^{X^{Y_{\lfloor r \rfloor_h}}_{\lfloor r \rfloor_h,r}} \right)^\frac{p}{q} \right] \leq 2^{p-1} C_f^p \left(1 + \E \left[ \left| X_{\lceil r \rceil_h,T}^{X^{Y_{\lfloor r \rfloor_h}}_{\lfloor r \rfloor_h,r}} \right|^p \right] \right).
\end{equation*}
Let $\Gamma \coloneqq \cup_{m \in \N \backslash \{ 0 \} } \{ jT/m, \ j=0,\ldots,m\}$. Since $\Gamma$ is countable, we may assume that $r \in [0,T] \backslash \Gamma$ (that makes no difference in the Lebesgue integral). 
This allows us to consider $i \in \{0, \ldots, n-1 \}$ such that $r \in (ih,(i+1)h)$, so that we have $\lceil r \rceil_h = (i+1)h$ and $\lfloor r \rfloor_h = ih$. Then, conditioning by $\mathbb{F}_{[ih,(i+1)h]}$, we have 
\begin{equation*}
     \E \left[ \left| X_{\lceil r \rceil_h,T}^{X^{Y_{\lfloor r \rfloor_h}}_{\lfloor r \rfloor_h,r}} \right|^p \right] = \E \left[ \left| X_{(i+1)h,T}^{X^{Y_{ih}}_{ih,r}} \right|^p \right] = \E \left[ \phi \left(  X^{Y_{ih}}_{ih,r} \right) \right],
\end{equation*}
where for every $x \in \R$, 
\[ \phi(x):= \E \left[ \left| X_{(i+1)h,T}^x \right|^p \right] \leqslant \sup_{0 \leqslant s \leqslant t \leqslant T } \E \left[ \left| X_{s,t}^x \right|^p \right]. \] 
It thus follows from \eqref{hyp:MomentsX} that 
\begin{equation*}
    \E \left[ \left| X_{\lceil r \rceil_h,T}^{X^{Y_{\lfloor r \rfloor_h}}_{\lfloor r \rfloor_h,r}} \right|^p \right] \leqslant C_{X,p} \left( 1+\E \left[ \left| X_{ih,r}^{Y_{ih}} \right|^{\mx} \right] \right). 
\end{equation*}
Then, by conditioning by $\mathbb{F}_{[0,ih]}$ and using \eqref{hyp:MomentsX} again, we obtain that
\[ \E \left[ \left| X_{ih,r}^{Y_{ih}} \right|^{\mx} \right] \leqslant C_{X,\mx} \left( 1 + \E \left[ \big| Y_{ih} \big|^{\mxx} \right] \right) \leqslant C_{X,\mx} \left( 1 + \sup_{u \in [0,T]} \E \left[ \big| Y_u \big|^{\mxx} \right] \right). \]
Finally, we conclude thanks to \eqref{hyp:MomentsY} that $\sup_{u \in [0,T]} \E[|Y_u|^{\mxx}] < + \infty$. Hence  
\begin{equation*}
    \sup_{h \in (0, T)} \sup_{r \in [0, T] \backslash \Gamma} \E \left[ \left| f' \left( X_{\lceil r \rceil_h,T}^{X^{Y_{\lfloor r \rfloor_h}}_{\lfloor r \rfloor_h,r}} \right) \right|^{p/q} \right]^{2q/p} < + \infty.
\end{equation*}
We proceed with similar estimates for the other expectation terms appearing in the right hand side of \eqref{eq:drift_estimate_holder_X} and \eqref{eq:drift_estimate_holder_Y}, by using \eqref{hyp:coeffs}, \eqref{hyp:coeffs}, \eqref{hyp:MomentsX}, \eqref{hyp:MomentsX} and \eqref{hyp:MomentsY}. The result follows afterwards. 

\subsection{Proof of Lemma \ref{lm:CroissancePolynomialeDeFhEnZ}}\label{subsec:Borne_integrale_dNu(z)dr}
Let $h \in (0, T)$, $r \in [0, T]$ and $z \in \R \backslash \{ 0 \}$. Using Taylor formula with integral leftover, we have
\begin{multline*}
    f \left( X_{\lceil r \rceil_h,T}^{X^{Y_{\lfloor r \rfloor_h}}_{\lfloor r \rfloor_h,{r}} + \sigma \left( r,X^{Y_{\lfloor r \rfloor_h}}_{\lfloor r \rfloor_h,{r}},z \right)} \right) - f \left( X_{\lceil r \rceil_h,T}^{X^{Y_{\lfloor r \rfloor_h}}_{\lfloor r \rfloor_h,{r}} } \right) - \sigma \left( r,X^{Y_{\lfloor r \rfloor_h}}_{\lfloor r \rfloor_h,{r}},z \right) f' \left(  X_{\lceil r \rceil_h,T}^{X^{Y_{\lfloor r \rfloor_h}}_{\lfloor r \rfloor_h,{r}}} \right) \partial_x X_{\lceil r \rceil_h,T}^{X^{Y_{\lfloor r \rfloor_h}}_{\lfloor r \rfloor_h,r}} \\
    =  \sigma \left( r,X^{Y_{\lfloor r \rfloor_h}}_{\lfloor r \rfloor_h,{r}},z \right)^2 \int_0^1 \Bigg\{(1 - \lambda) f''\left( X_{\lceil r \rceil_h,T}^{X^{Y_{\lfloor r \rfloor_h}}_{\lfloor r \rfloor_h,{r}} + \lambda \ \sigma \left( r,X^{Y_{\lfloor r \rfloor_h}}_{\lfloor r \rfloor_h,{r}},z \right)} \right) \left( \partial_x X_{\lceil r \rceil_h,T}^{X^{Y_{\lfloor r \rfloor_h}}_{\lfloor r \rfloor_h,r}+ \lambda \ \sigma \left( r,X^{Y_{\lfloor r \rfloor_h}}_{\lfloor r \rfloor_h,{r}},z \right)} \right)^2 \\ 
    + (1 - \lambda) f' \left( X_{\lceil r \rceil_h,T}^{X^{Y_{\lfloor r \rfloor_h}}_{\lfloor r \rfloor_h,{r}} +\lambda \ \sigma \left( r,X^{Y_{\lfloor r \rfloor_h}}_{\lfloor r \rfloor_h,{r}},z \right)} \right) \partial_x^2 X_{\lceil r \rceil_h,T}^{X^{Y_{\lfloor r \rfloor_h}}_{\lfloor r \rfloor_h,r} + \lambda \ \sigma \left( r,X^{Y_{\lfloor r \rfloor_h}}_{\lfloor r \rfloor_h,{r}},z \right)} \Bigg\} \ \d \lambda,
\end{multline*}
and
\begin{align*}
     & f \left( X_{\lceil r \rceil_h,T}^{Y_{r}+B_{r,z}} \right) - f \left( X_{\lceil r \rceil_h,T}^{Y_{r}} \right)  - B_{r,z} f' \left( X_{\lceil r \rceil_h,T}^{Y_r} \right) \partial_x X_{\lceil r \rceil_h,T}^{Y_r} \\
     ={} & B_{r,z}^2 \int_0^1  \bigg\{(1 - \lambda) f''\left( X_{\lceil r \rceil_h,T}^{Y_{r} + \lambda B_{r,z}} \right) \left( \partial_x X_{\lceil r \rceil_h,T}^{Y_{r} + \lambda B_{r,z}} \right)^2  + (1 - \lambda) f' \left( X_{\lceil r \rceil_h,T}^{Y_{r} + \lambda B_{r,z}} \right) \partial_{x}^2 X_{\lceil r \rceil_h,T}^{Y_{r} + \lambda B_{r,z}} \bigg\} \ \d \lambda.
\end{align*}
Applying Minkovski's integral inequality, we get (with the abuse of notation $\norm{f(r)}{L^2 (\P \otimes \lambda_{[0, T]})}$, for $f: \Omega \times [0, T] \longrightarrow \R$, to mean $\norm{f}{L^2 (\P \otimes \lambda_{[0, T]})}$) 
\begin{align*}
    & \Bigg\| \int_0^1 \sigma \left( r,X^{Y_{\lfloor r \rfloor_h}}_{\lfloor r \rfloor_h,{r}},z \right)^2 \Bigg\{ (1 - \lambda) f''\left( X_{\lceil r \rceil_h,T}^{X^{Y_{\lfloor r \rfloor_h}}_{\lfloor r \rfloor_h,{r}} + \lambda \ \sigma \left( r,X^{Y_{\lfloor r \rfloor_h}}_{\lfloor r \rfloor_h,{r}},z \right)} \right) \left( \partial_x X_{\lceil r \rceil_h,T}^{X^{Y_{\lfloor r \rfloor_h}}_{\lfloor r \rfloor_h,r}+\lambda \sigma \left( r,X^{Y_{\lfloor r \rfloor_h}}_{\lfloor r \rfloor_h,{r}},z \right)} \right)^2 
    \\ & + (1 - \lambda) f' \left( X_{\lceil r \rceil_h,T}^{X^{Y_{\lfloor r \rfloor_h}}_{\lfloor r \rfloor_h,{r}} + \lambda \sigma \left( r,X^{Y_{\lfloor r \rfloor_h}}_{\lfloor r \rfloor_h,{r}},z \right)} \right) \partial_x^2 X_{\lceil r \rceil_h,T}^{X^{Y_{\lfloor r \rfloor_h}}_{\lfloor r \rfloor_h,r}+\lambda \sigma \left( r,X^{Y_{\lfloor r \rfloor_h}}_{\lfloor r \rfloor_h,{r}},z \right)} \Bigg\} \ \d \lambda   \Bigg\|_{L^2(\P \otimes \lambda_{[0,T]})}  
    \\ \leqslant{} &  \int_0^1  \Bigg\| \sigma \left( r,X^{Y_{\lfloor r \rfloor_h}}_{\lfloor r \rfloor_h,{r}},z \right)^2 \Bigg\{ (1-\lambda) f''\left( X_{\lceil r \rceil_h,T}^{X^{Y_{\lfloor r \rfloor_h}}_{\lfloor r \rfloor_h,{r}} + \lambda \ \sigma \left( r,X^{Y_{\lfloor r \rfloor_h}}_{\lfloor r \rfloor_h,{r}},z \right)} \right) \left( \partial_x X_{\lceil r \rceil_h,T}^{X^{Y_{\lfloor r \rfloor_h}}_{\lfloor r \rfloor_h,r}+\lambda \sigma \left( r,X^{Y_{\lfloor r \rfloor_h}}_{\lfloor r \rfloor_h,{r}},z \right)} \right)^2 \\ 
    & +  (1-\lambda) f' \left( X_{\lceil r \rceil_h,T}^{X^{Y_{\lfloor r \rfloor_h}}_{\lfloor r \rfloor_h,{r}} + \lambda \sigma \left( r,X^{Y_{\lfloor r \rfloor_h}}_{\lfloor r \rfloor_h,{r}},z \right)} \right) \partial_x^2 X_{\lceil r \rceil_h,T}^{X^{Y_{\lfloor r \rfloor_h}}_{\lfloor r \rfloor_h,r}+\lambda \sigma \left( r,X^{Y_{\lfloor r \rfloor_h}}_{\lfloor r \rfloor_h,{r}},z \right)} \Bigg\} \Bigg\|_{L^2(\P \otimes \lambda_{[0,T]})} \ \d \lambda.
\end{align*}
Similarly, we get
\begin{align*}
    & \left\| B_{r,z}^2 \int_0^1  \bigg\{ (1-\lambda) f''\left( X_{\lceil r \rceil_h,T}^{Y_{r} + \lambda B_{r,z}} \right) \left( \partial_x X_{\lceil r \rceil_h,T}^{Y_{r} + \lambda B_{r,z}} \right)^2 \right. \\
    & \hspace{10em} \left. {} + (1-\lambda) f' \left( X_{\lceil r \rceil_h,T}^{Y_{r} + \lambda B_{r,z}} \right) \partial_{x}^2 X_{\lceil r \rceil_h,T}^{Y_{r} + \lambda B_{r,z}} \bigg\} \ \d \lambda \right\|_{L^2(\P \otimes \lambda_{[0,T]})} \\
    \leqslant{} & \int_0^1 \left\| B_{r,z}^2 \left\{ (1 - \lambda) f'' \left( X_{\lceil r \rceil_h,T}^{Y_{r} + \lambda B_{r,z}} \right) \left( \partial_x X_{\lceil r \rceil_h,T}^{Y_{r} + \lambda B_{r,z}} \right)^2 \right. \right. \\
    & \hspace{10em} \left. {}+ (1-\lambda) f' \left( X_{\lceil r \rceil_h,T}^{Y_{r} + \lambda B_{r,z}} \right) \partial_{x}^2 X_{\lceil r \rceil_h,T}^{Y_{r} + \lambda B_{r,z}} \bigg\} \right\|_{L^2(\P \otimes \lambda_{[0,T]})} \d \lambda.
\end{align*}
Let $\lambda \in [0, 1]$. By Hölder's inequality, 
\begin{multline*}
    \Bigg\| \sigma \left( r,X^{Y_{\lfloor r \rfloor_h}}_{\lfloor r \rfloor_h,{r}},z \right)^2 \Bigg\{ \lambda f''\left( X_{\lceil r \rceil_h,T}^{X^{Y_{\lfloor r \rfloor_h}}_{\lfloor r \rfloor_h,{r}} + \lambda  \sigma \left( r,X^{Y_{\lfloor r \rfloor_h}}_{\lfloor r \rfloor_h,{r}},z \right)} \right) \left( \partial_x X_{\lceil r \rceil_h,T}^{X^{Y_{\lfloor r \rfloor_h}}_{\lfloor r \rfloor_h,r}+\lambda \sigma \left( r,X^{Y_{\lfloor r \rfloor_h}}_{\lfloor r \rfloor_h,{r}},z \right)} \right)^2 \\ 
    +  \lambda f' \left( X_{\lceil r \rceil_h,T}^{X^{Y_{\lfloor r \rfloor_h}}_{\lfloor r \rfloor_h,{r}} + \lambda \sigma \left( r,X^{Y_{\lfloor r \rfloor_h}}_{\lfloor r \rfloor_h,{r}},z \right)} \right) \partial_x^2 X_{\lceil r \rceil_h,T}^{X^{Y_{\lfloor r \rfloor_h}}_{\lfloor r \rfloor_h,r}+\lambda \sigma \left( r,X^{Y_{\lfloor r \rfloor_h}}_{\lfloor r \rfloor_h,{r}},z \right)} \Bigg\} \Bigg\|_{L^2(\P \otimes \lambda_{[0,T]})}\\
    \leqslant \Bigg\| \sigma \left( r,X^{Y_{\lfloor r \rfloor_h}}_{\lfloor r \rfloor_h,{r}},z \right)^2 \Bigg\|_{L^{p/2}(\P \otimes \lambda_{[0,T]})} \ \Bigg\| \lambda f''\left( X_{\lceil r \rceil_h,T}^{X^{Y_{\lfloor r \rfloor_h}}_{\lfloor r \rfloor_h,{r}} + \lambda \ \sigma \left( r,X^{Y_{\lfloor r \rfloor_h}}_{\lfloor r \rfloor_h,{r}},z \right)} \right) \left( \partial_x X_{\lceil r \rceil_h,T}^{X^{Y_{\lfloor r \rfloor_h}}_{\lfloor r \rfloor_h,r}+\lambda \sigma \left( r,X^{Y_{\lfloor r \rfloor_h}}_{\lfloor r \rfloor_h,{r}},z \right)} \right)^2 \\ 
    +  \lambda f' \left( X_{\lceil r \rceil_h,T}^{X^{Y_{\lfloor r \rfloor_h}}_{\lfloor r \rfloor_h,{r}} + \lambda \sigma \left( r,X^{Y_{\lfloor r \rfloor_h}}_{\lfloor r \rfloor_h,{r}},z \right)} \right) \partial_x^2 X_{\lceil r \rceil_h,T}^{X^{Y_{\lfloor r \rfloor_h}}_{\lfloor r \rfloor_h,r} + \lambda \sigma \left( r,X^{Y_{\lfloor r \rfloor_h}}_{\lfloor r \rfloor_h,{r}},z \right)} \Bigg\|_{L^{\frac{2p}{p-4}}(\P \otimes \lambda_{[0,T]})}. 
\end{multline*}
One one hand we have, due to \eqref{hyp:coeffs}, 
\begin{equation*}
    \left| \sigma \left( r,X^{Y_{\lfloor r \rfloor_h}}_{\lfloor r \rfloor_h,{r}},z \right) \right| \leqslant C_\sigma |z|^k \left( 1 + \left| X^{Y_{\lfloor r \rfloor_h}}_{\lfloor r \rfloor_h,{r}} \right|^{\ms} \right),
\end{equation*}
leading to 
\begin{equation*}
    \Bigg\| \sigma \left( r,X^{Y_{\lfloor r \rfloor_h}}_{\lfloor r \rfloor_h,{r}},z \right)^2 \Bigg\|_{L^{p/2}(\P \otimes \lambda_{[0,T]})} \leqslant 2 C_\sigma^2 |z|^{2k} \ \left( T + \int_0^T \E \left[ \left| X^{Y_{\lfloor r \rfloor_h}}_{\lfloor r \rfloor_h,{r}} \right|^{\ms p} \right]^{2/p} \ \d r \right).
\end{equation*}
Let $i \in \{ 0, \ldots , n-1 \}$ such that $r \in (ih , (i+1)h)$. Conditioning by $\mathbf{F}_{i h}$, we get by using \eqref{hyp:MomentsX}
\begin{equation*}
    \E \left[ \left| X^{Y_{\lfloor r \rfloor_h}}_{\lfloor r \rfloor_h,{r}} \right|^{\ms p} \right] = \E \left[ \left| X^{Y_{i h}}_{i h, r} \right|^{\ms p} \right] \leqslant C_{X, p, \sigma} \left( 1 + \E \left[ \left| Y_{ih} \right|^{\msx} \right] \right) \leqslant C_{X, p, \sigma} \left( 1 + \sup_{r \in [0,T]} \E \left[ |Y_r|^{\msx} \right] \right), 
\end{equation*}
which gives
\begin{equation*}
    \Bigg\| \sigma \left( r,X^{Y_{\lfloor r \rfloor_h}}_{\lfloor r \rfloor_h,{r}},z \right)^2 \Bigg\|_{L^{p/2}(\P \otimes \lambda_{[0,T]})} \leqslant C_0 |z|^{2k}
\end{equation*}
where 
\[ C_0 \defeq 2 C_\sigma^2 \left( T + T  C_{X, p, \sigma} \left( 1 + \sup_{r \in [0,T]} \E \left[ |Y_r|^{\msx} \right] \right)^{\frac{2}{p}} \right) < + \infty. \]

Then applying Hölder's inequality again,
\begin{align*}
    & \left\|(1 - \lambda) f''\left( X_{\lceil r \rceil_h,T}^{X^{Y_{\lfloor r \rfloor_h}}_{\lfloor r \rfloor_h,{r}} + \lambda \ \sigma \left( r,X^{Y_{\lfloor r \rfloor_h}}_{\lfloor r \rfloor_h,{r}},z \right)} \right) \left( \partial_x X_{\lceil r \rceil_h,T}^{X^{Y_{\lfloor r \rfloor_h}}_{\lfloor r \rfloor_h,r}+\lambda \sigma \left( r,X^{Y_{\lfloor r \rfloor_h}}_{\lfloor r \rfloor_h,{r}},z \right)} \right)^2 \right\|_{L^{\frac{2p}{p-4}}(\P \otimes \lambda_{[0,T]})} \\ 
    \leqslant{} &  \left\|  f''\left( X_{\lceil r \rceil_h,T}^{X^{Y_{\lfloor r \rfloor_h}}_{\lfloor r \rfloor_h,{r}} + \lambda \sigma \left( r,X^{Y_{\lfloor r \rfloor_h}}_{\lfloor r \rfloor_h,{r}},z \right)} \right) \right\|_{L^{\frac{p}{q}}(\P \otimes \lambda_{[0,T]})}\times \left\|  \partial_x X_{\lceil r \rceil_h,T}^{X^{Y_{\lfloor r \rfloor_h}}_{\lfloor r \rfloor_h,r}+\lambda \sigma \left( r,X^{Y_{\lfloor r \rfloor_h}}_{\lfloor r \rfloor_h,{r}},z \right)} \right\|_{L^{\frac{4p}{p-2(q+2)}}(\P \otimes \lambda_{[0,T]})}^2
\end{align*}
and
\begin{align*}
    & \left\| (1 - \lambda) f' \left( X_{\lceil r \rceil_h,T}^{X^{Y_{\lfloor r \rfloor_h}}_{\lfloor r \rfloor_h,{r}} + \lambda \sigma \left( r,X^{Y_{\lfloor r \rfloor_h}}_{\lfloor r \rfloor_h,{r}},z \right)} \right) \partial_x^2 X_{\lceil r \rceil_h,T}^{X^{Y_{\lfloor r \rfloor_h}}_{\lfloor r \rfloor_h,r}+\lambda \sigma \left( r,X^{Y_{\lfloor r \rfloor_h}}_{\lfloor r \rfloor_h,{r}},z \right)} \right\|_{L^{\frac{2p}{p-4}}(\P \otimes \lambda_{[0,T]})} \\ 
    \leqslant{} & \left\|  f' \left( X_{\lceil r \rceil_h,T}^{X^{Y_{\lfloor r \rfloor_h}}_{\lfloor r \rfloor_h,{r}} + \lambda \sigma \left( r,X^{Y_{\lfloor r \rfloor_h}}_{\lfloor r \rfloor_h,{r}},z \right)} \right) \right\|_{L^{\frac{p}{q}}(\P \otimes \lambda_{[0,T]})}\times \left\|   \partial_x^2 X_{\lceil r \rceil_h,T}^{X^{Y_{\lfloor r \rfloor_h}}_{\lfloor r \rfloor_h,r}+\lambda \sigma \left( r,X^{Y_{\lfloor r \rfloor_h}}_{\lfloor r \rfloor_h,{r}},z \right)}\right\|_{L^{\frac{2p}{p-2(q+2)}}(\P \otimes \lambda_{[0,T]})}.
\end{align*}
Then combining \eqref{hyp:testfunction}, \eqref{hyp:MomentsX}, \eqref{hyp:MomentsX}, \eqref{hyp:MomentsX}, \eqref{hyp:coeffs}, \eqref{hyp:MomentsY} with the conditioning techniques already employed in Section \eqref{subsec:Borne_integrale_dr}, we obtain the existence of constants $C_i > 0$, $1 \leqslant i \leqslant 4$ such that:
\begin{align*}
    \left\|  f''\left( X_{\lceil r \rceil_h,T}^{X^{Y_{\lfloor r \rfloor_h}}_{\lfloor r \rfloor_h,{r}} + \lambda \ \sigma \left( r,X^{Y_{\lfloor r \rfloor_h}}_{\lfloor r \rfloor_h,{r}},z \right)} \right) \right\|_{L^{\frac{p}{q}}(\P \otimes \lambda_{[0,T]})} & \leqslant C_1 \left(1 + |z|^{ k \ \mx \frac{q}{p}} \right), \\
    \left\|  \partial_x X_{\lceil r \rceil_h,T}^{X^{Y_{\lfloor r \rfloor_h}}_{\lfloor r \rfloor_h,r}+\lambda \sigma \left( r,X^{Y_{\lfloor r \rfloor_h}}_{\lfloor r \rfloor_h,{r}},z \right)} \right\|_{L^{\frac{4p}{p-2(q+2)}}(\P \otimes \lambda_{[0,T]})}^2 & \leqslant C_2 \left( 1 +  |z|^{k \ \mdx \frac{p-2(q+1)}{2p} } \right), \\
    \left\|  f' \left( X_{\lceil r \rceil_h,T}^{X^{Y_{\lfloor r \rfloor_h}}_{\lfloor r \rfloor_h,{r}} + \lambda \sigma \left( r,X^{Y_{\lfloor r \rfloor_h}}_{\lfloor r \rfloor_h,{r}},z \right)} \right) \right\|_{L^{\frac{p}{q}}(\P \otimes \lambda_{[0,T]})} & \leqslant C_3 \left( 1 +  |z|^{ k \ \mx \frac{q}{p}} \right), \\
    \left\|   \partial_x^2 X_{\lceil r \rceil_h,T}^{X^{Y_{\lfloor r \rfloor_h}}_{\lfloor r \rfloor_h,r}+\lambda \sigma \left( r,X^{Y_{\lfloor r \rfloor_h}}_{\lfloor r \rfloor_h,{r}},z \right)}\right\|_{L^{\frac{2p}{p-2(q+2)}}(\P \otimes \lambda_{[0,T]})} & \leqslant C_4 \left(1 +  |z|^{k \ \mddx \frac{p-2(q+1)}{2p}} \right). 
\end{align*}
In the whole, we have shown that
\begin{align*}
    & \left\| \sigma \left( r,X^{Y_{\lfloor r \rfloor_h}}_{\lfloor r \rfloor_h,{r}},z \right)^2 \left\{ \lambda f''\left( X_{\lceil r \rceil_h,T}^{X^{Y_{\lfloor r \rfloor_h}}_{\lfloor r \rfloor_h,{r}} + \lambda \ \sigma \left( r,X^{Y_{\lfloor r \rfloor_h}}_{\lfloor r \rfloor_h,{r}},z \right)} \right) \left( \partial_x X_{\lceil r \rceil_h,T}^{X^{Y_{\lfloor r \rfloor_h}}_{\lfloor r \rfloor_h,r}+\lambda \sigma \left( r,X^{Y_{\lfloor r \rfloor_h}}_{\lfloor r \rfloor_h,{r}},z \right)} \right)^2 \right. \right.
    \\ & \left. \left. {} + \lambda f' \left( X_{\lceil r \rceil_h,T}^{X^{Y_{\lfloor r \rfloor_h}}_{\lfloor r \rfloor_h,{r}} + \lambda \sigma \left( r,X^{Y_{\lfloor r \rfloor_h}}_{\lfloor r \rfloor_h,{r}},z \right)} \right) \partial_x^2 X_{\lceil r \rceil_h,T}^{X^{Y_{\lfloor r \rfloor_h}}_{\lfloor r \rfloor_h,r}+\lambda \sigma \left( r,X^{Y_{\lfloor r \rfloor_h}}_{\lfloor r \rfloor_h,{r}},z \right)} \right\} \right\|_{L^2(\P \otimes \lambda_{[0,T]})} \\
    \leqslant{} & C_5 |z|^{2k} \left( 1 + |z|^{ k \ \mx \frac{q}{p}} \right) \left( 1 + |z|^{ k \ \mdx \frac{p-2(q+1)}{2p} } \right) \left(1 + |z|^{
     k \ \mx \frac{q}{p}} \right) \left( 1 + |z|^{ k \ \mddx \frac{p-2(q+1)}{2p}} \right),
\end{align*}
with $C_5 = C_0C_1C_2C_3C_4$. With similar arguments, replacing the use of \eqref{hyp:coeffs} by those of \eqref{hyp:coeffs}, we obtain
\begin{align*}
    & \left\| B_{r,z}^2 \bigg\{ \lambda f''( X_{\lceil r \rceil_h,T}^{Y_{r} + \lambda B_{r,z}}) \left( \partial_x X_{\lceil r \rceil_h,T}^{Y_{r} + \lambda B_{r,z}} \right)^2  + \lambda f'(X_{\lceil r \rceil_h,T}^{Y_{r} + \lambda B_{r,z}}) \partial_{x}^2 X_{\lceil r \rceil_h,T}^{Y_{r} + \lambda B_{r,z}} \bigg\} \right\|_{L^2(\P \otimes \lambda_{[0,T]})} \\
    \leqslant{}  & C_5' |z|^{2k} \left( 1 + |z|^{ k \ \mx \frac{q}{p}} \right) \left( 1 + |z|^{ k \ \mdx \frac{p-2(q+1)}{2p}} \right) \left( 1 + |z|^{
     k \ \mx \frac{q}{p}} \right) \left( 1 + |z|^{ k \ \mddx \frac{p-2(q+1)}{2p} } \right),
\end{align*}
for some constant $C'_5 > 0$. It follows, due to \eqref{hyp:mesure_levy}, that there exists $\psi \in L^1(\nu) \cap L^2(\nu)$ such that $\sup_{h \in (0, T)} \| F_h(\cdot,z) \|_{L^2(\P \otimes \lambda_{[0,T]})} \leqslant \psi(z)$.

Finally since 
\begin{align*}
      & F(r, z) \\
      ={} & f \left( X_{r,T}^{Y_r + \sigma \left( r,Y_r,z \right)} \right) - f \left( X_{r,T}^{Y_{r}+B_{r,z}} \right)
      +  B_{r,z} f' \left( X_{r,T}^{Y_r} \right) \partial_x X_{r,T}^{Y_r} - \sigma \left( r,Y_r,z \right) f'\left( X_{r,T}^{Y_r} \right) \partial_x X_{r,T}^{Y_r} \\ 
      ={} &  \left( f \left( X_{r,T}^{Y_r + \sigma \left( r,Y_r,z \right)} \right) - f(X_{r,T}^{Y_r}) -  \sigma \left( r,Y_r,z \right) f'\left( X_{r,T}^{Y_r} \right) \partial_x X_{r,T}^{Y_r} \right) \\
      & \qquad - \left( f \left( X_{r,T}^{Y_{r}+B_{r,z}} \right) - f \left( X_{r,T}^{Y_{r}} \right) -  B_{r,z} f' \left( X_{r,T}^{Y_r} \right) \partial_x X_{r,T}^{Y_r} \right),
\end{align*}
the same arguments apply also for $F$ and we obtain $ \| F (\cdot,z) \|_{L^2(\P \otimes \lambda_{[0,T]})} \leqslant \psi(z)$.

\subsection{Proof of the estimate \eqref{eq:Borne_Moment_Int_3eme_Terme}}\label{subsec:Borne_integrale_Sto}

The idea is to write the stochastic integral by using equality \eqref{eq:FormuleIntermediaireH}. Two of three terms are already bounded, by using \eqref{eq:borne_integrale_dr} and \eqref{eq:BorneMomentIntegraleNudZdR}. Moreover, we have by \eqref{hyp:testfunction}, that for every $x \in \R$, 
\[ |f(x)| \leqslant C_f (1+|x|^q) (1+|x|) \leqslant C_f C_q (1+|x|)^q, \]
where $C_q = \max\{ 2^{1-q}, 1 \}$, by distinguishing the cases $q<1$ or $q\geqslant1$ and using concavity or convexity of $x \mapsto x^q$. Hence, we get by triangle inequality
\[ \begin{split}
    \norm{f \left( X_{0, T}^{Y_0} \right) - f(Y_T)}{L^2 (\P)} & \leqslant C_f C_q \left( \norm{\left( 1 + \left| X_{0, T}^{Y_0} \right| \right)^{q+1}}{L^2 (\P)} + \norm{\Big( 1 + \left| Y_T \right| \Big)^{q+1}}{L^2 (\P)} \right) \\
    & = C_f C_q \left( \norm{ 1 + \left| X_{0, T}^{Y_0} \right| }{L^{2(q+1)} (\P)}^{q+1} + \norm{ 1 + \left| Y_T \right| }{L^{2(q+1)} (\P)}^{q+1} \right).
\end{split}  \]
Since $2q+2 < p$, we have 
\[ \begin{split}
    \norm{f \left( X_{0, T}^{Y_0} \right) - f(Y_T)}{L^2 (\P)} & \leqslant C_f C_q \left( \norm{ 1 + \left| X_{0, T}^{Y_0} \right| }{L^{p} (\P)}^{q+1} + \norm{ 1 + \left| Y_T \right| }{L^{p} (\P)}^{q+1} \right) \\
    & \leqslant C_f C_q \left( \left( 1 + \norm{  X_{0, T}^{Y_0} }{L^{p} (\P)} \right)^{q+1} + \left( 1 + \norm{  Y_T  }{L^{p} (\P)} \right)^{q+1} \right).
\end{split}  \]
Now, we use the same conditional technique used in subsection \ref{proof:controle_l2_second_terme}, to conclude with \eqref{hyp:MomentsX} that
\[ \E \left[ \left| X_{0, T}^{Y_0} \right|^p \right] \leqslant C_{X, p} \left( 1 + \E \left[ |Y_0|^{m_X^{(p)}} \right] \right) \leqslant C_{X, p} \left( 1 + \sup_{0 \leqslant r \leqslant T} \E \left[ |Y_r|^{m_X^{(p)}} \right] \right), \]
which is finite by hypothesis \eqref{hyp:MomentsY}. Similarly, we have 
\[ \E \left[ \left| Y_T \right|^p \right] = \E \left[ \left| X_{T, T}^{Y_T} \right|^p \right] \leqslant C_{X, p} \left( 1 + \sup_{0 \leqslant r \leqslant T} \E \left[ |Y_r|^{m_X^{(p)}} \right] \right). \]
We conclude that
\[  \norm{f \left( X_{0, T}^{Y_0} \right) - f(Y_T)}{L^2 (\P)} \leqslant 2 C_f C_q C_{X, p} \left( 1 + \sup_{0 \leqslant r \leqslant T} \E \left[ |Y_r|^{m_X^{(p)}} \right] \right)^{q+1} < \infty. \]
Then, by \eqref{eq:FormuleIntermediaireH}, \eqref{eq:borne_integrale_dr}, \eqref{eq:BorneMomentIntegraleNudZdR}, we conclude that 
\[ \begin{split}
    & \sup_{h \in (0, T)}   \E \left[ \left( \int_0^T \int_{\R \backslash \{ 0 \}}  \left\{ f \left( X_{\lceil r \rceil_h,T}^{X^{Y_{\lfloor r \rfloor_h}}_{\lfloor r \rfloor_h,{r}} + \sigma \left( r,X^{Y_{\lfloor r \rfloor_h}}_{\lfloor r \rfloor_h,{r}},z \right)} \right) -  f \Big( X_{\lceil r \rceil_h,T}^{Y_{r}+B_{r,z}} \Big) \right. \right. \right.    \\
    & \hspace{8em} \left. \left. \left. +{}  f \left( X_{\lceil r \rceil_h,T}^{X^{Y_{\lfloor r \rfloor_h}}_{\lfloor r \rfloor_h,{r}}} \right)  - f \Big( X_{\lceil r \rceil_h,T}^{Y_{r}} \Big)  \right\}  \ \tilde{N} (\delta r, \delta z) \right)^2 \right]^{1/2}  \\
    \leqslant {}  & \norm{f \left( X_{0, T}^{Y_0} \right) - f(Y_T)}{L^2 (\P)} + \sup_{h \in (0, T)} \norm{F_h}{L^2(\P \otimes \lambda_{[0, T]} \otimes \nu)}  \\
    & + \sup_{h \in (0, T)} \E \left[  \int_{0}^{T} \left|  f' \left( X_{\lceil r \rceil_h,T}^{X^{Y_{\lfloor r \rfloor_h}}_{\lfloor r \rfloor_h,r}} \right) \partial_x X_{\lceil r \rceil_h,T}^{X^{Y_{\lfloor r \rfloor_h}}_{\lfloor r \rfloor_h,r}}  \ \mu \left( r, X^{Y_{\lfloor r \rfloor_h}}_{\lfloor r \rfloor_h,r} \right) -  f' \left( X_{\lceil r \rceil_h,T}^{Y_r} \right) \partial_x X_{\lceil r \rceil_h,T}^{Y_r} A_r \right|^{2} \ \d r  \right]^{1/2} \\
    < {} & \infty. 
\end{split} \]
This proves \eqref{eq:Borne_Moment_Int_3eme_Terme}.

\subsection{Proof of Lemma \ref{lm:L2ContinuityExampleProcess}} \label{subsec:ProofL2Continuity}
We start by proving the $L^2$-continuity for the flow process $X$. Let $t \in [0,T]$, $(x,h) \in \R^2$, $s \in [0,t)$ and $\delta \in [0,t-s)$. By definition, one has
\begin{multline*}
  \| X_{s + \delta, t}^{x + h} - X_{s, t}^x \|_{L^2 (\Omega)} =  \Bigg\| x
  + h + \int_{s + \delta}^t b (X_{s + \delta, r}^{x + h}) \d r + \int_{s +
  \delta}^t \sigma (X_{s + \delta, r^-}^{x + h}) \d L_r  \\ - \bigg( x + \int_s^t b
  (X_{s, r}^x) \d r + \int_s^t \sigma (X_{s, r^-}^x) \d L_r \bigg)
  \Bigg\|_{L^2 (\Omega)}\\
   \leq  h + \left\| \int_{s + \delta}^t (b (X_{s + \delta, r}^{x + h}) - b
  (X_{s, r}^x)) \d r \right\|_{L^2 (\Omega)} + \left\| \int_{s + \delta}^t (b
  (X_{s + \delta, r^-}^{x + h}) - b (X_{s, r^-}^x)) \d L_r \right\|_{L^2
  (\Omega)} \\ + \left\| \int_s^{s + \delta} b (X_{s, r}^x) \d r \right\|_{L^2
  (\Omega)} + \left\| \int_s^{s + \delta} \sigma (X_{s, r^-}^x) \d L_r
  \right\|_{L^2 (\Omega)} .
\end{multline*}
We apply Minkovski's integral inequality and Ito's isometry to the Lévy process $L$ with Lévy measure $\nu(\d z) = \frac{\d z}{|z|^{\alpha+1}}$, leading to
\begin{equation*}
\begin{aligned}
  \left\| \int_{s + \delta}^t (b (X_{s + \delta, r}^{x + h}) - b (X_{s, r}^x))
  \d r \right\|_{L^2 (\Omega)} & \leq  \int_{s + \delta}^t \| b (X_{s +
  \delta, r}^{x + h}) - b (X_{s, r}^x) \|_{L^2 (\Omega)} \d r, \\
  \left\| \int_s^{s + \delta} b (X_{s, r}^x) \d r \right\|_{L^2 (\Omega)} &
  \leq \int_s^{s + \delta} \| b (X_{s, r}^x) \|_{L^2 (\Omega)} \d r, \\
  \left\| \int_{s + \delta}^t (\sigma (X_{s + \delta, r^-}^{x + h}) - \sigma
  (X_{s, r^-}^x)) \d L_r \right\|_{L^2 (\Omega)} & \leq  \left( \int_{s +
  \delta}^t \int_{\mathbb{R}\backslash \{ 0 \}} \| \sigma (X_{s + \delta,
  r}^{x + h}) - \sigma (X_{s, r}^x) \|_{L^2 (\Omega)}^2 z \nu (\d z) \d r
  \right)^{\frac{1}{2}}, \\
  \left\| \int_s^{s + \delta} \sigma (X_{s, r^-}^x) \d L_r \right\|_{L^2
  (\Omega)} & \leq \left( \int_s^{s + \delta} \int_{\mathbb{R}\backslash \{
  0 \}} \| \sigma (X_{s, r}^x) \|_{L^2 (\Omega)}^2 z \nu (\d z) \d r
  \right)^{\frac{1}{2}} .
\end{aligned}
\end{equation*}
Since $b$ and $\sigma$ are of class $C^1$ with bounded derivative, they are in
particular Lipschitz. Denoting respectively by $L_b$ and $L_{\sigma}$ their
Lipschitz constant, we have
\begin{equation*}
  \| b (X_{s + \delta, r}^{x + h}) - b (X_{s, r}^x) \|_{L^2 (\Omega)} \leq 
  L_b  \| X_{s + \delta, t}^{x + h} - X_{s, t}^x \|_{L^2 (\Omega)}
\end{equation*}
and
\begin{equation*}
  \| \sigma (X_{s + \delta, r^-}^{x + h}) - \sigma (X_{s, r^-}^x) \|^2_{L^2
  (\Omega)} \leq  L_{\sigma}^2  \| X_{s + \delta, t}^{x + h} - X_{s, t}^x
  \|^2_{L^2 (\Omega)},
\end{equation*}
In addition, we may compute
\begin{equation*}
  \int_{\mathbb{R}\backslash \{ 0 \}} z \nu (\d z) =  2 \int_0^1 z^{-
  \alpha} \d z \text{ } = \text{ } \frac{2}{\alpha + 1} .
\end{equation*}
It follows that
\begin{align*}
  \| X_{s + \delta, t}^{x + h} - X_{s, t}^x \|_{L^2 (\Omega)}  \leq &h + \delta \| b \|_{L^\infty(\R)} + \sqrt{\delta} \sqrt{2 (\alpha+1)^{-1}} \| \sigma
  \|_{L^\infty(\R)} + \int_{s + \delta}^t L_b
  \| X_{s + \delta, r}^{x + h} - X_{s, r}^x \|_{L^2 (\Omega)} \d r\\
  &   + \left( \int_{s + \delta}^t \frac{2}{\alpha + 1} L_{\sigma}^2 \| X_{s
  + \delta, r}^{x + h} - X_{s, r}^x \|^2_{L^2 (\Omega)} \d r
  \right)^{\frac{1}{2}} .
\end{align*}
This implies
\begin{align*}
  \underset{u \in [s + \delta, t]}{\sup}  \| X_{s + \delta, u}^{x + h} - X_{s,
  u}^x \|_{L^2 (\Omega)} \leq & h + \delta \| b \|_{L^\infty(\R)} + \sqrt{\delta} \sqrt{2 (\alpha+1)^{-1}} \| \sigma
  \|_{L^\infty(\R)} \\ & + \int_{s + \delta}^t L_b  \underset{u \in [s + \delta,
  r]}{\sup} \| X_{s + \delta, u}^{x + h} - X_{s, u}^x \|_{L^2 (\Omega)} \d r \\ &  + \left( \int_{s + \delta}^t \frac{2}{\alpha + 1} L_{\sigma}^2
  \underset{u \in [s + \delta, r]}{\sup} \| X_{s + \delta, u}^{x + h} - X_{s,
  u}^x \|^2_{L^2 (\Omega)} \d r \right)^{\frac{1}{2}},
\end{align*}
which allows to use the extended Gronwall's Lemma, e.g Lemma A.1 in \cite{bossymaurer2024} with $p = 2$, $f (t) = \underset{u
\in [s + \delta, t]}{\sup}  \| X_{s + \delta, u}^{x + h} - X_{s, u}^x \|_{L^2
(\Omega)}$, $g (t) = L_b$, $h (t) = \frac{2}{\alpha + 1} L_{\sigma}^2$, $k (t)
= 0$ and $\ell (t) = h + \delta \| b \|_{L^\infty(\R)} +  \sqrt{\delta} \sqrt{2 (\alpha+1)^{-1}}\| \sigma \|_{L^\infty(\R)}$,
where $f$ and $\ell$ are positive and non-decreasing and $g, h, k \in L^1 ([s
+ \delta, t])$.
We thus obtain the inequality
\begin{equation*}
  \underset{u \in [s + \delta, t]}{\sup}  \| X_{s + \delta, u}^{x + h} - X_{s,
  u}^x \|_{L^2 (\Omega)}  \leq  4 ( h + \delta \| b \|_{L^\infty(\R)} +  \sqrt{\delta} \sqrt{2 (\alpha+1)^{-1}}\| \sigma \|_{L^\infty(\R)}) e^{4 L_b + \frac{8}{\alpha + 1} L_{\sigma}^2} .
\end{equation*}
Hence, setting $C \defeq \max \left( 4 e^{4 L_b + \frac{8}{\alpha + 1} L_{\sigma}^2}, 4
(\| b \|_{L^\infty(\R)} + \sqrt{2 (\alpha+1)^{-1}} \| \sigma \|_{L^\infty(\R)}) e^{4 L_b + \frac{8}{\alpha + 1}
L_{\sigma}^2} \right)$ and assuming that $\delta < 1$, we get
\begin{equation*}
  \| X_{s + \delta, u}^{x + h} - X_{s, u}^x \|_{L^2 (\Omega)} \leq 
  \underset{u \in [s + \delta, t]}{\sup}  \| X_{s + \delta, u}^{x + h} - X_{s,
  u}^x \|_{L^2 (\Omega)}
   \leq \text{ } C (h + \sqrt{\delta}),
\end{equation*}
leading to
\begin{equation*}
  \underset{(\delta, h) \rightarrow (0, 0)}{\lim} \| X_{s + \delta, u}^{x + h}
  - X_{s, u}^x \|_{L^2 (\Omega)} =  0.
\end{equation*}

We now turn to the derivative of the flow $\partial_x X$. Due to \cite[Proposition 4.1]{BretonPrivault}, $\partial_x X$ is solution of the following SDE:
\begin{equation*}
    \partial_x X_{s,t}^x = 1 + \int_s^t \partial_x X_{s,r}^x \ b'(X_{s,r}^x) \d r + \int_s^t \partial_x X_{s,r^-}^x \ \sigma'(X_{s,r^-}^x) \d L_r.
\end{equation*}
The proof is then similar to the one of $X$, using the pivot term
\begin{equation*}
    \partial_{x} X_{s+\delta,r}^{x+h}  \ \varphi(X_{s+\delta,r}^{x+h}) -\partial_x X_{s,r}^x  \ \varphi(X_{s,r}^x) =  \big(\partial_x X_{s+\delta,r}^{x+h} - \partial_x X_{s,r}^x \big) \  \varphi(X_{s+\delta,r}^{x+h}) + \partial_x X_{s,r}^x \big( \varphi(X_{s+\delta,r}^{x+h}) - \varphi(X_{s,r}^x) \big) 
\end{equation*}
for $\phi \in \{b', \sigma'\}$, the fact that
\begin{equation*}
    \big\| \big(\partial_x X_{s+\delta,r}^{x+h} - \partial_x X_{s,r}^x \big) \  \varphi(X_{s+\delta,r}^{x+h}) \big\|_{L^2(\Omega)} \leq \| \varphi \|_{L^\infty(\R)} \ \| \partial_x X_{s+\delta,r}^{x+h} - \partial_x X_{s,r}^x \|_{L^2(\Omega)}
\end{equation*}
and, due to the boundedness of $\phi$ and the Cauchy--Schwarz inequality,
\begin{align*}
    \big\| \partial_x X_{s,r}^x \big( \varphi(X_{s+\delta,r}^{x+h}) - \varphi(X_{s,r}^x) \big) \|_{L^2(\Omega)} \leq  & \sqrt{2} \| \phi \|_{L^\infty(\R)}^{\frac{1}{2}} \ \E[\big(\partial_x X_{s,r}^x\big)^2 \big( \varphi(X_{s+\delta,r}^{x+h}) - \varphi(X_{s,r}^x) \big)]^{\frac{1}{2}}  \\
    & \leq \sqrt{2} \| \varphi \|_{L^\infty(\R)}^{\frac{1}{2}} \  \|\partial_x X_{s,r}^x\|_{L^4(\Omega)} \ \| \varphi(X_{s+\delta,r}^{x+h}) - \varphi(X_{s,r}^x)
    \|_{L^2(\Omega)}^{\frac{1}{2}} \\
    & \leq \sqrt{2} \| \varphi \|_{L^\infty(\R)}^{\frac{1}{2}} \ L_{\varphi}^{\frac{1}{2}} \  \|\partial_x X_{s,r}^x\|_{L^4(\Omega)} \ \| X_{s+\delta,r}^{x+h} - X_{s,r}^x
    \|_{L^2(\Omega)}^{\frac{1}{2}}  \\
    & \leq \sqrt{2} \| \varphi \|_{L^\infty(\R)}^{\frac{1}{2}} \ L_{\varphi}^{\frac{1}{2}} \  \|\partial_x X_{s,r}^x\|_{L^4(\Omega)} \ C^\frac{1}{2} (\sqrt{h} + \delta^{\frac
    {1}{4}}),
\end{align*}
where $L_{\varphi}$ designates the Lipschitz constant of $\varphi$.
Moreover due to \cite[Theorem 5.1]{BretonPrivault} we have that 
\begin{equation*}
    \sup_{x \in \R} \sup_{0 \leqslant s \leqslant r \leqslant T} \| \partial_x X_{s,r}^x \|_{L^4(\Omega)} < + \infty.
\end{equation*}
In the end, for $\varphi \in \{ b', \sigma' \}$, there exists a constant $C'_\varphi>0$ such that
\begin{equation*}
    \big\| \partial_{x} X_{s+\delta,r}^{x+h}  \ \varphi(X_{s+\delta,r}^{x+h}) -\partial_x X_{s,r}^x  \ \varphi(X_{s,r}^x) \big\|_{L^2(\Omega)} \leq C_{\varphi}' \ \big(  \|\partial_x X_{s+\delta,r}^{x+h} - \partial_x X_{s,r}^x \|_{L^2(\Omega)} + \sqrt{h} + \delta^{1/4} \big),
\end{equation*}
while by Cauchy--Schwarz inequality along with an other application of \cite[Theorem 5.1]{BretonPrivault} and the boundedness of the coefficients, for $\varphi \in \{ b', \sigma' \}$ there exists a constant $C''_\varphi>0$ such that
\begin{equation*}
    \sup_{x \in \R} \sup_{0 \leqslant s \leqslant t \leqslant T} \| \partial_x X_{s,r}^x \  \phi(X_{s,r}^x) \|_{L^2(\Omega)} \leq C''_{\varphi}.
\end{equation*}
Then similarly to the proof of the continuity of the flow $X$ just above, we conclude by using Minkowski's integral inequality, Itô isometry and the extended Gronwal's Lemma \cite[Lemma A.1]{bossymaurer2024} applied to the function $f(t) \defeq \sup_{u \in [s+\delta,t]} \| \partial_x X_{s+\delta,u}^{x+h} - \partial_x X_{s,u}^{x} \|_{L^2(\Omega)}$. 

\section{Acknowledgements}
This work benefited from stimulating discussions with our respective PhD supervisors, Mireille Bossy and Ciprian Tudor, who are warmly
acknowledged.
The authors acknowledge the support from the Société de Mathématiques Appliquées et Industrielles, as well as the Centre International de Rencontres Mathématiques which allow us to work on this project within the framework on the "Projet BOUM de la SMAI" during one week in May 2025.
Paul Maurer acknowledges the support of the French Agence Nationale de la Recherche
(ANR), under grant ANR-21-CE30-0040-01 (NETFLEX). Jérémy Zurcher acknowledges the support from Labex CEMPI (ANR-11-LABX-007-01).

\appendix
\section{Appendix} \label{sec:appendix}

\subsection{Some useful results from measure theory} \label{sec:MeasureTheory}

In the proof of Theorem \ref{maintheo}, we work with integrals against $\sigma$-finite Lévy measure $\nu$. For this reason, we need the general version of Vitali's Theorem to conclude about convergence of these integrals in Section \ref{subsec:cvg_h_zero}. We fix a $\sigma$-finite measure space $(X,\mathcal{A},\mu)$, and introduce the following definition.

\begin{Definition}
    We say that a family $\mathcal{F}$ of real valued measurable functions on $X$ has \emph{uniformly absolutely continuous integrals} if 
    \begin{enumerate}
        \item $\forall f \in \mathcal{F}, \quad f \in L^1(\mu)$,
        \item $\forall \eps > 0, \ \exists \delta > 0, \ \forall A \in \mathcal{A}, \quad \big( \mu(A) \leqslant \delta \implies \forall f \in \mathcal{F}, \ \int_A |f| \ \d \mu \leqslant \eps \big).$
    \end{enumerate}
    We also recall that a family $\mathcal{F}$ of real valued measurable functions on $X$ is said to be \emph{tight} if 
\begin{equation*}
    \forall \eps > 0, \ \exists K_{\varepsilon} \in \mathcal{A}, \ \mu(K_{\varepsilon}) < + \infty  \text{ and }  \forall f \in \mathcal{F}, \ \int_{X \backslash K_{\varepsilon}} |f| \ \d  \mu \leq \eps. 
\end{equation*}
\end{Definition}

The following theorem, that we will call \emph{Vitali's $L^1$-convergence Theorem}, is a consequence from \cite[Theorem 16.6]{Schilling_2005} and \cite[Theorem 16.8]{Schilling_2005}. 

\begin{Theoreme} \label{thm:vitali_general}
    Let $(f_n)_{n \in \N}$ be a sequence of real valued measurable functions on $X$. Assume that $(f_n)$ converges in measure to a measurable function $f$. Then $(f_n)_{n \in \N}$ converges in $L^1 (\mu)$ toward $f$ if and only if $(f_n)_{n \in \N}$  has uniformly absolutely continuous integrals and is tight. 
\end{Theoreme}

We will also prove the following useful criterium. 

\begin{Lemma} \label{lm:critere_vallee_poussin}
Let $\mathcal{F} \subset L^1(\mu)$ be a family of real valued integrable functions on $X$. Let $G: [0, \infty) \longrightarrow [0, \infty)$ increasing positive map such that $\frac{G(t)}{t} \to \infty$ when $t \to \infty$. Suppose moreover that $\sup_{f \in \Ff} \int_{X} G \circ |f| \ \d \mu < + \infty$. Then $\mathcal{F}$ has uniformly absolutely continuous integrals.
\end{Lemma}

\begin{Demo}
    Let $\eps >0$. Since $\lim_{t \to + \infty} G(t)/t = +\infty$, in particular $\underline{\lim}_{t \to + \infty} G(t)/t = +\infty$. It follows, by definition, that $\lim_{c \to + \infty} h(c) = +\infty$ where $h(c) \defeq \inf_{t > c} G(t) / t$. 
    Let $c_{\eps} > 0$ big enough so that $h(c_{\eps}) > \frac{2}{\eps} \sup_{f \in \Ff} \int_{X} G \circ |f| \d \mu$. We now set $\delta \defeq \frac{\eps}{2 c_{\eps}} > 0$. Let $A \in \mathcal{A}$ such that $\mu(A) < \delta$. Then, by writing
    \begin{equation*}
        \int_A |f| \d \mu = \int_A |f| \mathbf{1}_{\{|f| \leq c_\eps \}} \ \d \mu + \int_A |f|  \mathbf{1}_{\{|f| > c_\eps \}} \ \d \mu,
    \end{equation*}
    we have
    \begin{equation*}
         \int_A |f| \mathbf{1}_{\{|f| \leq c_\eps \}} \ \d \mu \leqslant c_\eps \mu(A) < \eps/2.
    \end{equation*}
    Moreover,
    \begin{align*}
        \int_A |f|  \mathbf{1}_{\{|f| > c_\eps \}} \ \d \mu & =  \int_A \big( G \circ |f| \big) \times \frac{|f|}{G \circ |f|} \mathbf{1}_{\{|f| > c_\eps \}} \ \d \mu \\
        & \leqslant \sup_{t > c_\eps} \frac{t}{G(t)} \times \int_A G \circ |f| \ \d \mu \\ 
        & \leqslant \frac{1}{h(c_\eps)} \sup_{f \in \Ff} \int_{X} G \circ |f| \ \d \mu < \eps/2. 
    \end{align*}
    It follows that $\int_A |f| \ \d \mu < \eps$, which proves the second point of the definition of the absolute uniform continuity of the integrals.  
\end{Demo}

Finally, for the convergence of the Skorohod integral, we will make use of the following proposition which is a slight adaptation of \cite[Theorem 4]{jakszto2010another}.

\begin{Proposition} \label{prop:measure_CVG_implies_weakL2_CVG}
    Let $(X,\mathcal{A},\mu)$ be a complete $\sigma$-finite measure space. Let $(f_n)_{n \in \N}$ be a sequence of real-valued measurable functions on $X$ such that $\sup_{n \in \N} \|f_n\|_{L^2(\mu)} < + \infty$. Under these assumptions, if $(f_n)_{n \in \N}$ converges in measure to a measurable function $f$, then it converges to $f$ under the $L^2(\mu)$-weak topology. In other words, for any $w \in L^2(\mu)$, we have
    \begin{equation*}
        \int_X f_n w \ \d \mu \xrightarrow[n \to +\infty]{} \int_X f w \ \d \mu.
    \end{equation*}
\end{Proposition}

\begin{Demo}
    The proof is similar to \cite[Theorem 4]{jakszto2010another}, itself adapted from \cite[Theorem 1]{jakszto2010another}. The difference is that we work in a general measure space $(X,\mathcal{A},\mu)$ instead of $(\R^d,\mathcal{B}(\R^d),\text{Leb}_{\R^d})$. The adaptation is straightforward for the following reasons:
    \begin{enumerate}
        \item Any Hilbert space is weakly sequentially compact (see e.g. \cite[Chapter 5, Lemma 1.4]{kato2013perturbation}). In particular, from any bounded sequence of elements of the Hilbert space $L^2(\mu)$, we can extract a subsequence that converges for the weak topology of $L^2(\mu)$. 
        \item The space $L^2(\mu)$ is uniformly convex, hence it possesses the weak Banach--Saks property.
        \item If a sequence of elements of $L^2(\mu)$ is convergent for the strong topology of $L^2(\mu)$, then we can extract a subsequence that converges $\P$-a.s to the same limit. This follows from \cite[Lemma 5.21]{abramovich2002invitation}, and from the fact that $L^2(\mu)$ is a Banach function space on $(X,\mathcal{A},\mu)$ (see e.g \cite[Section 3.1]{jakszto2012some}).
    \end{enumerate}
\end{Demo}

\begin{Remarque}
    Contrary to what is suggested in \cite{ItoAG} after equation (43), it is not necessary for the space $L^2(\mu)$ to be separable. In fact, if $(\Omega,\Ff,\P)$ is a probability space, it is not always true that $L^2(\Omega,\Ff,\P)$ is separable (this is only true if the $\sigma$-algebra $\Ff$ is assumed to be countably generated). But as mentioned, this has no impact on the final result.
\end{Remarque}

\begin{Remarque}
    Consider the measurable space $(\R,\mathcal{B}(\R))$. Recall that a Lévy measure $\nu$ on $(\R,\mathcal{B}(\R))$ is a Borel measure such that $\nu(\{0\}) = 0$ and $\int_{\R} (1 \wedge z^2) \nu(\d z) < + \infty$.  The measure $\nu$ is $\sigma$-finite, since $\R = \{0\} \cup \bigcup_{n \geqslant 1}^{\uparrow} \{ z \in \R \ ; \ |z| > 1/n \}$. Therefore, the result of Proposition \ref{prop:measure_CVG_implies_weakL2_CVG} applies for the measure space $(\R,\mathcal{B}(\R),\nu)$.
\end{Remarque}

\subsection{Itô formula for independent random fields with càdlàg semimartingales}

\begin{Proposition} \label{prop:ito_formula_indep_cadlag}
Let $(\Omega, \mathcal{F}, (\mathcal{F}_t)_{t \geq 0}, \mathbb{P})$ be a filtered
probability space. Let $T > 0$ be a fixed time horizon and $\mathcal{G}_T
\subset \mathcal{F}$ a $\sigma$-algebra independent of $\mathcal{F}_t$ for
every $t \in [0, T]$. \

Let $f: \Omega \times \mathbb{R} \rightarrow \mathbb{R}$ such that for
almost all $\omega \in \Omega$, $f (\omega, \cdot) \in \mathcal{C}^2
(\mathbb{R} \rightarrow \mathbb{R})$ and such that for every $x \in
\mathbb{R}$, $f (\cdot, x)$ is $\mathcal{G}_T$-measurable, and $(X_t)_{t \in
[0, T]}$ a $(\mathcal{F}_t)_{t \in [0, T]}$-càdlàg semimartingale\footnote{Here by "càdlàg semimartingale" we mean an adapted process $X$ with right-continuous and left-limited paths, such that the decomposition $X_t = X_0 + M_t + A_t$ holds, where $M_0 = A_0 = 0$, $M$ is a locally square integrable martingale, and $A$ is càdlàg, adapted, with paths of finite variation on compacts.}. Let
$\mathcal{H}_t =\mathcal{F}_t \vee \mathcal{G}_T = \sigma(\mathcal{F} , \mathcal{G})$. We consider the stochastic process $[0,T] \times \Omega \ni (t,\omega) \mapsto f(\omega,X_t(\omega))$, denoted by $(f(X_t))_{t \in [0,T]}$. 

Then $(X_t)_{t \in [0, T]}$ and $ f (X_t)_{t \in [0, T]}$ both are $(\mathcal{H}_t)_{t \in [0,
T]}$-semimartingale and we have, $\mathbb{P}$-p.s, for any $t \in [0, T],$
\begin{multline} \label{eq:ito_formula_indep_cadlag}
    f (X_t) - f (X_0) = \int_0^t f' (X_s) \ \d X_s ^{\mathcal{H}} + \frac{1}{2}  \int_0^t f''
   (X_s) \d [X]_s \\ + \sum_{0 < s \leq t} \{ f(X_s) - f(X_{s^-}) - f'(X_{s^-}) \Delta X_s - \frac{1}{2} f''(X_{s^-}) (\Delta X_s)^2 \},
\end{multline}
where $[X]$ designates the quadratic covariation of $X$, and the notation $\d X_s ^{\mathcal{H}}$ indicates that the Itô integration is taken with respect to the filtration $\mathcal{H}$. 
\end{Proposition}

\begin{Demo}
    The proof closely follows that of the standard Itô formula for càdlàg semimartingales (see \cite[Theorem 32]{protter}), but we present it here in full for the sake of completeness. Firstly, we check that the stochastic process $X$ is in fact a $\Hh$-semimartingale. Secondly, we prove that \eqref{eq:ito_formula_indep_cadlag} holds using a time-discretisation and Taylor's formula.
    \paragraph*{$X$ is a $\Hh$-semimartingale.}
    The process $(X_t)_{t \in [0,T]}$ is adapted to $(\Ff_t)_{t \in [0,T]}$ and we have $\Ff_t \subset \Hh_t$ for any $t \in [0,T]$, so $(X_t)_{t \in [0,T]}$ is adapted to $(\Hh_t)_{t \in [0,T]}$. Let $M,A$ be two stochastic processes on $[0,T]$ such that $M_0 = A_0 = 0$, $M$ is a  $(\Ff_t)_{t \in [0,T]}$-locally square integrable martingale, $A$ is a $(\Ff_t)_{t \in [0,T]}$-adapted càdlàg process with paths of finite variation on compacts, and $X_t = X_0 + A_t + M_t$ for every $t \in [0,T]$. Firstly, $A$ and $M$ are clearly adapted to $(\Hh_t)_{t \in [0,T]}$. Secondly, if $(\tau_n)_{n \in \N}$ is a localizing sequence for $M$, so that $M^{\tau_n}$ is a square integrable $(\Ff_t)_{t \in [0,T]}$-martingale for every $n \in \N$, then for $s \leq t \in [0,T]$ we have that
    \begin{equation}
        \E[M^{\tau_n}_t  |  \Hh_s ] = \E[M^{\tau_n}_t  |  \sigma(\Ff_s, \Gg_T) ] = \E[M^{\tau_n}_t  | \Ff_s ]
    \end{equation} 
    since $\Gg_T$ is independent from $\sigma(\sigma(M_t^{\tau_n}),\Ff_s)$. Since $M^{\tau_n}$ is a $(\Ff_t)_{t \in [0,T]}$-martingale it follows that $ \E[M^{\tau_n}_t  |  \Hh_s ] = M^{\tau_n}_s$, implying that $M^{\tau_n}$ is a $(\Hh_t)_{t \in [0,T]}$-(square integrable) martingale, so $M$ is a $(\Hh_t)_{t \in [0,T]}$-locally square integrable martingale.
    This finally shows that $(X_t)_{t \in [0,T]}$ is a $(\Hh_t)_{t \in [0,T]}$-semimartingale.

    \paragraph*{Application of Taylor's formula.}
    The proof of \eqref{eq:ito_formula_indep_cadlag} relies on the following Taylor's formula.

    \begin{Lemma} \label{lm:taylor_formula}
        Let $K \subset \R$ be a compact set and $\varphi:K \rightarrow \R$ be a function of class $\Cc^2$ on $K$. There exists an increasing function $r: \R_+ \rightarrow \R_+$ with $\lim_{u \to 0} r(u) = 0$ and a function $R: K^2 \rightarrow \R$ satisfying $|R(x,y)| \leq r(|x-y|)(y-x)^2$ for every $(x,y) \in K^2$ such that 
        \begin{equation*}
            \varphi(y) - \varphi(x) = \varphi'(x)(y-x) + \frac{1}{2} \varphi''(x) (y-x)^2 + R(x,y), \quad (x,y) \in K^2. 
        \end{equation*}
    \end{Lemma}
    Since $X$ has càdlàg paths, it is almost surely bounded on the compact $[0,T]$.  Therefore, for $\P$-almost every $\omega \in \Omega$ we may assume that there exists $k(\omega) \in \R_+$ such that the function $t \ni [0,T] \mapsto X_t(\omega)$ takes its values in $[-k(\omega),k(\omega)]$. This will allows us to apply \ref{lm:taylor_formula} on $x \mapsto f(\omega,x)$.

    We fix $t \in [0,T]$ and consider a refining partition $0 = t_0^n \leq \cdots \leq t_{\kappa_n}^n = t$ of $[0,t]$, such that $\lim_{n \to +\infty} \sup_{0 \leq i \leq \kappa_n-1} | t_{i+1}^n - t_i^n | = 0$. 
    Let $J \defeq \{ s \in [0,t] \: \ \Delta X_s \neq 0 \}$ be the set of the jump times of $X$ on $[0,t]$. Since $X$ is a $(\Hh_t)_{t \in [0,T]}$-semimartingale, the series $\sum_{0 < s \leq t} (\Delta X_s)^2 = \sum_{s \in J} (\Delta X_s)^2$ converges, hence for any $\eps > 0$ we may find a partition $J = A_{\eps} \sqcup B_{\eps}$ such that $\sum_{s \in B_\eps} (\Delta X_s)^2 \leq \eps^2$ and $A_\eps$ contain a finite number of elements. We now define the sets
    \begin{align*}
        L^\eps_n \defeq & \{i \in \{0, \ldots, \kappa_n -1 \} \: \ (t_i^n, t_{i+1}^n] \cap A_\eps \neq \varnothing   \} \\
        S^\eps_n \defeq & \{i \in \{0, \ldots, \kappa_n -1 \} \setminus L^\eps_n \: \ (t_i^n, t_{i+1}^n] \cap B_\eps \neq \varnothing   \} \\
        C^\eps_n \defeq & \{0, \ldots, \kappa_n -1 \}  \setminus \big( S^\eps_n \cup L^\eps_n \big), \\
    \end{align*}
    so that $L^\eps_n$ contains the intervals with large jumps of $X$, $S^\eps_n$, contains the intervals with small jumps of $X$, and $C^\eps_n$ contains the intervals where $X$ is continuous. Note that $L^\eps_n \sqcup S^\eps_n \sqcup C^\eps_n$ forms a partition of the indices $\{0, \ldots, \kappa_n -1 \}$, allowing to decompose $f(X_t) - f(X_0)$ into the three following sums:
    \begin{align*}
        f(X_t) - f(X_0) = & \sum_{i = 0}^{\kappa_n-1} \big(f(X_{t_{i+1}^n}) - f(X_{t_i^n}) \big) \\
        = & \sum_{i \in L^\eps_n} \big(f(X_{t_{i+1}^n}) - f(X_{t_i^n}) \big) + \sum_{i \in S^\eps_n} \big(f(X_{t_{i+1}^n}) - f(X_{t_i^n}) \big) + \sum_{i \in C^\eps_n} \big(f(X_{t_{i+1}^n}) - f(X_{t_i^n}) \big).
    \end{align*} 
    Applying Lemma \ref{lm:taylor_formula} with $\varphi = f(\omega, \cdot)$ for $\P$-almost every $\omega$,  we obtain that the existence of an increasing stochastic process $r: \Omega \times \R_+ \rightarrow \R_+$ with $\lim_{u \to 0} r(u) = 0$ $\P$-a.s and a stochastic process $R: \Omega \times \R^2$ satisfying $|R(x,y)| \leq r_{\omega}(|x-y|) (y-x)^2$ for every $(x,y) \in [-k,k]$  $\P$-a.s, such that we have
    \begin{align*}
        f(X_t) - f(X_0) = & \sum_{i \in L^\eps_n} \big(f(X_{t_{i+1}^n}) - f(X_{t_i^n}) \big) \\
        & + \sum_{i = 0}^{\kappa_n-1} \big( f'(X_{t_i^n}) \big(X_{t_{i+1}^n} - X_{t_i^n} \big) + \frac{1}{2} f''(X_{t_i^n}) \big(X_{t_{i+1}^n} - X_{t_i^n} \big)^2 \big) \\
        & \quad - \sum_{i \in L^\eps_n} \big( f'(X_{t_i^n}) \big(X_{t_{i+1}^n} - X_{t_i^n} \big) + \frac{1}{2} f''(X_{t_i^n}) \big(X_{t_{i+1}^n} - X_{t_i^n} \big)^2 \big) \\
        & \qquad + \sum_{i \in S^\eps_n \sqcup C^\eps_n} R(X_{t_i^n},X_{t_{i+1}^n}) \quad \P\text{-a.s}.
    \end{align*}
    Since $X$ is a $(\Hh_t)_{t \in [0,T]}$-semimartingale and the maps $(\omega,t) \mapsto f'(\omega,X_t(\omega))$ and $(\omega,t) \mapsto f''(\omega,X_t(\omega))$ are càdlàg $(\mathcal{H})$-adapted stochastic processes, we may apply respectively \cite[Theorem 21]{protter} and \cite[Theorem 30]{protter} to show that $\sum_{i = 0}^{\kappa_n-1} f'(X_{t_i^n}) (X_{t_{i+1}^n} - X_{t_i^n} ) $ converges in probability to $\int_0^t f'(X_{s^-}) \d X_s ^{\mathcal{H}}$, and that $\sum_{i = 0}^{\kappa_n-1} \frac{1}{2} f''(X_{t_i^n}) (X_{t_{i+1}^n} - X_{t_i^n} )^2$ converges in probability to $\frac{1}{2} \int_0^t f''(X_{s^-}) \d [X]_s$ respectively. We will demonstrate that $\P$-a.s the following convergences hold:
    \begin{equation} \label{eq:first_cvg_jumpsito}
        \sum_{i \in L^\eps_n} \big(f(X_{t_{i+1}^n}) - f(X_{t_i^n}) \big) \xrightarrow[n \to +\infty]{} \sum_{s \in A_\eps} \big(f(X_{s}) - f(X_{s^-}) \big)
    \end{equation}
    and
    \begin{multline} \label{eq:second_cvg_jumpsito}
        \sum_{i \in L^\eps_n} \big( f'(X_{t_i^n}) \big(X_{t_{i+1}^n} - X_{t_i^n} \big) + \frac{1}{2} f''(X_{t_i^n}) \big(X_{t_{i+1}^n} - X_{t_i^n} \big)^2 \big) \\ \xrightarrow[n \to +\infty]{} \sum_{s \in A_\eps} \big( f'(X_{s^-}) \Delta X_s + \frac{1}{2} f''(X_{s^-}) (\Delta X_s)^2 \big).
    \end{multline}
    More precisely, we will show in detail that \eqref{eq:first_cvg_jumpsito} holds, and \eqref{eq:second_cvg_jumpsito} will follow by similar arguments.
    \paragraph*{Proof of \eqref{eq:first_cvg_jumpsito}.}
    Observe that if $A_\eps$ is empty, then $\sum_{i \in L_n^\eps} (f (X_{t_{i + 1}^n}) - f (X_{t_i^n})) = 0$ and there is nothing to do. We now assume that $A_\eps$ is not empty. We let $m \in
     \mathbb{N} \backslash\{0 \}$ and $s_1 < \cdots < s_m$ be elements of $[0, t]$ such that
    $A_\eps = \{ s_1, \ldots, s_m \}$ and we fix $\eta > 0$. Since $f$ is a.s continuous
    and $X$ is a.s right-continuous with left limits, there exists $\delta > 0$ such that:
    \begin{equation*}
        \forall j \in \{ 1, \ldots, m \}, \text{\quad} \forall s \in [0, t],
    \text{\quad} s_j < s \leq s_j + \delta \text{ } \Longrightarrow \text{ }
    | f (X_s) - f (X_{s_j}) | \leq \eta
    \end{equation*}
    and
    \begin{equation*}
        \forall j \in \{ 1, \ldots, m \}, \text{\quad} \forall s \in [0, t],
    \text{\quad} s_j - \delta \leq s < s_j \text{ } \Longrightarrow \text{ }
    | f (X_s) - f (X_{s_j^-}) | \leq \eta.
    \end{equation*}
    We set $\delta' = \min \left( \delta, \underset{0 \leq j \leq m - 1}{\inf} | s_{j
    + 1} - s_j | \right)$. Using that $\underset{n \rightarrow \infty}{\lim} 
    \underset{1 \leq i \leq \kappa_n - 1}{\sup}  | t_{i + 1}^n - t_i^n | = 0$,
    there exists $N \in  \mathbb{N} \backslash\{0 \}$ such that for any $n \geq N$, we have
    that
    \begin{equation*}
    \forall i \in \{ 1, \ldots, \kappa_n - 1 \}, \text{\quad} | t_{i + 1}^n -
    t_i^n | \leq \frac{\delta'}{2}.
    \end{equation*}
    Let $n \geq N$. In view of the definition of $\eta$, we have that for any $j
    \in \{ 1, \ldots, m \}$, there exists a unique $\tau (j) \in \{ 1, \ldots,
    \kappa_n \}$ such that $s_j \in (t_{\tau (j)}^n, t_{\tau (j) + 1}^n]$, and we
    can rewrite the sum as
    \begin{equation*}
    \sum_{i \in L_n^\eps} (f (X_{t_{i + 1}^n}) - f (X_{t_i^n})) = \sum_{j = 1}^m (f
    (X_{t_{\tau (j) + 1}^n}) - f (X_{t_{\tau (j)}^n} )).
    \end{equation*}
    Moreover, the fact that $s_j \in (t_{\tau (j)}^n, t_{\tau (j) + 1}^n]$ implies
    $s_j \leq t_{\tau (j) + 1}^n < s_j + \delta'$ and $s_j - \delta' < t_{\tau (j)}^n <
    s_j$, hence $\text{ } | f (X_{t_{\tau (j) + 1}^n}) - f (X_{s_j}) | \leq
    \eta$ and $| f (X_{t_{\tau (j)}^n}) - f (X_{s_j}) | \leq \eta$.
    It follows that
    \begin{equation*}
    \left| \sum_{j = 1}^m (f (X_{t_{\tau (j) + 1}^n}) - f (X_{t_{\tau (j)}^n}
    )) - (f (X_{s_j}) - f (X_{s_j^-})) \right| \leq 2 m \eta,
    \end{equation*}
    leading to
    \begin{equation*}
    \underset{n \rightarrow + \infty}{\lim}  \sum_{i \in L^\eps_n} (f (X_{t_{i + 1}^n})
    - f (X_{t_i^n})) = \sum_{s \in A_\eps} (f (X_s) - f (X_{s^-})).
    \end{equation*}
    \paragraph*{Control of the leftover term.} We will now estimate the term $\sum_{i \in S^\eps_n \sqcup C^\eps_n} R(X_{t_i^n},X_{t_{i+1}^n})$. By construction, we have that $\P$-a.s,
    \begin{align*}
        \sum_{i \in S^\eps_n \sqcup C^\eps_n} R(X_{t_i^n},X_{t_{i+1}^n}) \leq & \sum_{i \in S^\eps_n \sqcup C^\eps_n} r(|X_{t_{i+1}^n}-X_{t_i^n}|) (X_{t_{i+1}^n}-X_{t_i^n})^2 \\
        & \leq \sup_{i \in S^\eps_n \sqcup C^\eps_n} r(|X_{t_{i+1}^n} - X_{t_i^n}|) \times \sum_{i=0}^{\kappa_n -1} (X_{t_{i+1}^n}-X_{t_i^n})^2. 
    \end{align*}
    When considering a continuous semimartingale $X$, the term $\sup_{i \in S^\eps_n \sqcup C^\eps_n} r(|X_{t_{i+1}^n} - X_{t_i^n}|)$ is managed by the a.s. uniform continuity of the trajectories of $X$ on the interval $[0,t]$. However, in our case where $X$ is only assumed to be càdlàg, we have to rely on the analogue of uniform continuity for right-continuous processes. To that extend, we adopt the formalism from \cite[Section 12]{billingsley2013convergence}. We define, for $A \subset [0,t]$, 
    \begin{equation*}
        w(A) \defeq \sup_{s,r \in A} |X_r - X_s|.
    \end{equation*}
    Let $\delta > 0$. We call $\delta$-sparse set any partition $0 = s_0 < \ldots < s_v = t$ of $[0,t]$ with mesh greater than $\delta$, i.e. such that  
    \begin{equation*}
        \min_{1 \leq i \leq v} (s_i - s_{i-1}) > \delta.
    \end{equation*}
    We denote by $\mathcal{S}(\delta)$ the set containing all the $\delta$-sparse sets. 
    Finally, we define 
    \begin{equation*}
        w'(\delta) \coloneqq \inf_{\{s_i\}_{0 \leq i \leq \nu} \in \mathcal{S}(\delta) } \max_{1 \leq i \leq v} w([s_{i-1},s_i)).
    \end{equation*}
    From \cite{billingsley2013convergence} (beginning of page 123 and consequence of Lemma 12.1 therein), we know that $\lim_{\delta \to 0} w'(\delta) = 0$. This result can be considered as the analogue of Heine's uniform continuity Theorem for càdlàg functions. In particular, there exists a $\delta > 0$ such that $ w'(\delta) < \eps$. Since we have $\lim_{n \to +\infty} \sup_{0 \leq i \leq \kappa_n} |t_{i+1}^n - t_i^n| = 0$, we can choose an integer $N \in \mathbb{N} \backslash\{0 \}$ such that for every $n \geq N$, we have $|t_{i+1}^n - t_i^n| < \delta$ for every $i = 0, \ldots, \kappa_n -1$. We fix such a $n \geq N$.
    By definition of the infimum, we can choose a $\delta$-sparse set $\{ s_j \}_{0 \leq j \leq \nu}$ such that $w([s_{j-1},s_j)) < w'(\delta) + \eps$ for every $j \in \{0, \ldots, \nu \}$. We consider afterwards an integer $i \in S^\eps_n \sqcup C^\eps_n$. Since  $|t_{i+1} - t_{i}| \leq \delta$, two cases may appear.
    \begin{itemize}
        \item Either $t_{i+1}$ and $t_i$ lie in the same interval $[s_{j-1},s_j)$ and $|X_{t_{i+1}} - X_{t_i}|$ is bounded above by $w([s_{j-1},s_j)) < w'(\delta) + \eps$,
        \item either $t_{i+1}$ and $t_i$ lie in adjacent intervals $[s_{j-1},s_j]$ and $[s_j,s_{j+1}]$ and $|X_{t_{i+1}} - X_{t_i}|$ is bounded above by $|X_{ {s_j}^-} - X_{t_i}| + |X_{s_j} - X_{{s_j}^-}| + |X_{t_{i+1}} - X_{s_j}| \leq  2 w'(\delta) + 2 \eps + \eps$ since $|\Delta X_{r}| \leq \eps$ for any $r \in (t_i,t_{i+1}]$.
    \end{itemize}
    In both cases we end up with the upper-bound
    \begin{equation*}
        |X_{t_{i+1}} - X_{t_i}| \leq 5 \eps.
    \end{equation*}
    Since $r$ is  $\P$-a.s increasing, it follows that
    \begin{equation*}
        r(|X_{t_{i+1}} - X_{t_i}|) \leq r(5 \eps),
    \end{equation*}
    and this is true for any $i \in S^\eps_n \sqcup C^\eps_n$ with $n \geq N$, so we deduce that $\P$-a.s, we have:
    \begin{equation*}
        \sup_{i \in S^\eps_n \sqcup C^\eps_n} r(|X_{t_{i+1}^n} - X_{t_i^n}|) \leq r(5 \eps),
    \end{equation*}
    and it follows that 
    \begin{equation*}
       \sum_{i \in S^\eps_n \sqcup C^\eps_n} R(X_{t_i^n},X_{t_{i+1}^n}) \leq r(5 \eps) \sum_{i=0}^{\kappa_n -1} (X_{t_{i+1}^n}-X_{t_i^n})^2.
    \end{equation*}
    Then we have that $\sum_{i=0}^{\kappa_n -1} (X_{t_{i+1}^n}-X_{t_i^n})^2$ converges in probability to $[X]_t$. Consequently, in probability: 
    \begin{equation} \label{eq:limsup_bound_for_leftover}
        \limsup_{n \to + \infty} \sum_{i \in S^\eps_n \sqcup C^\eps_n} R(X_{t_i^n},X_{t_{i+1}^n}) \leq r(5 \eps) [X]_t. 
    \end{equation}
    \paragraph{Conclusion.}
    Letting $\eps \to 0$ in \eqref{eq:limsup_bound_for_leftover}, since $r(\eps) \to 0$ we obtain that in probability,
    \begin{equation*}
         \limsup_{n \to + \infty} \sum_{i \in S^\eps_n \sqcup C^\eps_n} R(X_{t_i^n},X_{t_{i+1}^n}) \xrightarrow[\eps \to 0]{} 0.
    \end{equation*}
    On the other hand, since $\mathbf{1}_{A_\eps} \xrightarrow[\eps \to 0]{} \mathbf{1}_{J}$, the dominated convergence Theorem implies that
    \begin{multline*}
        \sum_{s \in A_\eps} \big(f(X_{s}) - f(X_{s^-}) -  f'(X_{s^-}) \Delta X_s - \frac{1}{2} f''(X_{s^-}) (\Delta X_s)^2  \big) \\
        \xrightarrow[\eps \to 0]{} \sum_{s \in J} \big(f(X_{s}) - f(X_{s^-}) -  f'(X_{s^-}) \Delta X_s - \frac{1}{2} f''(X_{s^-}) (\Delta X_s)^2  \big) 
    \end{multline*}
    providing the latter series is absolutely convergent, but this is the case since $X$ has value in $[-k,k]$ almost surely, $f(\omega, \cdot) \in C^2(\R \rightarrow \R)$ almost surely, and $\sum_{s \in J} (\Delta X_s)^2 =  [X]_t < + \infty$. This ends the proof of Proposition \ref{prop:ito_formula_indep_cadlag}.   
\end{Demo}

\printbibliography

\end{document}